\input amstex
\documentstyle{amsppt}
\nologo
\TagsOnRight
\magnification 1200
\voffset -1cm

\def\Z{\Bbb Z}  \def\Q{\Bbb Q} \def\R{\Bbb R} \def\N{\Bbb N}
\let\phi\varphi \let\x\times \let\tl\tilde
\def\on#1{\expandafter\def\csname#1\endcsname{{\operatorname{#1}}}}
\on{Hom} \on{Ext} \on{Int} \on{lf} \on{Tel} \let\but\setminus \on{lk} \on{Lk}
\let\eps\varepsilon \on{coker} \on{Tor} \on{c} \on{id} \on{Cl} \on{THom}
\def\dirlim{\lim\limits_\rightarrow{}} \def\imp{$\Rightarrow$} \on{rel}
\def\invlim{\lim\limits_\leftarrow{}} \def\iff{$\Leftrightarrow$} \on{st}
\def\derlim{\lim\limits_\leftarrow{^1}{}} \on{Fr} \on{ker} \on{im} \on{sgn}
\let\ssm\rightsquigarrow \let\emb\hookrightarrow \def\when{$\Leftarrow$}
\font\ssym=cmbsy5 at 2.5pt 
\font\sssym=cmbsy5 at 2pt
 
\def\striangle{\text{\ssym\char"34}}
\def\sstriangle{\text{\sssym\char"34}}
\def\sing#1#2{\smash{\overset\striangle\to{\smash{#1}
\vphantom{\vrule height#2pt}}}}
\def\ssing#1#2{\smash{\overset\sstriangle\to{\smash{#1}
\vphantom{\vrule height#2pt}}}}
\def\sH{\sing{H\,}{6}\!} \def\spi{\sing\pi{3.5}} \def\stau{\sing\tau{3.5}}
\def\sstau{\ssing\tau{2}} \def\sC{\sing{\text{C}}{6}} \def\sitC{\sing{C\,}{6}\!}
\def\ssC{\sing{\text{C}}{5}}
\let\q\backslash \def\tr{{\sssize\triangle}}

\topmatter 
\address Steklov Mathematical Institute, Division of Geometry and Topology;
Gubkina ~8, Moscow 119991, Russia \endaddress
\email melikhov\@mi.ras.ru \endemail

\thanks Supported by Russian Foundation for Basic Research Grant No.
08-01-00663-a; Mathematics Branch of the Russian Academy of Sciences Program
``Theoretical Problems of Contemporary Mathematics''; and President Grant
MK-5411.2007.1.
\endthanks

\title Steenrod homotopy \endtitle

\author Sergey A. Melikhov \endauthor

\abstract Steenrod homotopy theory is a natural framework for doing algebraic
topology on general spaces in terms of algebraic topology of polyhedra;
or from a different viewpoint, it studies the topology of the lim$^1$ functor
(for inverse sequences of groups).
This paper is primarily concerned with the case of compacta, in which Steenrod
homotopy coincides with strong shape.
An attempt is made to simplify the foundations of the theory and to clarify and
improve some of its major results.

We use geometric tools such as Milnor's telescope compactification, comanifolds
(=mock bundles) and the Pontryagin--Thom Construction, to obtain new simple
proofs for results by Barratt--Milnor; Cathey; Dydak; Dydak--Segal; Eda--Kawamura;
Edwards--Geoghegan; Fox; Geoghegan--Krasinkiewicz; Jussila; Krasinkiewicz--Minc;
Marde\v si\'c; Mittag-Leffler--Bourbaki; and for three unpublished results by
Shchepin.
An error in Lisitsa's proof of the ``Hurewicz theorem in Steenrod homotopy''
is corrected.
It is shown that over compacta, R. H. Fox's overlayings are equivalent to
I.~ M.~ James' uniform covering maps.
Other results include:

$\bullet$
A morphism between inverse sequences of countable (possibly non-abelian) groups that
induces isomorphisms on $\invlim$ and $\derlim$ is invertible in the pro-category.
This implies the ``Whitehead theorem in Steenrod homotopy'', thereby answering
two questions of A. Koyama.

$\bullet$
If $X$ is an LC$_{n-1}$ compactum, $n\ge 1$, its $n$-dimensional Steenrod
homotopy classes are representable by maps $S^n\to X$, provided that $X$
is simply connected.
The assumption of simple connectedness cannot be dropped, by a well-known
example of Dydak and Zdravkovska.

$\bullet$
A connected compactum is Steenrod connected (=pointed $1$-movable) iff every
its uniform covering space has countably many uniform connected components.
\endabstract

\endtopmatter

\vskip-15pt
\centerline{\bf Contents}
{\eightrm
1. Introduction

\ \ A. Motivations

\ \ B. Basic notions

\ \ C. Steenrod theories versus singular ones

\ \ D. Comparison with literature

\ \ E. Synopsis

2. Steenrod homotopy category

3. Homotopy groups

4. Homology and cohomology

5. Some examples

6. Comparison of theories

7. Covering theory

8. Zero-dimensional homotopy
}

\head 1. Introduction \endhead

Rephrasing a well-known quote of A. N. Whitehead, one can declare that
the ultimate goal of Steenrod homotopy theory is to eliminate any need for
general topology --- by reducing it to combinatorics.
In less figurative language, Steenrod homotopy theory is a natural framework
for doing algebraic topology on quite general spaces in terms of algebraic
topology of polyhedra.
The principal purpose of this memoir is to provide a transparent
treatment of Steenrod homotopy theory of {\it compacta}, i.e.\ compact
metrizable spaces.
An extension to residually finite-dimensional (=finitistic) metrizable
complete uniform spaces is outlined at the end of \S7.
The exposition is essentially self-contained.

The specific content of the paper is summarized at the end of the introduction.

\subhead A. Motivations \endsubhead
The author's interest in Steenrod homotopy was motivated by questions arising
from other subjects.

1) D. Rolfsen's 1972 problems on topological isotopy (=homotopy through
embeddings) in $\R^3$: Does there exist a knot that is not isotopic to
the unknot?
(Such a knot can only be wild, since every PL knot is isotopic to the unknot
by a non-locally-flat PL isotopy.)
Are isotopic PL links necessarily PL isotopic?
An affirmative answer to the latter problem would follow from the analogue
of the Vassiliev Conjecture (on completeness of finite type invariants of
knots) for ``links modulo knots'' (i.e.\ PL links up to PL isotopy) \cite{M1}.

In general, it is natural to view Rolfsen's problems as more intuitive
counterparts of the Vassiliev Conjecture, which clarify the geometric meaning
(see \cite{M1}) of a $\derlim$-type obstruction to validity of the latter.
Such an obstruction $\theta(f,g)$, defined for any pair of knots $f,g$ that
are indistinguishable by finite type invariants, arises from a reformulation
of finite type invariants in terms of isovariant homotopy theory of
appropriately compactified configuration spaces.%
\footnote{Here is a sketch.
If $M$ is a smooth manifold, the configuration space $M^{(r)}:=M^r\but$(all
diagonals) has the ``linear'' Fulton--MacPherson compactification $M^{[r]}$
(see \cite{Si}) and also a ``polynomial'' compactification $M^{\{r\}}$, which
is obtained from $M^{[r]}$ by blowing up submanifolds of increasingly
degenerate (as measured by the coranks of representative collections of
vectors) limit configurations \cite{M3}.
Every PL or piecewise-smooth knot $f\:S^1\emb S^3$ induces an $S_r$-equivariant
stratum preserving map $f^{[r]}\:(S^1)^{[r]}\to (S^3)^{[r]}$, and if
$f$ is smooth, $f^{[r]}$ is ``aligned'' in the sense of Sinha
(see \cite{BCSS}), or {\it equivalently}, lifts to an $S_r$-equivariant
stratum preserving map $f^{\{r\}}\:(S^1)^{\{r\}}\to (S^3)^{\{r\}}$,
where $(S^1)^{\{r\}}=(S^1)^{[r]}$ since $S^1$ is one-dimensional.
It now follows from the work of Sinha \cite{Si} and Voli\'c \cite{V1}, \cite{V2}
that smooth knots $f$ and $g$ are indistinguishable by rational
Vassiliev invariants if $f^{\{r\}}$ and $g^{\{r\}}$ are homotopic in the space
$V_r$ of $S_r$-equivariant stratum preserving maps
$(S^1)^{\{r\}}\to (S^3)^{\{r\}}$ for each $r$.
(It is likely that the proof can be refined to work over the integers,
cf.\ \cite{Koyt}; to this one should add that it is still unknown whether
general Vassiliev invariants reduce to rational ones.)
Conversely, if $f$ and $g$ are indistinguishable by general Vassiliev invariants
of types $<r$, by a result of Goussarov \cite{Gu} and Habiro \cite{Ha} they are
related by smooth isotopy and $C_r$-moves.
It is easy to see that each $C_r$-move looks just like an isotopy to
any $r$-tuple of points on the knot; this means that $f^{\{r\}}$ and
$g^{\{r\}}$ are related by a homotopy $h_r\:I\to V_r$.
Finally, if $f$ and $g$ are indistinguishable by all Vassiliev invariants,
each $h_r$ combines with $p_{r+1}h_{r+1}$, where
$p_r\:V_{r+1}\to V_r$ is the forgetful map, into a loop
$\ell_r\:(S^1,pt)\to (V_r,f^{\{r\}})$.
Let $\theta(f,g)$ be the class of $([\ell_1],[\ell_2],\dots)$ in
$\derlim\pi_1(V_r,f^{\{r\}})$.
If $f$ and $g$ are related by a smooth isotopy $\phi_t$, then each
$p_{r+1}\phi_t^{\{r+1\}}=\phi_t^{\{r\}}$, whence $\theta(f,g)$ is trivial.}
Besides, the author learned from different sources that some obstruction
of this kind seems to have already been considered by R. Milgram in the 90s
(unpublished).
Incidentally, the knot theorist Ralph Fox must have had very similar ideas in
mind when he switched to Steenrod homotopy in his last two papers (see
\S7).

2) Borsuk's problem of embeddability of $n$-dimensional absolute retracts in
$\R^{2n}$.
By a classical result of Shchepin--Shtan'ko, embeddability of a given
$n$-dimensional compactum $X$ in $\R^m$ in codimension $m-n\ge 3$ is equivalent
to proper embeddability of an appropriate infinite mapping telescope
$P_{[0,\infty)}$ (see \S2) into $\R^m\x [0,\infty)$; see \cite{MS} for
a precise statement.
In the metastable range $m>\frac{3(n+1)}2$, there is a complete obstruction
to such an embeddability problem in equivariant stable cohomotopy, which
reduces to an obstruction in ordinary cohomology when $m=2n$ \cite{MS}.
In particular, it turns out that the impact of the LC$_\infty$ condition
for $X$ on the embeddability of $X$ in $\R^{2n}$ is closely related to
a certain subgroup in the \v Cech $2n$-cohomology of $X\x X\but\Delta$,
measuring the non-commutativity of $\invlim$
(arising because $X\x X\but\Delta$ is non-compact) and $\dirlim$ (arising
because $X$ is not a polyhedron).
This algebra suffices to construct a counterexample to the parametric version
of the Borsuk problem, that is, an $n$-dimensional absolute retract, $n\ge 2$,
admitting non-isotopic embeddings into $\R^{2n+1}$.
The Borsuk problem itself turns out to be much harder, however.

3) R. D. Edwards' approach to the Hilbert--Smith Problem.
Does there exist a free action of the group $\Z_p$ of $p$-adic integers on
an LC$_\infty$ compactum with a finite dimensional orbit space?
It is not hard to reformulate the property of $X$ to be the orbit space of
a free action of $\Z_p$ on an LC$_\infty$ compactum in terms of local
Steenrod homotopy of $X$.
This gives us a purely topological (with no group actions involved) ---
even combinatorial topological (due to the nature of Steenrod homotopy)
problem to think about.
\medskip

\subhead B. Basic notions \endsubhead
Let $X\i\R^m$ be a compactum.
The {\it Steenrod homotopy set} $\pi_n(X,x)$ is the set of equivalence classes
of level-preserving maps
$$F\:(S^n\x [0,\infty),\,b\x [0,\infty))\to (\R^m\x[0,\infty),\,x\x [0,\infty))$$
such that the closure of $F(S^n\x [0,\infty))$ in $\R^m\x[0,\infty]$ meets
$\R^m\x\{\infty\}$ in a subset of $X\x\{\infty\}$ --- up to homotopy through
maps of the same kind.
As one would expect, $\pi_0(X,x)$ is a pointed set, $\pi_1$ is a group
(using the base ray $b\x[0,\infty)$), and $\pi_2,\pi_3,\dots$ are abelian
groups.

Let $\spi_n(X,x)$ stand for the {\it singular} homotopy set, i.e.\ the set of
homotopy classes of maps $(S^n,b)\to (X,x)$.
There is a homomorphism (pointed map when $n=0$)
$$\spi_n(X,x)@>\sstau>>\pi_n(X,x)$$
assigning to a representing spheroid $f\:(S^n,b)\to (X,x)$
the Steenrod spheroid
$f\x\id_{[0,\infty)}\:S^n\x[0,\infty)\to X\x[0,\infty)\i \R^m\x [0,\infty)$.
It not hard to see that if $X$ is a polyhedron, $\stau$ is an isomorphism.

The definition of $\pi_n(X)$ can be formulated in a purely combinatorial
way, without mention of the closure.
Let us pick a sequence
$\dots\i P_1\i P_0$ of closed polyhedral neighborhoods of $X$ in $\R^m$ with
$P_0=\R^m$ and $\bigcap P_i=X$, and consider their mapping telescope
$P_{[0,\infty)}=\bigcup P_i\x [i,i+1]$, which lies in $\R^m\x [0,\infty)$.
Then $\pi_n(X,x)$ is simply the set of level-preserving
homotopy classes of level-preserving maps
$$G\:(S^n\x [0,\infty),\,b\x [0,\infty))\to (P_{[0,\infty)},\,x\x\{0,\infty\}).$$
(Indeed, every $G$ is an $F$, and every $F$ can be made into a $G$ by
an appropriate reparametrization of $[0,\infty)$, cf.\ the proof of Lemma 2.5.)
This is a special case of the definition in \S3, which will not presuppose
an embedding into $\R^m$ nor that $X$ is finite-dimensional.

Restricting a $G$ as above to the integer levels of $S^n\x [0,\infty)$, we
get a family of maps $g_i\:(S^n,b)\to (P_i,x)$ such that each $g_{i+1}$
is homotopic to $g_i$ with values in $P_i$.
This describes an epimorphism (pointed surjection when $n=0$)
$$\pi_n(X,x)@>\check\tau>>\invlim\pi_n(P_i,x).$$
The inverse limit%
\footnote{The {\it inverse limit} $\invlim S_i$ of a sequence
$\dots@>f_1>>S_1@>f_0>>S_0$ of sets and maps is the subset of the product
$\prod_i S_i$ consisting of all sequences $(\dots,s_1,s_0)$ with
$s_i=f_i(s_{i+1})$ for all $i$.
Such sequences of points are called {\it threads}, and a sequence of sets and
maps as above is called an {\it inverse sequence}; its maps may be referred to
as the {\it bonding maps}.
An inverse limit of groups (and homomorphisms) is naturally a group, and
an inverse limit of topological spaces (and continuous maps) is naturally
a topological space.
For instance, $X$ is homeomorphic to $\invlim(\dots\i P_1\i P_0)$ in
the above situation.}
on the right hand side is also known as the ``\v Cech homotopy group''
$\check\pi_n(X,x)$.
While being an invariant of $(X,x)$, it is of limited interest as such, mainly
because it fails the homotopy exact sequence of a pair --- which does hold for
Steenrod homotopy groups.
(The relative Steenrod and \v Cech homotopy groups are defined in the
straightforward way.)

What is the kernel of $\check\tau$?
An element of the kernel would be represented by a $G$ whose
restriction to every integer level is a constant map.
Each $[i,i+1]$ layer of $G$ then gives rise to a map
$G_i\:(S^{n+1},b)\to (P_i,x)$.
The homotopy classes of these spheroids $G_1,G_2,\dots$ are not entirely
well defined, however.
For we could subtract some piece (i.e., a spheroid) from some $G_i$ and
add it to $G_{i+1}$, without changing the class of $G$ in $\pi_n(X)$.
Limiting ourselves to the abelian situation ($n>0$) for simplicity of notation,
we conclude that $\ker\check\tau$ is precisely the cokernel of the homomorphism
$$\phi\:\prod_i\pi_{n+1}(P_i)\to\prod_i\pi_{n+1}(P_i)$$
given by $(\dots,g_1,g_0)\mapsto (\dots,g_1-f_1(g_2),g_0-f_0(g_1))$, where
$f_i\:\pi_{n+1}(P_{i+1})\to\pi_n(P_i)$ are induced by the inclusions.
This cokernel is known as the {\it derived limit} $\derlim\pi_{n+1}(P_i)$ of
the inverse sequence of our abelian groups.
(Incidentally, note that the kernel of $\phi$ is nothing but
$\invlim\pi_{n+1}(P_i)$.)

So, our findings can be summarized (for $n>0$) by the short exact sequence
$$0\to\derlim\,\pi_{n+1}(P_i)\to\pi_n(X)\to\invlim\pi_n(P_i)\to 0.\tag{$*$}$$
Due to the straightforward nature of the inverse limit on the right,
the derived limit on the left is undoubtedly the central algebraic object
of Steenrod homotopy theory.

Let us now consider a (possibly unbounded) closed subset $X\i\R^m$.
Let $\dots\i P_1\i P_0$ be closed polyhedral neighborhoods of $X$ in $\R^m$
with $P_0=\R^m$ and $\bigcap P_i=X$ and such that each $P_i$ is contained in
the $\frac1{2^i}$-neighborhood of $X$ and contains
the $\frac1{2^{i+1}}$-neighborhood of $X$, with respect to the usual
Euclidean metric on $\R^m$.
The {\it Steenrod homotopy set} $\pi_n(X,x)$ is again the set of
level-preserving homotopy classes of level-preserving maps
$$G\:(S^n\x [0,\infty),\,b\x [0,\infty))\to (P_{[0,\infty)},\,x\x\{0,\infty\})$$
into the mapping telescope $P_{[0,\infty)}=\bigcup P_i\x [i,i+1]$.
If $X$ is a uniform polyhedron (see \cite{Is}), for instance, a union of
cubes of the form $[i_1,j_1]\x\dots\x [i_m,j_m]$, where $i_k,j_k\in\Z$
(possibly $i_k=j_k$), then the natural homomorphism $\spi_n(X)\to\pi_n(X)$
is again an isomorphism.
So in this case, $\pi_n(X)$ is a topological invariant of $X$.
In general, it is not: if $X$ is the pair of co-asymptotic hyperbolas
$\{(x,y)\mid y=\pm\frac1x\}$ in the plane, then $\pi_0(X)=0$.
However, $\pi_n(X)$ is always a {\it uniform} invariant of $X$, i.e.\ it is
invariant under homeomorphism that is uniformly continuous in both directions
(see \S7, where a more general definition of $\pi_n(X)$ is given, which does
not presuppose an embedding into $\R^m$).
In return, our $\pi_n(X)$ remains feasibly computable, as it still fits
into the short exact sequence ($*$).
\medskip

The said also applies to homology and cohomology.
Recall that if $K$ is a locally finite simplicial complex, in addition to
the usual simplicial chains $C_*(K)$ and cochains
$C^*(K)=\Hom(C_*(K),\Z)$ there are also the {\it compactly supported
cochains} $C^*_\c(K)=\dirlim C^*(K,K_i)$,
where the subcomplexes $\dots\i K_1\i K_0$ have complements consisting of
finitely many simplices, and $\bigcap K_i=\emptyset$,\ %
\footnote{The {\it direct limit} $\dirlim G_i$ of a sequence
$G_0@>f_0>>G_1@>f_1>>\dots$ of groups and homomorphisms is a group, whose
elements are equivalence classes of {\it threads}, i.e.\ sequences of
the form $g=(g_k,g_{k+1},\dots)$, where each $g_{i+1}=f_i(g_i)$ for $i\ge k$.
Threads $g$ and $h=(h_l,h_{l+1},\dots)$ are considered equivalent if
$g_i=h_i$ for some $i\ge\max(k,l)$ (and hence also for $i+1,i+2,\dots$).
The product of the classes of $g$ and $h$ is the class of
$(g_mh_m,g_{m+1}h_{m+1},\dots)$, where $m=\max(g,h)$.}
as well as the {\it locally-finite chains}
$C_*^\lf(K)=\invlim C_*(K,K_i)\simeq\Hom(C^*_\c(K),\Z)$.
The {\it locally-finite homology} $H_n^\lf(P)=H_n(C_*^\lf(K))$
of the locally compact polyhedron $P$, triangulated by $K$, and the
{\it compactly supported cohomology} $H^n_\lf(P)=H_{-n}(C^*_\c(K))$
are invariants of {\it proper} homotopy equivalence.
The Poincare duality for possibly non-compact orientable $m$-manifolds
reads:
$H^n(M)\simeq H_{m-n}^\lf(M)$, $H_n(M)\simeq H^{m-n}_\c(M)$
(for manifolds with boundary there will be $4$ different isomorphisms).

The groups $H^n_\c$ and $H_n^\lf$ can be very easily generalized to closed
subsets of a Euclidean space, if we use the above notation.
The {\it locally finite Steenrod homology} $H_n^\lf(X)$ is defined to be
$H_{n+1}^\lf(P_{[0,\infty)})$, and the {\it Pontryagin cohomology with compact
supports} $H^n_\c(X)$ is $H^{n+1}_\c(P_{[0,\infty)})$.
It is not hard to show that although they are defined using the uniform
neighborhoods $P_i$ of $X$, these groups turn out to be topological (rather than
just uniform) invariants of $X$.
However, since they are not homotopy invariants of $X$ (but only proper homotopy
ones), their usage is limited.

The usual {\it Steenrod homology} $H_n(X)$ is defined to be the $(n+1)$st
homology group of the ``$[0,\infty)$-locally finite'' chain complex
$C_*^{[0,\infty)}(P_{[0,\infty)})=\invlim C_*(P_{[0,\infty)},P_{[i,\infty)})$,
where $P_{[0,\infty)}$ is triangulated so that each $P_i$ is triangulated by
a subcomplex, and we write $P_J$ for the preimage of the subset
$J\i [0,\infty)$ under the projection $P_{[0,\infty)}\to [0,\infty)$.
Similarly to the above, there is a natural short exact sequence
$$0\to\derlim\, H_{n+1}(P_i)\to H_n(X)\to\invlim H_n(P_i)\to 0.\tag{$**$}$$
This exact sequence (which may be identified as the universal coefficient
formula for the left exact functor $\invlim$) is due to Milnor; by analogy,
($*$) is also known as the Milnor exact sequence.

The {\it Pontryagin cohomology} $H^n(X)$ is defined to be the $-(n+1)$st
homology group of the ``$[0,\infty)$-compactly supported'' cochain complex
$C^*_{[0,\infty)}(P_{[0,\infty)})=\dirlim C^*(P_{[0,\infty)},P_{[i,\infty)})$.
Dually to ($**$), there is a natural isomorphism
$$H^n(X)\simeq\dirlim H^n(P_i).\tag{$**'$}$$
Instead of presenting a proof (it is easier than the above proof of ($*$)
and ($**$) and can be found in \S4), let us clarify to some extent the content
of the duality.
If $C_*$ is the chain complex $$0@<<<\bigoplus G_i@<\psi<<\bigoplus
G_i@<<<0$$
where $G_i$ are abelian groups and $\psi$ sends every $g_i\in G_i$ to
$g_i-f_i(g_i)$, then obviously $H_0(C_*)\simeq\dirlim G_i$,
$H^0(C_*)\simeq\invlim G_i^\star$ and $H^1(C_*)\simeq\derlim G_i^\star$, where
$G^\star=\Hom(G,\Z)$.
There is no surprise here: like any derived functor, $\derlim$ may be viewed as
a topological object.
We shall take advantage of this viewpoint in the proofs of Theorem 3.1(d) and
Lemma 3.7(a).

\subhead C. Steenrod theories versus singular ones \endsubhead
A compactum $X$ is said to be {\it LC$_\infty$} if it is locally $n$-connected
for all $n$.
For such compacta singular homotopy and homology coincide with Steenrod
ones (see \S6), and singular cohomology coincides with Pontryagin cohomology
(see \cite{Br}).
For general compacta, however, it turns out that singular groups become somewhat
``anomalous'' and difficult to compute, whereas Steenrod groups still possess
a number of reasonable properties which make them feasibly computable and
applicable in geometric topology.
Here are a few standard illustrations of this well-known thesis:
\smallskip

1) The Alexander duality $H^n(X)\simeq\tl H_{m-n-1}(\R^m\but X)$ for compact
$X$ can use any kind of homology, but must use Pontryagin (not singular)
cohomology.
The Alexander duality $H_n(X)\simeq\tl H^{m-n-1}(\R^m\but X)$ for
compact $X$ can use any kind of cohomology, but must use Steenrod (not
singular) homology. See Theorem 4.3.
Note that an immediate corollary of either of these isomorphisms is
the Jordan curve theorem.

2) Finite-sheeted%
\footnote{and, more generally, arbitrary coverings in the category of uniform
spaces (rather than topological spaces)} coverings over Steenrod connected
(i.e.\ with trivial $\pi_0$) compacta are classified by the topologized
Steenrod fundamental group $\pi_1$ (see Corollaries 7.5 and 7.8).
They are not classified by the topologized singular fundamental group
$\spi_1$ (by considering the Warsaw circle, which is both Steenrod connected
and path connected).

3) The group $[X,\,K(G,n)]$, defined either as the group of homotopy classes
of maps $X\to K(G,n)$ or the group of Steenrod homotopy classes
(see Proposition 2.3), is isomorphic to the Pontryagin (not singular)
cohomology $H^n(X;\,G)$ for any countable abelian $G$ (see Proposition 4.2).

4) Steenrod homology and Steenrod homotopy, not to mention Pontryagin cohomology,
behave reasonably with respect to inverse limits (see ($*$), ($**$) and ($**'$)
above as well as Theorems 3.1(c) and 4.1(iii)), whereas singular homology and
homotopy of simplest inverse limits is very difficult to compute (see
Theorem 1.1 and Example 5.6).

5) Steenrod homology and Pontryagin cohomology satisfy a stronger form of
excision and are characterized on pairs of compacta as being the only
ordinary homology and cohomology theory satisfying this ``map excision'' axiom
along with an additivity axiom \cite{Mi} (see \S4 below for the precise statement).

6) Pontryagin cohomology coincides with Alexander--Spanier cohomology \cite{Sp}
and with sheaf cohomology (with constant coefficients) \cite{Br}.
Steenrod homology coincides with Massey homology (which is the dual theory
of Alexander--Spanier cohomology) \cite{Mas} and, if the coefficients are
finitely generated%
\footnote{If the coefficients are not finitely generated or not locally
constant, the Borel--Moore theory is known to have certain ``defects'',
which are discussed in detail by E. G. Sklyarenko in the Editor's comments
to the Russian translation of the first edition of Bredon's book \cite{Br}.
Instead, Sklyarenko proposes to construct a {\it co}sheaf Steenrod homology
(not to be confused with the cosheaf \v Cech homology, defined in the second
edition of Bredon's book).
Such a theory can indeed be constructed along the lines of Steenrod's second
definition of his homology (see \S4 of the present paper); the details will
appear elsewhere.}%
, with the sheaf homology of Borel and Moore (with constant coefficients)
\cite{Br}, \cite{BoM}.
Concerning all of these assertions see \cite{Sk}.

7) An $n$-dimensional compactum has trivial Steenrod homology and Pontryagin
cohomology in dimensions $>n$; whereas singular homology tells nothing about
dimension:

\proclaim{Theorem 1.1} {\rm (Barratt--Milnor \cite{BM})} Let $E^n$ be
the $n$-dimensional Hawaiian earring, that is the one-point
compactification of $\R^n\x\N$.
Then the singular $(2n-1)$-homology $\sH_{2n-1}(E^n)\ne 0$.
\endproclaim

The original proof was in the language of homotopy theory; below is a short
proof based on the Pontryagin--Thom construction and linking numbers.
It works to show that $\sH_{2n-1}(E^n;\,\Q)$ is uncountable.
There seems to be no difficulty in reproving the entire result of
\cite{BM} --- that $\sH_{kn-k+1}(E^n;\,\Q)$ is uncountable for each $k>0$
--- along the same lines by using higher Massey products.

\demo{Proof}
Let $K_0,K_1,\dots\i S^{2n-1}$ be a null-sequence of
disjoint PL copies of $S^{n-1}$ such that each $K_{2i+1}$ is
linked with $K_{2i}$ with linking number $l_i>0$ and unlinked with all
other $K_j$.
Let $T_0,T_1,\dots$ be disjoint regular neighborhoods of the $K_i$'s.
Define a map $f\:S^{2n-1}\to E^n$ by projecting each $T_i\cong S^{n-1}\x\R^n$
onto $\R^n\x\{i\}$ and sending the complement of $\bigcup T_i$ onto
the compactifying point.
We claim that $[f]\ne 0\in\sH_{2n-1}(E^n)$.

Suppose that $S^{2n-1}=\partial M$ for some compact oriented $2n$-pseudo-manifold%
\footnote{We recall that a [pseudo-connected] (orientable)
{\it $m$-pseudo-manifold} $M$ with boundary $\partial M$ is a polyhedral pair
$(M,\partial M)$ admitting a triangulation $(K,\partial K)$ such that
$[K\but K^{(m-2)}]\cup[\partial K\but (\partial K)^{(m-3)}]$ is a [connected]
(orientable) $m$-manifold with boundary $\partial K\but (\partial K)^{(m-3)}$.}
$M$ such that $f$ extends to a continuous map $F\:M\to E^n$.
Let us perturb $F$ so as to make it simplicial with respect to some
triangulations of $[-1,1]^n\x\N$ and $F^{-1}([-1,1]^n\x\N)$ where each
$x_i=(0,i)\in\R^n\x\N$ lies in the interior of an $n$-simplex.
Then each $M_i:=F^{-1}(x_i)$ is a compact oriented $n$-pseudo-manifold with
boundary $\partial M_i=K_i$, whose regular neighborhood in $M$ is homeomorphic
to $M_i\x I^n$, see \cite{Wi; Theorem 1.3.1}.
Furthermore, $M_i$ is co-oriented in $M$ as long as an orientation of
$\R^n\x\N\i E$ is fixed.
Let $z_{2i}$ be the generator of
$H^n(M,\,M\but M_{2i})\simeq H^n(M_{2i}\x I^n,\,M_{2i}\x\partial I^n)\simeq
H^0(M_{2i})$ corresponding to this co-orientation, and let $c_{2i}$ be its image
in $H^n(M)$.
Suppose that some nontrivial linear combination
$m_0c_0+\dots+m_rc_{2r}=0\in H^n(M)$.
Let $R=M_0\cup\dots\cup M_{2r}$.
Then $m_0z_0+\dots+m_rz_{2r}\in H^n(M,\,M\but R)$ is the image of some
$b\in H^{n-1}(M\but R)$.
The restriction of $b$ over
$X:=\partial M\but (K_0\cup\dots\cup K_{2r})$ equals
$m_0D[K_0]+\dots+m_rD[K_{2r}]$, where
$D\:H_{n-1}(K_0\sqcup\dots\sqcup K_{2r})\to H^{n-1}(X)$ is the Alexander
duality.
Hence $b(m_0[K_1]+\dots+m_r[K_{2r+1}])=m_0^2l_0+\dots+m_r^2l_r>0$.
On the other hand, each $[K_{2i+1}]=0\in H_{n-1}(M,\,M\but R)$
since $K_{2i+1}=\partial M_{2i+1}$.
The obtained contradiction shows that $c_{2i}$ are all linearly independent
in $H^n(M)$.
Thus $H^n(M)$ is not finitely generated, which is a contradiction. \qed
\enddemo

\subhead D. Comparison with literature \endsubhead
It should be clarified (at the referee's request) that, to the best of
the author's knowledge, the combinatorial treatment of Steenrod homotopy
theory, as presented in this paper, has no ``original exposition'' in
the literature, apart from certain fragments.%
\footnote{These include Milnor's proof of Theorem 2 \cite{Mi}, Kodama's proof
of Theorem 4 \cite{Ko2} and the first paragraph of Ferry's proof of Theorem 4
\cite{Fe2}.}
At the same time, there is no lack in expositions of various closely related
constructions, whose interrelations and relations with Steenrod homotopy
are very confusing at the beginning; the basic relations are briefly
summarized below.
We shall not dwell on these technicalities in what follows; in particular, most
of the ``known results'' mentioned in this paper are, strictly speaking, only
reformulations of such, while their equivalence with the original
formulations will sometimes be evident only to the experts.

The principal object of study in this paper is the Steenrod homotopy category,
which in the case of compacta is defined in \S2.
In this case it coincides with the strong shape category, also known as
the fine shape category \cite{EH}, \cite{DS2}, \cite{KO} (see also
\cite{Fe2; \S5}, \cite{C1}, \cite{IS} and a survey in \cite{MaS; \S III.9}).
Moreover, in the case of compacta, the original shape category of Borsuk (see
\cite{F1}, \cite{DS1}, \cite{MaS}) has the same isomorphism classes of objects
(see Proposition 2.7) and, at least in the finite dimensional case, the same
isomorphisms (see Theorem 3.9); its morphisms are generally less informative
(see examples in \S5).

In general, the strong shape category has topological spaces as its objects and
thus differs from the Steenrod homotopy category, whose objects are complete
uniform spaces.%
\footnote{This disagrees with the terminology of \cite{EH} and \cite{Po}, whose
``Steenrod homotopy'' is only used as yet another synonym of strong shape.}
There are many reasons for the latter; the most obvious one is that
the homotopy exact sequence of a fibration fails for strong shape groups,
already in the case of coverings over compacta.%
\footnote{Let $X$ be the Hawaiian snail in Example 5.7, that is the mapping
torus of the self-homeo\-morphism $\sigma_n$, $(r,i)\mapsto(r,i+1)$, of
the $n$-dimensional Hawaiian earring $E^n=(\R^n\x\Z)^+$, where $n\ge 2$.
Then by the Milnor exact sequence ($*$), $\pi_{n+1}(X)\simeq\Z[\Z]$
as a module over the group ring of the Steenrod fundamental group
$\pi_1(X)\simeq\Z$.
The universal cover $\tl X$ is (non-uniformly) homeomorphic to $E^n\x\R$, and
therefore has the strong shape of $E^n$, for which $\pi_{n+1}(E^n)=0$.
Since $E^n$ is compact, the $(n+1)$st strong shape group of $E^n$ is isomorphic
to $\pi_{n+1}(E^n)$.
Thus the $(n+1)$st strong shape group of $\tl X$ vanishes and so the universal
covering $\tl X\to X$ does not induce an isomorphism of the $(n+1)$st strong
shape groups.
At the same time, by (the non-compact case of) the same exact sequence ($*$),
the $(n+1)$st Steenrod homotopy group $\pi_{n+1}(\tl X)\simeq\Z[\Z]$; moreover,
by the naturality of $(*)$, $\pi_{n+1}(\tl X)\to\pi_{n+1}(X)$ is an isomorphism.}
In view of the results of \S7 this should not be surprising: one should rather
expect to have a homotopy exact sequence of a uniform covering --- which one
indeed has in Steenrod homotopy.
(The shape theorist might simply reply that not every infinite-sheeted
covering over a compactum would then be considered a strong shape fibration;
but for the geometric topologist, this only makes it more evident that strong
shape theory cries out to be improved.)
Another important reason is that a non-compact topological space as simple as
$[0,\infty)$ admits no cofinal sequence in the directed
set of all coverings, and so the strong shape theory of metrizable topological
spaces cannot do without non-sequential inverse spectra (indexed by directed
sets other than the positive integers).
Doing anything geometric with non-sequential mapping telescopes seems to be
barely feasible; anyway it is not needed since the Steenrod homotopy
category of metrizable complete uniform spaces manages well with inverse
sequences.

Due to constraints of space and time, the extension of the Steenrod homotopy
category to all separable metrizable complete uniform (SMCU) spaces has to be
postponed to a separate paper, as it involves a solution to Isbell's problem
\cite{Is; Research Problem B$_2$} asking, essentially, whether
infinite-dimensional uniform polyhedra can be defined in a meaningful way.
Residually finite-dimensional SMCU spaces, i.e.\ inverse limits of
finite-dimensional uniform polyhedra, are well understood after Isbell's work;
the extension of the Steenrod homotopy category to such spaces is briefly
outlined in \S7 and is applied to Steenrod homotopy of compacta in \S8.
We could have as well defined the {\it uniform} Steenrod homotopy category for
such spaces (whose restriction to uniform polyhedra is the uniform homotopy
category, rather than the usual homotopy category) but we shall not need it.
In \cite{SSG}, a peculiar {\it semi-uniform} Steenrod homotopy category has
been defined for residually finite-dimensional SMCUs, based on the notion of
a semi-uniform homotopy that is a (possibly non-uniformly-continuous)
homotopy through uniformly continuous maps.
We note that $\R$ is semi-uniformly contractible, but not uniformly contractible,
which should not be surprising since $\R$ is not a uniform absolute retract
(a notion whose definition does not involve any homotopies, cf.\ \cite{Is}).
\medskip

Steenrod homotopy and homology $\pi_i(X)$, $H_i(X)$ and Pontryagin cohomology
$H^i(X)$, defined above in the case of closed subsets $X$ of a Euclidean space,
are invariants of Steenrod homotopy type (which in particular implies that they
are well-defined as uniform invariants).
The groups $H_i^\lf(X)$ and $H^i_\c(X)$ are invariants of proper strong shape
(so in particular, of homeomorphism and of proper Steenrod homotopy type).

Steenrod homotopy groups of compacta were introduced by Christie \cite{Ch}
and rediscovered 30 years later by Quigley \cite{Q1}.
They are sometimes called ``Quigley approaching groups'' in the literature.
Steenrod homology and Pontryagin cohomology of compacta have been introduced
by Steenrod \cite{St} and (in the finite-dimensional case) Pontryagin \cite{P}
themselves.
Pontryagin also considered arbitrary closed subsets $X$ of a Euclidean space
and defined their Pontryagin cohomology with compact supports $H^i_\c(X)$ ---
by means of equation ($**'$), where the cohomology of $P_i$'s is also taken
to be with compact supports.
Since these groups are topological invariants of $X$, Pontryagin did not have
to worry about taking the neighborhoods $P_i$ uniform in order to prove
the well-definedness; but he was consistent in using only direct and inverse
sequences (indexed by the natural numbers rather than general directed sets).
Pontryagin cohomology $H^i(X)$ was considered in \cite{Do} and \cite{Miy}
(see also \cite{AgS}) under different names.
As for Steenrod homology $H_i(X)$ (and homotopy $\pi_i(X)$) in the non-compact
case, whose discussion Steenrod prudently avoided in the indicated paper, they
might be not in the literature altogether.
The locally finite Steenrod homology $H_i^\lf(X)$ and Pontryagin cohomology
with compact supports $H^i_\c(X)$ were studied by Sklyarenko (see \cite{Sk})
and others for arbitrary {\it $\sigma$-compacta}, i.e.\ locally compact
separable metrizable spaces $X$.

Pontryagin's cohomology of compacta is better known as ``\v Cech cohomology''
in the literature.
Our primary concern is to avoid confusion of three different entities:
(Alexandroff--)\v Cech(--Dowker) homology, which is not a homology theory as
it does not satisfy the exactness axiom; \v Cech(--Dowker) cohomology, which is
a cohomology theory of topological spaces; and Pontryagin cohomology, which
is a cohomology theory of complete uniform spaces.
In this paper we shall only need Pontryagin cohomology in the case of compacta,
where it coincides with the \v Cech--Dowker cohomology.
Even in this case it still seems worthwhile to call it ``Pontryagin cohomology''
as it helps eliminating the misleading (and historically unjustified!) link
with the inexact Alexandroff--\v Cech homology.%
\footnote{As for the history, it goes as follows.
In his 1931 paper \cite{P} in ``Mathematische Annalen'', Pontryagin
introduced direct limit of a direct sequence of groups, and used it to define
{\it co}homology of a finite-dimensional compactum \cite{P; Ch.\ III, \S II}.
(He did not use the term ``cohomology'', which only appeared several years
later.)
He also considered inverse sequences, but was able to define homology of
a compactum only over the rationals (in Appendix II), proving incidentally
that every inverse sequence of finite dimensional vector spaces over $\Q$
is equivalent to an inverse sequence of epimorphisms.
\v Cech was apparently unaware of this paper of Pontryagin when he wrote
his 1932 paper \cite{\v C}, which was his earliest on the subject.
Anyway he did not attempt to consider cohomology or direct limits in
that paper, but instead introduced {\it inverse} limit of an inverse spectrum
of groups and used it to clarify the meaning of Alexandroff's 1929 definition
of homology of a compactum and to extend it to additional coefficient groups
and to non-metrizable compact spaces.
(Independent equivalent definitions of homology for metric compacta were
also given by Vietoris (1927) and in the finite-dimensional case by
Lefschetz (1930), though the equivalence of Alexandroff's and Lefschetz's
approaches with that of Vietoris was only to be established much later
\cite{Le}.)
In the same paper \v Cech also attempted to define homology of non-compact
spaces, but instead defined what in the locally compact case turns out to be
the Alexandroff--\v Cech homology of the Stone--\v Cech compactification of
the original space (cf.\ \cite{ES; X.9.12}) --- which is not a homotopy
invariant of the original space.
It was not until 1950 that Dowker ``corrected'' this definition of \v Cech.}

The literature gives the fullest account of two kinds of extension of Steenrod
homology and Pontryagin cohomology of compacta to non-compact spaces: firstly,
Steeenrod--Sitnikov homology and the aforementioned \v Cech--Dowker cohomology,
defined by Sitnikov and Dowker in the 50s using the direct limit (of a possibly
non-sequential spectrum); and secondly, ``strong'' (or ``coherent'') homology
and cohomology, defined by Lisitsa--Marde\v si\v c and Miminoshvili in the 80s
using homotopy inverse limit (of a possibly non-sequential spectrum).
All these groups are strong shape invariants, as opposed to $H_i$ and $H^i$,
which are only Steenrod homotopy invariants.
This means that $H_i$ and $H^i$ can detect finer phenomena (even though they
are often easier to compute); specific examples are postponed to a subsequent
paper by the author, which will also discuss interrelations of
the aforementioned homology and cohomology theories.

\subhead E. Synopsis \endsubhead
We briefly review the content of the sections.
\medskip

{\it \S2. Steenrod homotopy category.} The section starts with our basic
setup for the Steenrod theory and its justification, due to Milnor (2.1).
Next we prove its equivalence with Christie's original 1944 formulation (2.4) by
a technique (2.5) that will be often used throughout the paper.

The leitmotif of further results of \S2 is consideration of special situations
that provide alternative viewpoints towards Steenrod homotopy as a whole.
In brief: Steenrod homotopy vs.\ pro-homotopy (2.6, 2.8b) --
Pointed vs.\ unpointed Steenrod homotopy (2.9) -- Steenrod connected vs.\ path
connected compacta (2.10) -- \v Cech extension of functors
(2.3) -- Cell-like compacta (2.7).
Overall, \S2 aims to provide a new simple approach to foundations of
the theory; we also give a short proof of the Dydak--Segal ``Fox theorem''
in Steenrod homotopy (2.8a).
\medskip

{\it \S3. Homotopy groups.} Steenrod homotopy groups are defined, and
basic tools of dealing with them are presented (3.1--3.4) and illustrated with
the aid of UV$_n$ compacta (3.5).
The ``exact continuity''%
\footnote{ñf.\ the ``Continuity versus Exactness'' section in the Eilenberg--Steenrod
book \cite{ES}}
for Steenrod homotopy groups (3.1c) seems to be a new result.
The main result of \S3 is the ``Whitehead theorem'' in Steenrod homotopy (3.6),
answering a 1983 question of A. Koyama.
We also give two simple geometric proofs of the Edwards--Geoghegan stability
theorem (3.10) and a simple proof of Dydak and Kodama's theorem on stability
of a union (3.14).

Finally, Steenrod fibrations are introduced and are found to be a simplified but
equivalent version of the shape fibrations of Marde\v si\'c--Rushing.
This enables us to simplify Cathey's proofs of the exact sequence in Steenrod
homotopy (3.15b) and Steenrod homotopy equivalence of fibers (3.15c) for
Steenrod fibrations; and
Marde\v si\'c and Rushing's proof that cell-like maps between finite-dimensional
compacta are Steenrod fibrations (3.15d).
We also show that a Steenrod fibration over a connected compactum induces
a surjection on $\pi_0$ (3.15b$'$).
\medskip

{\it \S4. Homology and cohomology.} We briefly recall the definition and
simplest properties of Steenrod homology and Pontryagin cohomology groups
(4.1, 4.2).
As an illustration, we present an elementary proof of the Alexander duality
theorems (4.3).
A clean proof of the (absolute) ``Hurewicz theorem'' in Steenrod homotopy is
given (4.4), along with a counterexample to a previously available ``proof''.
\medskip

{\it \S5. Some examples.} With some exceptions, this and further sections
presuppose only 2.1--2.5, 3.1--3.4, 4.1--4.2 and are independent of each other.
\S5 is a little zoo of compacta highlighting various features of Steenrod
homotopy and homology as opposed to singular theories.
These include the solenoid (5.3), the Hawaiian earring (5.6) and an alternative
implementation of the Dydak and Zdravkovska example (5.8).
\medskip

{\it \S6. Comparison of theories.}
This section is devoted to comparison of singular, Steenrod and \v Cech
theories under various local connectivity hypotheses, as well as to comparison
of these very hypotheses.
The main result of \S6 is a homotopy realization theorem (6.5).
We also prove a parallel LC$_{n-1}$ homology realization theorem (6.7), which is
originally due to Shchepin (unpublished) and implies a 1999 result of Eda and
Kawamura.
Next, we give geometric proofs of Marde\v si\'c's (6.8d) and Jussila's (6.8b)
comparison theorems and prove two further unpublished results by Shchepin:
realization of Steenrod cycles by fractal pseudo-manifolds for HLC$_{n-1}$
compacta (6.10) and the HLC$_n$ = HL\v C$_n$ theorem (6.8a).
On our way to these results, we reprove the homotopy comparison theorem of
Hurewicz, Borsuk, et al.\ (6.1) and Ferry's ``Eventual Hurewicz theorem'' (6.3).

As an application, we propose a geometric ``explanation'' of the finite
generation of Steenrod homology of HLC$_\infty$ compacta (6.11).
\medskip

{\it \S7. Covering theory.}
We revisit Fox's classification of overlayings, simplifying its statement and
proof (7.4) and observe that overlayings coincide with uniform covering
maps (7.6).
The notion of universal generalized overlaying is illustrated with
the ``Nottingham compactum'' (7.12).
\medskip

{\it \S8. Zero-dimensional homotopy.}
This section is devoted to Steenrod connected (=pointed 1-movable) compacta.
We obtain a simple proof of the Geoghegan--Krasinkiewicz theorem that virtually
Steenrod disconnected compacta have empty Steenrod components (8.4).
Brin's counterexample (apparently unpublished) to its converse is apparently
recovered (8.6).
Next, we give a relatively simple proof of the Krasinkiewicz--Minc theorem
(8.7), which implies, in particular, that path connected compacta are Steenrod
connected.
As a consequence, a characterization of Steenrod connected compacta in terms of
the cardinality of uniform components of their uniform covering spaces is
obtained (8.9).
We also include simplified proofs of the McMillan--Krasinkiewicz theorem
on continuous images (8.2) and Krasinkiewicz's theorem on unions (8.3) of
Steenrod connected compacta.
\medskip

{\it ``Algebraic topology'' content of the paper.}
A ``classification'' of sequential countable pro-groups in terms of $\invlim$
and $\derlim$ is obtained (3.8), answering a question of Koyama.
This is based on a lemma that ``the vanishing of $\derlim$ implies
the Mittag-Leffler condition for sequential countable pro-cosets'' (3.7b).
The paper includes simple expositions of known results about $\derlim$
and the Mittag-Leffler condition (3.1bd, 3.2, 3.3, 3.4, 3.7a, 3.11).
There also is a simple geometric construction of J. H. C. Whitehead's long exact
sequence into which the Hurewicz homomorphisms in different dimensions fit (4.5).
\medskip

{\it ``General topology'' content of the paper.}
We introduce a notion of convergent inverse sequence of metrizable uniform spaces,
which generalizes that of a Mittag-Leffler inverse sequence of sets (viewed as
discrete uniform spaces), and presents an ``explanation'' of Bourbaki's
``Mittag-Leffler Theorem'' (7.10) and consequently of the Baire Category
Theorem (7.11).

\definition{Acknowledgements}
This work would have never appeared without (i) a conversation with J. Dydak,
who threw me into confusion by asking whether I ever gave a serious
thought to the Hilbert--Smith Problem; (ii) S. P. Novikov's 70th Anniversary,
which caused me to begin writing a 5-page note on the subject of
the Hilbert--Smith Problem, intended for a volume dedicated to that event;
(iii) a correspondence with
A. N. Dranishnikov, which helped me to realize that I needed to stop thinking
about the Hilbert--Smith Problem until I had a better understanding of Steenrod
homotopy; (iv) the prudence of the Editors of the aforementioned volume:
V. M. Buchstaber, O. K. Sheinman and E. V. Shchepin --- whose wisely calculated
shifting of deadlines magically transformed the aforementioned short note into
the present paper.

I'm grateful to the topology group of the University of Tennessee for their
hospitality and encouragement during my visit there; and specifically to
Jurek Dydak for very many helpful remarks and references to the literature and
to Kolya Brodskiy, Pepe (Jos\'e) Higes, Brendon LaBuz and Conrad Plaut for
useful discussions and remarks on uniform spaces and uniform covering maps.
Conversations with P. Akhmetiev, A. N. Dranishnikov, A. V. Chernavskij,
R. Mikhailov, Yu.\ B. Rudyak, Yu.\ Turygin and E. V. Shchepin taking place
over the years helped me to better understand Steenrod homology.
I am particularly indebted to Misha Skopenkov for attentive listening that helped
reveal two errors in time and for remarks on the text of the paper.
In addition, I am grateful to N. Mazurenko, who helped translate the original
English text into Russian.
\enddefinition

\head 2. Steenrod homotopy category \endhead

Let $X$ be a compactum.
Then $X$ is the inverse limit of an inverse sequence
$P=(\dots@>p_1>>P_1@>p_0>>P_0)$ of compact polyhedra and PL maps (H. Freudenthal,
1937; a concise and clear exposition of the proof can be found in
\cite{Is; V.33}).
The {\it mapping telescope} $P_{[0,\infty)}$ is the infinite union
$\dots\cup_{P_2}MC(p_2)\cup_{P_1}MC(p_1)$ of the mapping cylinders of
the bonding maps.
It comes endowed with the projection $\pi_P\:P_{[0,\infty)}\to [0,\infty)$,
sending each $MC(p_i)$ onto $[i,i+1]$ in a level-preserving fashion.
For each $J\i [0,\infty)$, we denote $\pi_P^{-1}(J)$ by $P_J$.
We shall assume that $P_0$ is a point; this does not restrict generality
as we can always augment the given inverse sequence by the map $P_0\to pt$ and
shift the indices by one (these operations do not affect the inverse limit).
This makes each $P_{[0,k]}$ contractible.

The inverse limit $P_{[0,\infty]}$ of the retractions
$r_k\:P_{[0,k+1]}\to P_{[0,k]}$ is easily seen to be a compactification of
$P_{[0,\infty)}$ by $X$.
The projections $\Pi_k\:P_{[0,\infty]}\to P_{[0,k]}$ combine with
the homotopies $r_k\simeq\id_{P_{[0,k+1]}}$ into a homotopy
$\Pi_t\:P_{[0,\infty]}\to P_{[0,\infty]}$, $t\in [0,\infty]$, where
$\Pi_\infty$ is the identity.
$P_{[0,\infty]}$ is LC$_\infty$ since it $\eps$-deformation retracts
(via $\Pi_{1/\eps}$) onto a compact polyhedron for each $\eps>0$.%
\footnote{In fact, it is well-known (but not needed for our purposes) that
$P_{[0,\infty]}$ is an ANR (and consequently an AR) by using either
Lefschetz's or Hanner's characterization of ANRs (found respectively in
Borsuk's and Hu's ``Theories of Retracts'').}

Let $Y$ be another compactum, and let
$Q=(\dots@>q^2_1>>Q_1@>q^1_0>>Q_0)$ be an inverse sequence of compact polyhedra
and PL maps ($Q_0=pt$) with inverse limit $Y$.

\proclaim{Lemma 2.1} (a) {\rm (Milnor)} Every map $f\:X\to Y$ extends to a
continuous map $f_{[0,\infty]}\:P_{[0,\infty]}\to Q_{[0,\infty]}$ that

\smallskip
(i) agrees with $f$ on $X$ and

\smallskip
(ii) restricts to a proper map
$f_{[0,\infty)}\:P_{[0,\infty)}\to Q_{[0,\infty)}$.

\smallskip
\noindent
(b) Every two maps satisfying (i) and (ii) are homotopic through maps
satisfying (i) and (ii).
\endproclaim

Part (a) is essentially \cite{Mi; Theorem 2}.
While Milnor's paper remained unpublished, the ANR $P_{[0,\infty]}$ and Lemma 2.1
were rediscovered in mid-70s \cite{Ko1}, \cite{K1}, \cite{CS}; resp.\
\cite{Ko2}, \cite{KOW} (see also \cite{DS2; 3.4}).
The ``discrete telescope'' $P_{\N\cup\infty}\i P_{[0,\infty]}$ and the infinite
mapping telescope $P_{[0,\infty)}$ were considered already by
H. Freudenthal (cf\. \cite{K1}) and S. Lefschetz (cf\. \cite{Mi}) in 1930s.
The following proof of Lemma 2.1 simplifies Krasinkiewicz's approach in \cite{K1}.

\demo{Proof} Let us triangulate $P_{[0,\infty)}$ so that the diameters
of simplices tend to zero as they approach $X$, and let
$r\:X\cup P_{[0,\infty)}^{(0)}\to X$ be a retraction, where
$P_{[0,\infty)}^{(0)}$ denotes the $0$-skeleton of the triangulation.
Let $r_J$ be the restriction of $r$ to $P_J^{(0)}$.
Since $q^\infty_i$, $q^{i+1}_i$ and $r$ are uniformly continuous, for each $i$
there exists a $j=j(i)$ such that the composition
$P_{[j(i),\infty)}^{(0)}@>r_{[j(i),\infty)}>>X@>q^\infty_i>>Q_i$ extends
``by linearity'' to a map $\phi_i\:P_{[j(i),\infty)}\to Q_i$, and moreover
the composition
$P_{[j(i+1),\infty)}@>\phi_{i+1}>>Q_{i+1}@>q^{i+1}_i>>Q_i$
is homotopic ``by linearity'' to the restriction of $\phi_i$.

Combining the appropriate restrictions of these maps and homotopies, we obtain
an extension of $f$ to a continuous map $F\:P_{[k,\infty]}\to Q_{[0,\infty]}$,
where $k=j(1)$, moreover $F(P_{[k,\infty)})\i Q_{[0,\infty)}$ by construction.
The compact set $F(P_k)$ lies in some $Q_{[0,l]}$, which is contractible.
Hence $F$ extends to a map $f_{[0,\infty]}\:P_{[0,\infty]}\to Q_{[0,\infty]}$
satisfying (i) and (ii). \qed
\enddemo

\demo{(b)} Since $P_{[0,\infty]}$ is LC$_0$, the retraction $r$ in the proof
of (a) is homotopic to the identity keeping $X$ fixed.
Since $(P_{[0,\infty]},\,X\cup P_{[0,\infty)}^{(0)})$ is a Borsuk pair (=the
inclusion $X\cup P_{[0,\infty)}^{(0)}\emb P_{[0,\infty]}$ is a cofibration),
$r$ extends to a continuous self-map
$R\:P_{[0,\infty]}\to P_{[0,\infty]}$, homotopic to the identity keeping $X$
fixed.
Then the given maps $F,G\:P_{[0,\infty]}\to Q_{[0,\infty]}$ are homotopic to
$FR$ and $GR$ keeping $X$ fixed.
Given an increasing sequence $k_1,k_2,\dots$, the telescopic bonding map
$[q^{k_0}_0,\,q^{k_\infty}_\infty]\:Q_{[k_0,\infty]}\to Q_{[0,\infty]}$ is
homotopic to the identity keeping $X$ fixed.
Finally, $[q^{k_0}_0,\,q^{k_\infty}_\infty]FR$ is homotopic $\rel X$ to
$[q^{k_0}_0,\,q^{k_\infty}_\infty]GR$ for an appropriate sequence $(k_i)$ by
the proof of part (a), applied to $f\x\id_I\:X\x I\to Y$. \qed
\enddemo

\proclaim{Corollary 2.2} When $Y=X$,
$(\id_X)_{[0,\infty)}\:P_{[0,\infty)}\to Q_{[0,\infty)}$ is a proper homotopy
equivalence.
\endproclaim

\definition{Steenrod homotopy category}
We define a {\it Steenrod homotopy class} $X\ssm Y$ to be the proper homotopy
class of a proper map $f\:P_{[0,\infty)}\to Q_{[0,\infty)}$.
By virtue of the proper homotopy equivalences $(\id_X)_{[0,\infty)}$ and
$(\id_Y)_{[0,\infty)}$, this definition does not depend on the choice of
the inverse sequences $P$ and $Q$.
The composition $X\ssm Y\ssm Z$ is defined by composing appropriate proper
maps.
Assigning to a map $f\:X\to Y$ the proper homotopy class, denoted $[f]$, of
$f_{[0,\infty)}$, yields a map $\stau$ from the set $[X,Y]^\tr$
of homotopy classes of maps $X\to Y$ into the set $[X,Y]$ of Steenrod homotopy
classes $X\ssm Y$, which is natural in both variables.
A Steenrod homotopy class $f\:X\ssm Y$ is a {\it Steenrod homotopy equivalence}
if there exists a $g\:Y\ssm X$ such that $gf=[\id_X]$ and $fg=[\id_Y]$; in
other words, if $f$ is represented by a proper homotopy equivalence.
When such an $f$ exists, $X$ and $Y$ are said to be of the same {\it shape}
(\underbar{S}teenrod \underbar{h}omotopy ty\underbar{pe}).

Steenrod homotopy classes were introduced (under a different name) by
D. E. Christie in his 1944 dissertation \cite{Ch}.
He gave several equivalent definitions and found uncountably many Steenrod
homotopy classes from a point into the dyadic solenoid (concerning these
see Example 5.3).
Christie's versions of the definition (see Proposition 2.4(a) below) as well
as Ferry's version \cite{Fe2; \S5} (see Proposition 2.4(b)) are slightly
simpler in that they partially bypass Lemma 2.1, but on the other hand their
asymmetry does not allow to define composition of Steenrod homotopy classes.
It appears that Christie himself was unaware of Steenrod's 1940 paper \cite{St},
where Steenrod homology was introduced (and, among other things, uncountably many
$0$-dimensional Steenrod homology classes of the dyadic solenoid were found).
However, the paper of Steenrod is found in the reference list of the 1942 book
\cite{Le} by Christie's thesis advisor Lefschetz.

On several occasions (in 3.10, 4.2, 7.4) we will need Steenrod homotopy classes
from a compactum to a possibly non-compact polyhedron, such as the classifying
space%
\footnote{Recall that every countable CW-complex is homotopy equivalent to
a locally finite polyhedron.
For instance, the homotopy type of $K(\Z/2,\,1)$ is represented by the direct
mapping telescope of the sequence of inclusions $\R P^1\i\R P^2\i\dots$.}
of a countable group.
If $P=(\dots\to P_1\to P_0)$ and $Q=(\dots\to Q_1\to Q_0)$ are inverse
sequences of (possibly non-compact) polyhedra, we call a map
$Q_{[0,\infty)}\to R_{[0,\infty)}$ between the two infinite mapping
telescopes {\it semi-proper} if for each $k$ there exists an $l$ such that
$f^{-1}(R_{[0,k]})\i Q_{[0,l]}$.
Now if all the $P_i$ are compact and $X=\invlim P_i$; and $Q_0=pt$ and $Q_i=K$
for all $i>0$, where $K$ is a non-compact polyhedron, we define
a {\it Steenrod homotopy class} $X\ssm K$ to be the semi-proper homotopy class
of a semi-proper map $F\:P_{[0,\infty)}\to Q_{[0,\infty)}$.
An equivalent definition will be given in \S7 in the case where $K$ is
finite-dimensional.
We write $[X,K]$ for the set of Steenrod homotopy classes $X\ssm K$.
Lemma 2.1 yields a map $\stau\:[X,K]^\tr\to[X,K]$.
\enddefinition

\proclaim{Proposition 2.3} (a) If $K$ is a polyhedron,
$\stau\:[X,K]^\tr\to [X,K]$ is a bijection.

(b) If $K$ is a polyhedron and $X$ is the inverse limit of polyhedra $P_i$,
then there is a natural bijection $[X,\,K]\to\dirlim [P_i,\,K]$.
\endproclaim

Part (a) is found, for instance, in \cite{DS2; 4.4}, and (b) is also well-known.

\demo{Proof} (a) Since $K$ is a polyhedron, $Q_{[0,\infty)}$ can be chosen to
be the open cone $CK\cup_{K=K\x 0}K\x [0,\infty)$.
Every semi-proper map of $P_{[0,\infty)}$ to this open cone restricts to a map
$P_{[k,\infty)}\to K\x [0,\infty)\simeq K$ for some $k$.
On the other hand, every $f\:P_{[k,\infty)}\to K$ is homotopic to an $f'$
satisfying $f'|_{P_{i+1}}=(f'|_{P_i})p_i$ for each $i$, and therefore
extending by continuity over $X$.
By repeating the same construction for homotopies, we obtain a composition
$[X,K]@>\phi>>\dirlim S_i@>\psi>>\dirlim [P_{[k,\infty)},\,K]@>\chi>>
[X,K]^\tr$, where $S_k$ is the set of semi-proper homotopy classes
of semi-proper maps $P_{[k,\infty)}\to K\x [0,\infty)$.
On the other hand, $\stau$ factors as
$[X,K]^\tr@>\chi'>>\dirlim [P_{[k,\infty)},\,K\x[0,\infty]]
@>\psi'>>\dirlim S_k @>\phi'>>[X,K]$ by construction.
It is easy to check that $\phi'=\phi^{-1}$, $\psi'=\psi^{-1}$ and
$\chi'=\chi^{-1}$.
In particular, to see that $\chi'\chi=\id$, note that if two maps
$P_{[k,\infty]}\to K$ agree on $X$, their restrictions to some
$P_{[l,\infty]}$ are homotopic. \qed
\enddemo

\demo{(b)} Each $P_i$ is a deformation retract of $P_{[i,\infty)}$, so
$\dirlim [P_i,K]=\dirlim [P_{[i,\infty)},K]$.
The latter coincides with $[X,K]$ by the proof of (a). \qed
\enddemo

Let us retreat to the case of Steenrod homotopy classes from a compactum to
a compactum.

\proclaim{Proposition 2.4} Let $Y$ be the limit of an inverse sequence of
compact polyhedra $Q_i$, where $Q_0$ is not assumed to be a point.
Steenrod homotopy classes $X\ssm Y$ correspond bijectively to:

(a) classes of level preserving maps $X\x[0,\infty)\to Q_{[0,\infty)}$
up to level preserving homotopy;

(b) proper homotopy classes of proper maps $X\x[0,\infty)\to Q_{[0,\infty)}$.
\endproclaim

Note that the (most useful) case where $X$ is a polyhedron is obvious,
modulo Lemma 2.5(a) below, since $P_{[0,\infty)}$ can be chosen to be the open
cone $CX\cup_{X=X\x\{0\}}X\x [0,\infty)$ and then $CX$ can be discarded
along with $Q_{[-1,0]}:=CQ_0$.

The statement of part (b) is implicit in the definitions of \cite{Fe2; \S5},
and it might be possible to extract a proof from \cite{DS2} (see in particular
Lemma 3.4 there).

\demo{Proof} If $Q_0=pt$, a Steenrod homotopy class immediately yields a proper
homotopy class as in (b), and by Lemma 2.5(a) below, a proper homotopy class
as in (b) yields a level-preserving homotopy class as in (a).
The general case reduces to the case $Q_0=pt$ similarly to Lemma 2.5(a).

Let $f_{[0,\infty)}\:X\x[0,\infty)\to Q_{[0,\infty)}$ be level-preserving.
It follows from Lemma 2.1(a) that each $f_i\:X\to Q_i$ factors through
a map $F_i\:P_{[k_i,\infty]}\to Q_i$ for some $k_i$.
From Lemma 2.1(b), the composition
$P_{[k_{i+1},\infty]}@>F_{i+1}>>Q_{i+1}@>q^{i+1}_i>>Q_i$ restricted to
$P_{[l_i,\infty]}$ for some $l_i\le\min(k_i,k_{i+1})$ is homotopic
to the restriction of $F_i$.
Without loss of generality, the sequence $l=(l_1,l_2,\dots)$ is non-decreasing.
Then $f_{[0,\infty)}$ factors through a level-preserving map
$R_{[0,\infty)}\to Q_{[0,\infty)}$, where $R_i=P_{[l_i,\infty]}$ with bonding
maps $R_{i+1}\to R_i$ being the inclusions.
Then $f_{[0,\infty)}$ also factors through a level-preserving map
$f'_{[0,\infty)}\:P^l_{[0,\infty)}\to Q_{[0,\infty)}$, where $P^l_i=P_{l_i}$.

A similar argument shows that if $f_{[0,\infty)}$ is level-preserving homotopic
to $g_{[0,\infty)}$, then $f'_{[0,\infty)}$ and the similarly constructed
$g'\:P^m_{[0,\infty)}\to Q_{[0,\infty)}$ are level-preserving homotopic
when precomposed with the level-preserving maps
$P^n_{[0,\infty)}\to P^l_{[0,\infty)}$ and
$P^n_{[0,\infty)}\to P^m_{[0,\infty)}$ for some sequence $n=(n_i)$ majorizing
both $l$ and $m$. \qed
\enddemo

\proclaim{Lemma 2.5}
(a) Every proper map $X\x [0,\infty)\to Q_{[0,\infty)}$ is
proper homotopic to a level-preserving map.
Proper homotopic level-preserving maps $X\x [0,\infty)\to Q_{[0,\infty)}$
are level-preserving homotopic.

(b$_0$) Every proper map $f\:P_{[0,\infty)}\to Q_{[0,\infty)}$ is properly
homotopic to an $f'$ such that for some infinite subsequence $P^k_i=P_{k_i}$,
the composition of the ``reindexing'' proper homotopy equivalence
$k_{[0,\infty)}\:P^k_{[0,\infty)}\to P_{[0,\infty)}$ and $f'$ is
level-preserving.

(b$_1$) If level-preserving maps $f,g\:P_{[0,\infty)}\to Q_{[0,\infty)}$
are properly homotopic, there exists an infinite increasing sequence $k=(k_i)$
such that the compositions of
$f^k,g^k\:P^k_{[0,\infty)}\to Q^k_{[0,\infty)}$ with the level-preserving map
$[p^{k_0}_0,p^{k_\infty}_\infty)\:Q^k_{[0,\infty)}\to Q_{[0,\infty)}$ are
level-preserving homotopic.
\endproclaim

Parts (b$_0$) and (b$_1$) seem to be parallel to \cite{EH; 3.7.20}.

\demo{Proof. (a)} If $\Phi\:X\x[0,\infty)\to Q_{[0,\infty)}$ is proper,
$\phi\:[0,\infty)\to [0,\infty)$ defined by
$\phi(t)=\inf (p(\Phi(X\x [t,\infty))))$ is a proper nondecreasing continuous
function, so can be approximated by a homeomorphism $\tl\phi$ such that
$\tl\phi(t)\le\phi(t)$ for every $t$.
Then $\Phi$ is properly homotopic to $\Phi'$ defined by
$\Phi'(x,t)=\Pi_{\tl\phi(t)}(\Phi(x,t))$.
Finally, $\Phi'$ is homotopic to $\Phi''$ defined by
$\Phi''(x,t)=\Phi'(x,\tl\phi^{-1}(t))$, which is level-preserving.

To prove the second assertion it suffices to observe that the homotopy from
$\Phi$ to $\Phi''$ may be assumed to depend continuously on $\Phi$ and to be
the identity whenever $\Phi$ is level-preserving. \qed
\enddemo

\demo{(b)} This is similar to (a).
The proper map $k_{[0,\infty)}$ replaces the passage from $\Phi'$ to
$\Phi''$, and the level-preserving map $[p^{k_0}_0,p^{k_\infty}_\infty)$
arises on passing from $\Phi$ to $\Phi'$. \qed
\enddemo

\proclaim{Proposition 2.6} Let $X$ and $Y$ be the limits of inverse sequences
of compact polyhedra $P=(\dots@>p_1>>P_1@>p_0>>P_0)$ and
$Q=(\dots@>p_1>>P_1@>p_0>>P_0)$.
The following are equivalent:

(i) $X$ and $Y$ are of the same shape;

(ii) $P\simeq Q$ in pro-homotopy, that is, there exist an increasing sequence
$k\:\N\to\N$ and level-preserving maps
$f_{[0,\infty)}\:P^k_{[0,\infty)}\to Q^k_{[0,\infty)}$ and
$g_{[0,\infty)}\:Q^{ks}_{[0,\infty)}\to P^k_{[0,\infty)}$, where $P^k_i$
stands for $P_{k(i)}$ and $s(i)=i+1$, such that
$g_if_{i+1}\:P_{k(i+1)}\to P_{k(i)}$ and $f_ig_i\:Q_{k(i+1)}\to Q_{k(i)}$
are homotopic to the bonding maps for all $i$;

(iii) $P$ and $Q$ are equivalent under the equivalence relation generated by
passing to subsequences and homotoping the bonding maps.

(iv) $P_{[0,\infty)}$ and $Q_{[0,\infty)}$ have the same proper simple
homotopy type.
\endproclaim

The proof that (ii)\imp(i) traces back to \cite{EH; 5.2.9}.
The equivalence relation of (iii) was introduced by Siebenmann \cite{S2}.

\demo{Proof} (iii)\imp (iv) is straightforward (see details in \cite{S2}),
and (iv)\imp (i) is immediate from definitions.
(i)\imp(ii) follows using Lemma 2.5(b).
Finally, assuming (ii), both $P$ and $Q$ are equivalent in the sense of (iii)
to the combined inverse sequence
$\dots@>g_1>>P_{k(1)}@>f_1>>Q_{k(1)}@>g_0>>P_{k(0)}@>f_0>>Q_{k(0)}$. \qed
\enddemo

\proclaim{Proposition 2.7} {\rm (Lacher \cite{La})} A compactum $X$ has the shape
of a point iff it is the limit of an inverse sequence of maps between PL balls.
\endproclaim

Compacta that satisfy any of the equivalent conditions of Proposition 2.7 are
called {\it cell-like} in the literature.
Maps with cell-like point-inverses will be discussed in Theorem 3.15(d,e).

\demo{Proof} Suppose that $X$ is the inverse limit of polyhedra $P_i$.
To prove the ``if'' direction, it suffices to note that if $P_i$ is
contractible, $f_i\:P_{i+1}\to P_i$ is the bonding map and $h_i$ is
a null-homotopy of $\id_{P_i}$, the compositions $h_if_i$ and $f_ih_{i+1}$ are
homotopic.

Conversely, suppose that $X$ has the shape of a point.
Then by Proposition 2.6, for each $i$ there exists a $j=j(i)>i$ such that
the bonding map $P_j\to P_i$ is null-homotopic.
Then it factors through the inclusion of $P_j$ into the cone $CP_j$.
Let $j_0=0$ and $j_{k+1}=j(j_k)$, and set $Q_k=CP_{j_k}$.
Then $X$ is the inverse limit of the maps $Q_{k+1}\to P_{j_k}\i Q_k$.
Each cone $Q_k$ is collapsible, so a regular neighborhood of a copy of it
in some Euclidean space is a PL ball $B^{n_k}$.
Then $X$ is the inverse limit of the maps
$B^{n_{k+1}}\to Q_{k+1}\to Q_k\i B^{n_k}$. \qed
\enddemo

Note that Proposition 2.6 does not assert that $f\:X\ssm Y$ is a Steenrod homotopy
equivalence iff it is represented by an isomorphism in pro-homotopy.
This is known to be true in the finite dimensional case (see Theorem 3.9);
the general case remains open (see \cite{DS2} and \cite{G\"u} for partial
results).
Another result in this direction is the following version ``for maps'' of
Proposition 2.6.

\proclaim{Proposition 2.8} Let $X$ and $Y$ be compacta.
\medskip
(a) {\rm (Dydak--Segal)} Every Steenrod homotopy class $f\:X\ssm Y$ factors
into the composition $h[i]$ of an embedding $i\:X\emb Z$ and a Steenrod strong
deformation retraction (so, in particular, a Steenrod homotopy equivalence)
$h\:Z\ssm Y$.
\medskip

(b) Let $f\:X\emb Y$ be an inclusion map.
The following are equivalent:

\medskip\noindent
(i) $f$ is a Steenrod homotopy equivalence;
\smallskip\noindent
(ii) $X$ is a Steenrod strong deformation retract of $Y$;
\smallskip\noindent
(iii) $(X,Y)$ is an inverse limit of compact polyhedral pairs $(Q_i,\,P_i)$
such that, writing $f_{[0,\infty)}\:P_{[0,\infty)}\to Q_{[0,\infty)}$ for
the inclusion map, there exists a level-preserving map
$g_{[0,\infty)}\:Q^s_{[0,\infty)}\to P_{[0,\infty)}$, where $Q^s_i=Q_{i+1}$,
such that each $g_if_{i+1}\:P_{i+1}\to P_i$ equals the bonding map, and each
$f_ig_i\:Q_{i+1}\to Q_i$ is homotopic $\rel P_{i+1}$ to the bonding map.
\endproclaim

A Steenrod homotopy class $h\:Z\ssm Y$ is a {\it Steenrod strong deformation
retraction} if $Y\i Z$, and the pair $(Z,Y)$ is an inverse limit of compact
polyhedral pairs $(Q_i,\,P_i)$ such that $h$ is represented by a proper strong
deformation retration $Q_{[0,\infty)}\to P_{[0,\infty)}$.

Part (a) was proved in \cite{DS2} by a different construction, occupying
pp. 24--26 there.
The implication (i)\imp(ii) is related to \cite{C1; 1.15}
(see also \cite{DS2; 4.5, 6.3$'$}).
The idea of the ``algebraic'' proof of (iii)\imp(i) goes back to
\cite{DS2; 6.2}, and the construction used in the final part of the ``geometric''
proof is found in \cite{EH; proof of 3.7.4}.

\demo{Proof. (a)} By Lemma 2.5(b), $f$ can be represented by
a level-preserving map $f_{[0,\infty)}\:P_{[0,\infty)}\to Q_{[0,\infty)}$,
where $\invlim P_i=X$ and $\invlim Q_i=Y$.
Every bonding map $P_{i+1}\to P_i$ extends to a map
$MC(f_{i+1})\to MC(f_i)\cup Q_{[i,i+1]}\to MC(f_i)$ between the mapping
cylinders, which also restricts to the bonding map $Q_{i+1}\to Q_i$.
Writing $R_i=MC(f_i)$, clearly $R_{[0,\infty)}$ proper strong deformation
retracts onto $Q_{[0,\infty)}$.
The retraction $h_{[0,\infty)}\:R_{[0,\infty)}\to Q_{[0,\infty)}$ represents
a Steenrod strong deformation retraction $h\:Z\ssm Y$, where
$Z=\invlim R_i$.
Clearly, $f_{[0,\infty)}$ is properly homotopic to the composition of
the inclusion map $P_{[0,\infty)}\emb R_{[0,\infty)}$ with $h_{[0,\infty)}$.
\qed
\enddemo

\demo{(b). (i)\imp(ii)\imp(iii)} Let us represent $f$ by a proper map
$f_{[0,\infty)}\:P_{[0,\infty)}\to Q_{[0,\infty)}$.
Assuming (i), $f_{[0,\infty)}$ is a proper homotopy equivalence.
Then by the proofs of standard results of homotopy extension theory
\cite{Sp; 1.4.10, 1.4.11, 1.D.2}, there exists a proper strong deformation
retraction of $Q_{[0,\infty)}$ onto $P_{[0,\infty)}$.
This proves (ii).
The implication (ii)\imp(iii) follows using Lemma 2.5(b). \qed
\enddemo

\demo{(iii)\imp (i). Algebraic proof}
Writing $(Q,P)=(\dots\to(Q_1,P_1)\to(Q_0,P_0))$, the hypothesis implies that
$(Q,P)$ is isomorphic to $(P,P)$ in pro-homotopy.
Then by the proof of Proposition 2.6, $(Y,X)$ has the shape of $(X,X)$.
If $\phi\:(Y,X)\ssm(X,X)$ and $\psi\:(X,X)\ssm(Y,X)$ are such that
$\phi\psi=[\id_{(X,X)}]$ and $\psi\phi=[\id_{(Y,X)}]$, then
$\psi\:X\ssm Y$, $\phi\:X\ssm X$ and the inclusion $i\:Y\to X$ satisfy
$\psi\phi[i]=[\id_Y]$ and $[i]\psi\phi=[\id_X]$.
Thus $[i]$ is a Steenrod homotopy equivalence. \qed
\enddemo

\demo{Geometric proof} Since $g_i|_{P_{i+1}}$ equals the bond, the mapping
cylinder $MC(g_i)$ contains $P_{[i,i+1]}$.
The homotopy between $f_ig_i$ and the bonding map yields a homotopy
equivalence $\phi\:MC(g_i)\to Q_{[i,i+1]}$ restricting to the identity on
$Q_i\cup P_{[i,i+1]}$.
Now $Q'_i:=MC(g_i)$ and $P'_i:=P_{[i,i+1]}$ can be arranged into inverse
sequences with bonding maps $Q'_{i+1}@>\pi\phi>>Q_{i+1}\i Q'_i$ and
$P'_{i+1}@>\pi>>P_{i+1}\i P'_i$, where each $\pi$ projects the mapping cylinder
onto its target space.
Since $\dots\to(Q'_1,\,P'_1)\to(Q'_0,P'_0)$ is equivalent to
$\dots\to(Q_1,\,P_1)\to(Q_0,P_0)$ under the equivalence relation of
Proposition 2.6(iii), the inclusion map
$f'_{[0,\infty)}\:P'_{[0,\infty)}\to Q'_{[0,\infty)}$ represents the Steenrod
homotopy class of $f$.
On the other hand, each $Q'_i$ obviously collapses onto $P'_i$.

The strong deformation retraction of $Q'_i$ onto $P'_i$ yields a deformation
retraction of $Q'_i\x I$ onto $P'_i\x I\cup Q'_i\x\partial I$.
Similarly, $Q'_{[i,i+1]}$ deformation retracts onto
$P'_{[i,i+1]}\cup Q_{i+1}'\cup Q'_i$.
Composing the latter retractions with the retractions $Q'_i\to P'_i$,
we obtain a proper deformation retraction
$r\:Q'_{[0,\infty)}\to P'_{[0,\infty)}\cup Q'_{\N}\to P'_{[0,\infty)}$.
Thus $r$ is a proper homotopy equivalence inverting the proper homotopy
class of the inclusion $f'_{[0,\infty)}$. \qed
\enddemo

A map $f\:X\to Y$ (respectively, a Steenrod homotopy class $[F]\:X\ssm Y$)
between compacta with basepoints $x\in X$, $y\in Y$ is called {\it pointed}
if it is basepoint preserving, i.e.\ $f(x)=y$ (respectively, if the proper map
$F\:P_{[0,\infty)}\to Q_{[0,\infty)}$ is base ray preserving, i.e.\
$F(p_{[0,\infty)})\i q_{[0,\infty)}$, where $p_i$ and $q_i$ are the images of
$x$ and $y$ in $P_i$ and $Q_i$).
The set $[(X,x),\,(Y,y)]$ of pointed Steenrod homotopy classes receives
the natural structure of a pointed set (i.e.\ a set with a distinguished point).

\proclaim{Proposition 2.9}
Let $X$ and $Y$ be connected compacta with basepoints $x$ and $y$, and let
$f\:(X,x)\ssm(Y,y)$ be a Steenrod homotopy class.
If the underlying unpointed class $\bar f\:X\ssm Y$ of $f$ is a Steenrod
homotopy equivalence then $f$ is a Steenrod homotopy equivalence.
\endproclaim

This was proved in \cite{DS2; 4.6} and \cite{C1; 1.16 (see also 1.14, 1.15)}.

\demo{Proof}
By the pointed version of Proposition 2.8(a), we may assume without loss of
generality that $f$ is an inclusion map (then, in particular, $x=y$).
Then by (i)\imp(ii) in Proposition 2.8(b), $X$ is a Steenrod strong
deformation retract of $Y$.
But this implies that $(X,x)$ is a Steenrod strong deformation
retract of $(Y,y)$. \qed
\enddemo

A compactum $X$ is {\it Steenrod connected} if there is only one Steenrod
homotopy class from a point into $X$.
Proposition 2.9 implies that there is no difference between pointed and
unpointed shape for Steenrod connected compacta.

\example{Example 2.10 (topologist's sine curve)} Let $X\i [0,1]\x[-1,1]$ be
the closure of the graph of $\sin\frac1x$, where $x$ runs over $(0,1]$.
Then $X$ is Steenrod connected (in fact, even cell-like), but not path connected.
\endexample

Steenrod connected compacta will be considered in more detail in \S8.
For instance, Theorem 8.7 (due originally to Krasinkiewicz and Minc) implies
that a compactum $X$ is Steenrod connected if and only if every two maps
$pt\to X$ represent the same Steenrod homotopy class.
In particular, path connected compacta are Steenrod connected.

On the other hand, as observed in \cite{GK; Remark 10.1}, the set of Steenrod
homotopy classes $pt\ssm X$ that are representable by maps is not a shape
invariant: if $\alpha\:pt\ssm X$ is not representable by a map, represent it by
a level-preserving map $f\:[0,\infty)\to P_{[0,\infty)}$, where $X=\invlim P_i$.
Then $Y:=X\cup f([0,\infty))$ is Steenrod homotopy equivalent to $X$, but
$\alpha\:pt\ssm Y$ is now represented by $f|_{\{0\}}$.

\head 3. Homotopy groups \endhead

The definition of the set $[X,Y]$ of Steenrod homotopy classes $X\ssm Y$
generalizes to pairs of compacta in the obvious way, based on the fact that
a pair of compacta is an inverse limit of pairs of compact polyhedra.
We define the {\it Steenrod homotopy set} $\pi_n(X;\,x)=[(S^n,pt),\,(X,x)]$.
Thus (see Proposition 2.4(a)), $\pi_n(X;\,x)$ can be identified with
the set of level-preserving homotopy classes of level-preserving maps
$$\left(S^n\x[0,\infty),\,pt\x[0,\infty)\vphantom{_[}\right)\to
\left(P_{[0,\infty)},\,p_{[0,\infty)}\right)$$
where $P=(\dots@>p_1>>P_1@>p_0>>P_0)$ is an inverse sequence of compact
polyhedra and PL maps with inverse limit $X$; and $p_{[0,\infty)}$ is
the mapping telescope of the images $p_i$ of the basepoint $x$ in the $P_i$'s.
By Proposition 2.3(a), if $X$ is a polyhedron,
$\stau\:\spi_n(X;\,x)\to\pi_n(X;\,x)$ is a bijection.

The group operations (for $n\ge 1$) and the induced maps $f_*$, where $f$
is a pointed (i.e.\ basepoint preserving) map or a pointed (i.e.\ base ray
preserving) Steenrod homotopy class, are defined in the obvious way.
We suppress the basepoint from notation and write $\pi_n(X)$ when it does not
lead to a confusion.
Next, let $\pi_n(X,A;\,x)=[(D^n,\partial D^n,pt),\,(X,A,x)]$ for any pair of
compacta $(X,A)$ with a basepoint $x\in A$ and any $n\ge 1$.
More generally, we allow any Steenrod homotopy class $x\:pt\ssm A$ as a basepoint.
We also write $\pi_n(X,A)$ when the basepoint is not explicitly used.
The group operations on $\pi_n(X,A)$ are defined in the obvious way for $n\ge 2$.
There is an action of $\pi_1(A)$ on $\pi_n(X)$ for $n\ge 1$ and on $\pi_n(X,A)$
for $n\ge 2$, and an action of $\pi_1(X)$ on $\pi_1(X,A)$.
We leave $\pi_0(X,A)$ undefined.

\definition{Topology on $\pi_n(X)$}
A base of the topology on $\pi_n(X)$ is given by the point-inverses of
the maps $f_*\:\pi_n(X)\to\pi_n(P)$ induced by all maps $f\:X\to P$, where $P$
is a compact polyhedron.
Clearly, group operations (when $n>0$) as well as all induced maps
are continuous in this topology.
If $X$ is the inverse limit of polyhedra $P_i$, by Lemma 2.1(a) every map $f$
from $X$ to a polyhedron $Q$ extends to a map $P_{[k,\infty]}\to Q$ for some $k$,
and so $f_*$ factors through $(p^\infty_k)_*$.
Consequently, a base of the topology on $\pi_n(X)$ can also be taken to
consist of the point-inverses of all maps $(p^\infty_i)_*$ induced by
the projections $p^\infty_i\:X\to P_i$.

Clearly, the topological group (pointed space when $n=0$) $\pi_n(X,x)$ is
a shape invariant of the pair $(X,x)$.
When $X$ is Steenrod connected (i.e.\ $\pi_0(X)=pt$), it is
also a shape invariant of $X$ by Proposition 2.9.
In addition, the underlying unpointed space of $\pi_0(X)$ is obviously
a shape invariant of $X$.
\enddefinition

\definition{Derived limit}
Let $G=(\dots@>p_1>>G_1@>p_0>>G_0)$ be an inverse sequence of groups.
The group $\prod G_i$ acts on itself on the right by
$$(g_0,\,g_1,\dots)\cdot(x_0,x_1,\dots)=
(x_0^{-1}g_0p_0(x_1),\,x_1^{-1}g_1p_1(x_2),\dots).$$
The stabilizer of $1$ under this action can be identified with the inverse limit
$\invlim G_i$.
We define $\derlim G_i$ to be the orbit space of this action, thought of as
a pointed set.
Note that the map $f\:\prod G_i\to\prod G_i$ of pointed sets, defined by
$f(g_0,g_1,\dots)=(g_0^{-1}p_0(g_1),g_1^{-1}p_1(g_2),\dots)$, is equivariant
with respect to the right regular action of $\prod G_i$ on the domain and the
above-defined action on the target.
We have $\invlim G_i=\ker f$, and $\derlim G_i$ vanishes iff $f$ is surjective.
When the $G_i$ are abelian, $f$ is a homomorphism, moreover, $\derlim G_i$ can be
identified with $\coker f$ and thus becomes endowed with an abelian group
structure.%
\footnote{If for each $i$ and each $x_i\in G_i$ there exists a $j>i$ such that
the image of every $y_j\in G_j$ in $G_i$ commutes with $x_i$, then there is
a monoid structure on $\derlim G_i$; and if additionally each $G_i$ is solvable,
it extends to a group structure \cite{G\"u}.}

We endow $\invlim G_i$ with the inverse limit topology, which is induced from
the product topology on the product $\prod G_i$ of the discrete groups;
and $\derlim G_i$ with the quotient topology, which is always anti-discrete
since changing finitely many components of any $(g_1,g_2,\dots)\in\prod G_i$
does not change its orbit.

When the groups $G_i$ are abelian, $G\mapsto\derlim G_i$ is the derived functor
of the inverse limit functor $G\mapsto\invlim G_i$ in the sense of
homological algebra.
For it can be shown that $\derlim G_i$ vanishes if all bonging maps in $G$
are epimorphisms, and therefore
$$\CD
\dots@>>>G_2            @>>>G_1      @>>>G_0\\
@.@V{\left(\smallmatrix 1\\ p_1\\ p_0p_1\endsmallmatrix\right)}VV
@V{\left(\smallmatrix 1\\ p_0\endsmallmatrix\right)}VV@V1VV\\
\dots@>>>G_2\x G_1\x G_0@>>>G_1\x G_0@>>>G_0\\
@.@V{\left(\smallmatrix
p_1\, & -1        & \,0\\
0\,     & \,p_0\, & \,-1
\endsmallmatrix\right)}
VV@V{\left(\smallmatrix p_0\, & \,-1\endsmallmatrix\right)}VV@VVV\\
\dots@>>>G_1\x G_0      @>>>G_0      @>>>0\\
\endCD$$
is a $\invlim$-acyclic resolution of $G$ (it is due to J. E. Roos).
In the abelian case, Theorem 3.1(d) below is the usual long exact sequence
for derived functors, and Theorem 4.1(iii) below (the homology version of
Theorem 3.1(b)) is the universal coefficients formula for derived functors.
\enddefinition

\proclaim{Theorem 3.1}
(a) \cite{Q2} If $(X,A)$ is a pair of compacta with a basepoint $x\in A$ (or
more generally $x\:pt\ssm A$), there is an exact sequence of pointed sets
$$\dots\to\pi_2(X,A)\to\pi_1(A)\to\pi_1(X)\to\pi_1(X,A)\to\pi_0(A)\to\pi_0(X),$$
whose maps to the left of $\pi_1(X)$ (resp.\ $\pi_2(X)$) are homomorphisms
of groups (resp.\ right $\Z\pi_1(A)$-modules).
In addition, $\pi_2(X)\to\pi_2(X,A)$ is also $\pi_1(A)$-equivariant,
$\partial\:\pi_2(X,A)\to\pi_1(A)$ is a crossed module%
\footnote{That is, $\partial(s\cdot g)=g^{-1}(\partial s)g$ for $g\in\pi_1(A)$,
$s\in\pi_2(X,A)$ and $s^{-1}ts=t\cdot(\partial s)$ for $s,t\in\pi_2(X,A)$.}%
, $\pi_1(X)\to\pi_1(X,A)$ is $\pi_1(X)$-equivariant with respect to the right
regular action on $\pi_1(X)$, and the non-trivial point-inverses of
$\pi_1(X,A)\to\pi_0(A)$ coincide with the orbits of $\pi_1(X)$.

\medskip
(b) \cite{Q1}, \cite{Q2}, \cite{Gr}, \cite{EH; 5.2.1}, \cite{Wa}
If $(X,A)$ is the inverse limit of pairs $(P_i,Q_i)$ of connected compact
polyhedra, there is an exact sequence
$$1\to\derlim\pi_{n+1}(P_i,Q_i)\to\pi_n(X,A)\to\invlim\pi_n(P_i,Q_i)\to 1$$
of pointed sets when $n\ge 1$ or $A=Q_i=pt$ and $n=0$, whose maps are group
homomorphisms when $n\ge 2$ or $A=Q_i=pt$ and $n=1$.

\medskip
(c) The conclusion of (b) holds when $(P_i,Q_i)$ are pairs of connected compacta.

\medskip
(d) \cite{BK; IX.2.3 (see also Ch.\ XI)}, \cite{MaS; p.\ 168} If
$1\to K_i@>j_i>>G_i@>f_i>>Q_i\to 1$ is an exact sequence of inverse
sequences of groups (i.e.\ $j_i$ and $f_i$ commute with the bonding maps),
there is an exact sequence of pointed sets
$$1\to\invlim K_i@>\!\lim\! j_i\!>>\invlim G_i@>\!\lim\! f_i>>\invlim Q_i@>\delta>>
\derlim K_i@>\!\lim^1\!\! j_i>>\derlim G_i@>\!\lim^1\!\! f_i>>\derlim Q_i\to 1,$$
where $\delta$ is $\invlim Q_i$-equivariant with respect to the right regular
action on $\invlim Q_i$ and the right action
$$[(k_0,k_1,\dots)]\cdot(g_0K_0,g_1K_1,\dots)=
[(g_0^{-1}k_0p_0(g_1),\,g_1^{-1}k_1p_1(g_2),\dots)]$$
on $\derlim K_i$, and the non-trivial point-inverses of
$\derlim j_i$ coincide with the orbits of $\invlim Q_i$.

\medskip
(d$'$) \cite{Ir} If $1\to H_i@>j_i>>G_i@>f_i>>G_i/H_i\to 1$ is an exact sequence
of inverse sequences of pointed sets, where $G_i$ and $H_i$ are groups and
$j_i$ are homomorphisms, there is an exact sequence of pointed sets
$$1\to\invlim H_i@>\lim j_i>>\invlim G_i@>\lim f_i>>\invlim (G_i/H_i)@>\delta>>
\derlim H_i@>\lim^1 j_i>>\derlim G_i,$$
where the non-trivial point-inverses of $\delta$ coincide with the orbits of
the obvious action of $\invlim G_i$ on $\invlim (G_i/H_i)$.
\endproclaim

The proof shows that all homomorphisms in Theorem 3.1 are natural and continuous.

\demo{Proof. (a)} This follows from the homotopy exact sequence for proper
maps, which is verified directly. \qed
\enddemo

\demo{(b)} A direct verification shows that the following sequence of
pointed sets of level-preserving pointed homotopy classes is exact to
the left of the $a$ arrow with $n=1$ or with $Q_i=pt$ and $n=0$:
$$\multline\dots\to
[(D^{n+1},S^n)\x\N,\,(P_\N,Q_\N)]^*_\ell@>b>>
[(D^n,S^{n-1})\x[0,\infty),\,(P_{[0,\infty)},\,Q_{[0,\infty)})]^*_\ell\\
@>c>>[(D^n,S^{n-1})\x\N,\,(P_\N,Q_\N)]^*_\ell@>a>>\dots
\endmultline$$
``Pointed'' means that $pt\x[0,\infty)$ always maps to $p_{[0,\infty)}$ and
$pt\x\N$ always maps to $p_\N$.
Here $c$ is defined by restriction, and $b$ using the obvious identification
of $[(D^{n+1},S^n)\x\N,\,(P_\N,Q_\N)]^*_\ell$ with
$$[(D^n\x [0,\infty),\,S^{n-1}\x[0,\infty)\cup
D^n\x\N),\,(P_{[0,\infty)},Q_{[0,\infty)})]^*_\ell$$
for $n>0$ or when $n=0$ and $Q_i=pt$.

Finally, $a$ may be identified with $f\:\prod\pi_n(P_i,Q_i)\to
\prod\pi_n(P_i,Q_i)$, given by
$(g_1,g_2,\dots)\mapsto(p_{1*}(g_2)-g_1,\,p_{2*}(g_3)-g_2,\dots)$ in
the abelian case (for $n\ge 3$ or $Q_i=pt$ and $n=2$) and
by $(g_1,g_2,\dots)\mapsto(g_1^{-1}p_{1*}(g_2),\,g_2^{-1}p_{2*}(g_3),\dots)$
in the group case (for $n\ge 2$ or $Q_i=pt$ and $n=1$).
In the pointed set case, $a$ can be similarly defined
if its target $[(D^1,S^0)\x\N,\,(P_\N,Q_\N)]^*_\ell$ for $n=1$ (respectively
$[S^0\x\N,\,P_\N]^*_\ell$ for $Q_i=pt$ and $n=0$) is amended by moving
the basepoint into $D^1\but S^0$, contrary to the standard convention
(respectively, is replaced with the basepoint free homotopy set
$[S^0\x\N,\,P_\N]_\ell$).

The maps in the long exact sequence are group homomorphisms to the left of
$a\:\prod\pi_n(P_i,Q_i)\to\prod\pi_n(P_i,Q_i)$ (non-inclusively!), where $n=2$,
or $Q_i=pt$ and $n=1$.
The assertion now follows since $\invlim\pi_n(P_i,Q_i)=\ker f$ and in
the abelian case also $\derlim\pi_{n+1}(P_i,Q_i)=\coker f$.
In the non-abelian case, writing the rightmost homomorphism $a$ as $a\:G\to S$,
we note that $G$ acts on $S$ by the formula in the definition of $\derlim$,
and $a$ is equivariant with respect to this action on $S$ and the right regular
action on $G$.
The orbits of this action are precisely the non-trivial point-inverses of
the subsequent (to the right) map $b$. \qed
\enddemo

\demo{(c)} We suppress the basepoints.
For the sake of clarity, the inverse sequence of compacta
$\dots@>p_1>>P_1@>p_0>>P_0$ will be denoted by $\dots@>q_1>>C_1@>q_0>>C_0$.

Let us represent each $C_i$ as the inverse limit of
$P^{(i)}=(\dots@>p^{(i)}_1>>P^{(i)}_1@>p^{(i)}_0>>P^{(i)}_0)$, where
$P^{(i)}_j$ are compact polyhedra and $p^{(i)}_j$ are PL maps.
The maps $q^i$ induce Steenrod homotopy classes represented by proper maps
$q^i_{[0,\infty)}\:P^{(i+1)}_{[0,\infty)}\to P^{(i)}_{[0,\infty)}$.
By Lemma 2.5(b$_0$), we may assume after passing to a subsequence of $P^{(1)}$
that $q^0_{[0,\infty)}$ is level-preserving.
Continuing this procedure, we may assume that all $q^i_{[0,\infty)}$ are
level-preserving.
Let $T_J=(\dots@>q^1_J>>P^{(1)}_J@>q^0_J>>P^{(i)}_J)$, and let us write
$T_{J\x K}$ for $(T_J)_K$.
By the proof of (b), the exact sequence of (c) holds provided that
$\pi_n(X)$ is replaced by the group $G$ of fibered homotopy classes of maps
$F\:S^n\x[0,\infty)\x [0,\infty)\to T_{[0,\infty)\x [0,\infty)}$ that are
fibered over $[0,\infty)\x[0,\infty)$.

Consider the diagonal ``staircase'' mapping telescope
$$T^\Delta=T_{[0,1]\x 0}\cup T_{1\x [0,1]}\cup T_{[1,2]\x 1}\cup T_{2\x [1,2]}
\cup\dots =Q_{[0,\infty)},$$ where $Q_{2i}=P^{(i)}_i$ and
$Q_{2i+1}=P^{(i)}_{i+1}$.
Clearly, $X$ is homeomorphic to $\invlim Q_i$.
On the other hand, restricting the fiberwise map $F$ over $T^\Delta$ yields
a level-preserving $F^\Delta\:S^n\x[0,\infty)\to Q_{[0,\infty)}$, and
conversely, every $F^\Delta$ extends to an $F$ using the homotopies similar
to $\Pi_t$.
It follows that $G\simeq\pi_n(X)$. \qed
\enddemo

\demo{(d)} For variety, we include a geometric proof in the case where each $G_i$
has a classifying space with compact $3$-skeleton $BG_i$, and each $Q_i$ has
a classifying space with compact $4$-skeleton $BQ_i$; the general countable group
case can be done similarly using the definitions at the end of \S7.
By considering the mapping cylinder, we may assume that $BG_i\i BQ_i$, which
inclusion induces the given surjection $G_i\to Q_i$ on $\pi_1$.
We have $\pi_2(BG_i)=0=\pi_2(BQ_i)$ and $\pi_3(BQ_i)=0$, implying
$\pi_3(BQ_i,BG_i)=0$.
Let $\dots\to(BG_1,BQ_1)\to(BG_0,BQ_0)$ be an inverse sequence whose bonding
maps induce the given homomorphisms $G_{i+1}\to G_i$ (and consequently also
$Q_{i+1}\to Q_i$).
Let $(X,A)=\invlim (BQ_i,BG_i)$.
Then part (b) identifies the terms from the exact sequence of (a) with those
from the statement of (d).
The surjectivity of the rightmost map follows since in the notation of (b),
if $\pi_0(P_i)$, $\pi_0(Q_i)$ and $\pi_1(P_i,Q_i)$ vanish, then
$\pi_0(A)\to\pi_0(X)$ is surjective by the proof of (b). \qed
\enddemo

\demo{(d$'$)} Similar to (d); see also the proof of Lemma 3.7(a) below. \qed
\enddemo

\example{Example 3.2}
As a sample application of Theorem 3.1(d), let us compute the group
$\derlim(\dots@>p>>\Z@>p>>\Z)$.
The given inverse sequence fits into a short exact sequence with
$G=(\dots@>\id>>\Z@>\id>>\Z)$ and $Q=(\dots@>\bmod p^2>>\Z/p^2@>\bmod
p>>\Z/p@>>>1)$.
Obviously $\derlim G$ is trivial, so
$\derlim(\dots@>p>>\Z@>p>>\Z)\simeq\invlim Q/\invlim G=\Z_p/\Z$.
\endexample

The similarity between this calculation and the proof of Lemma 3.3 below
will be analyzed and applied in the proof of Theorem 8.4.

An inverse sequence $\dots\to G_1\to G_0$ satisfies the {\it Mittag-Leffler
condition}, if for each $i$ there exists a $j>i$ such that for each $k>j$
the image of $G_k\to G_i$ equals that of $G_j\to G_i$.
When the $G_i$ are groups, it is an easy but enlightening exercise to verify that
the Mittag-Leffler condition implies the vanishing of $\derlim G_i$.
The converse holds if each $G_i$ is countable:

\proclaim{Lemma 3.3 \cite{G}, \cite{Ge1}, \cite{DS1}, \cite{MM}} Let $G_i$ be
an inverse sequence of countable groups.
If $\derlim G_i$ is countable (for example, trivial), then $G_i$ satisfy
the Mittag-Leffler condition.
\endproclaim

We include a proof since it will be needed later on several occasions.
Our presentation is closest to \cite{G}, where Lemma 3.3 is stated under
an additional hypothesis.

\demo{Proof} Suppose that the $G_i$ are not Mittag-Leffler.
Then there exists a $k$ such that the images $A_i$ of $G_{k+i}$ in $G_k$ do not
stabilize.
Since each $G_i\to A_i$ is onto, $\derlim G_i\to\derlim A_i$ is onto, so
$\derlim A_i$ is countable.
If $(g_0,g_1,\dots)$ and $(h_0,h_1,\dots)$ are elements of $\prod A_i$
representing the same element of $\derlim A_i$, each $g_i=x_ih_ix_{i+1}^{-1}$
for some $(x_0,x_1,\dots)\in\prod A_i$.
Then each $g_0\dots g_n=x_0h_0\dots h_nx_{n+1}^{-1}$.
Hence $(g_0A_1,\,g_0g_1A_2,\dots)$ and $(h_0A_1,\,h_0h_1A_2,\dots)$ lie in
the same orbit of the left action of $A_0$ on $A:=\invlim A_0/A_i$.
As this argument is reversible, we obtain a bijection between $\derlim A_i$
and the orbit set $A_0\q A$.
However, $A$ admits a bijection with the uncountable set
$(A_0/A_1)\x (A_1/A_2)\x\dots$, hence $A_0\q A$ must be
uncountable. \qed
\enddemo

An inverse sequence of pointed sets $G_i$ is said to satisfy the {\it dual
Mittag-Leffler condition}, if there exists a $k$ such that for every $j>k$
there exists an $i>j$ for which the kernel from $G_i$ to $G_k$ equals that
to $G_j$ (cf.\ \cite{Ge2}).
If this condition holds, it is easy to see that $\invlim G_i$ injects into
some $G_k$.

\proclaim{Lemma 3.4} Let $\dots\to G_1\to G_0$ be an inverse sequence of pointed
sets. For (b) and (c), assume additionally that either it is an inverse sequence
of groups or each $G_i$ is finite.

(a) If the $G_i$ satisfy the Mittag-Leffler condition and $\invlim G_i=0$, then
for each $i$ there exists a $j>i$ such that $G_j\to G_i$ is trivial.

(b) If $\invlim G_i$ is discrete, it injects into some $G_k$.

(c) Let the $G_i$ satisfy the Mittag-Leffler condition.
If $\invlim G_i$ is discrete, then the $G_i$ satisfy the dual Mittag-Leffler
condition.
\endproclaim

The converse statements to (a), (b) and (c) are true, and obvious.
The additional hypothesis for (b) and (c) can be dropped if $\invlim G_i$,
understood as an inverse limit of discrete uniform spaces, is assumed to be
discrete as a uniform space.

Part (b) will be needed a little later (in Theorem 3.10 and also in \S6).
Part (c) will be used in Theorem 3.12.

\demo{Proof. (a)}
Let $L_i$ be the stable image of the $G_j$ with $j>i$ in $G_i$.
Then $L_{i+1}$ maps onto $L_i$, and consequently $\invlim L_j$ maps onto $L_i$.
The composition $\invlim L_j\to L_i\i G_i$ factors through $\invlim G_i$,
hence $L_i$ is trivial. \qed
\enddemo

\demo{(b)} If $g_1,g_2,\dots$ is a sequence of elements of $\invlim G_i$ such
that each $g_i$ maps trivially to $G_i$, then this sequence converges to
the identity in the inverse limit topology.
It remains to consider the case of finite sets.
Suppose that $g_1,h_1,g_2,h_2,\dots$ is a sequence of elements of $\invlim G_i$
such that $g_i$ and $h_i$ map to the same element of $G_i$.
Since $\invlim G_i$ is a compactum, there exists an infinite increasing
sequence $n_i$ such that $g_{n_i}$ and $h_{n_i}$ converge to some $g$ and $h$
respectively.
Then $g=h$, and since $\invlim G_i$ is discrete, $g_{n_i}=g$ and $h=h_{n_i}$
for sufficiently large $i$. \qed
\enddemo

\demo{(c)} Let $L_i$ be the stable image of the $G_j$ with $j>i$ in $G_i$.
Then $\invlim G_j$ maps onto $L_i$.
By the hypothesis, there exists a $k$ such that $\invlim G_j\to L_i$
is a bijection for all $i\ge k$.
Given a $j>k$, let $i>j$ be such that $G_i$ maps into $L_j$.
Since $L_j\to L_k$ is a bijection, $\ker(G_i\to G_j)=\ker(G_i\to G_k)$.
\qed
\enddemo

\proclaim{Proposition 3.5} Let $X$ be a compactum, and fix some $m\ge 0$.

(a) $X$ is the limit of an inverse sequence of $m$-connected polyhedra iff
$\pi_n(X)=0$ for all $n<m$ and $\pi_m(X)$ is anti-discrete.

(b) $X$ is the limit of an inverse sequence of $(m+1)$-connected maps between
$m$-connected polyhedra iff $\pi_n(X)=0$ for all $n\le m$.
\endproclaim

We recall that a map $f\:X\to Y$ is called {\it $m$-connected} if its
mapping cylinder $MC(f)$ is relatively $m$-connected, i.e.\ $\pi_i(MC(f),X)=0$
for $i\le m$.
By the homotopy exact sequence of a pair, a map $f\:X\to Y$ between
$m$-connected polyhedra is $(m+1)$-connected iff $\pi_{m+1}(X)\to\pi_{m+1}(Y)$
is a surjection.

Compacta that satisfy any of the equivalent conditions of Proposition 3.5(a)
are called {\it UV$_m$} or {\it approximatively $m$-connected} or
{\it $m$-shape connected} in the literature.
This notion is already found in \cite{Ch; \S5}, along with the ``only if''
directions of both parts of Proposition 3.5.
Part (a) is well-known (see e.g.\ \cite{DS1; proof of 8.3.2}).

\demo{Proof. (a)} The ``only if'' part is obvious.
Conversely, suppose that $X$ is the inverse limit of polyhedra $P_i$.
By Theorem 3.1(b), $\invlim\pi_n(P_i)=0=\derlim\pi_n(P_i)$ for $n\le m$.
By Lemmas 3.3 and 3.4(a), for each $i$ and each $n\le m$ there exists a
$j=j_n(i)$ such that $\pi_n(P_j)\to\pi_n(P_i)$ is trivial.
If $j=j_0(j_1(\dots j_k(P_i)))$, then, by induction on $k$, the restriction
of $P_j\to P_i$ to the $k$-skeleton $P_j^{(k)}$ of some fixed triangulation
of $P_j$ is null-homotopic.
Set $i_0=0$ and $i_{k+1}=j_0(j_1(\dots j_m(i_k)))$ and let
$Q_k$ be the union of $P_{i_k}$ with the cone $CP_{i_k}^{(m)}$.
Then the bonding map $p^{i_{k+1}}_{i_k}\:P_{i_{(k+1)}}\to P_{i_k}$ extends to
a map $Q_{k+1}\to P_{i_k}$.
Then $X$ is the inverse limit of the maps $Q_{k+1}\to P_{i_k}\i Q_k$. \qed
\enddemo

\demo{(b)} The ``only if'' part follows from Theorem 3.1(b).
Conversely, suppose that $X$ is the inverse limit of $m$-connected polyhedra
$P_i$ (using (a)) and that $G_i:=\pi_{m+1}(P_i)$ is Mittag-Leffler (using
Lemma 3.3).
Thus for each $i$ there exists a $j(i)>i$ such that for each $k>j$ the image
of $p^k_i\:G_k\to G_i$ equals that of $p^j_i$.
Now $G_j$ is finitely generated, using the Hurewicz theorem when $m>0$.
Let $g_1,\dots,g_r$ be a set of generators.
Then $p^j_i(g_l)=p^{j+1}_i(h_l)$ for some $h_1,\dots,h_r\in G_{j+1}$.
Let $f_l\:(S^{m+1},pt)\to (P_j,p_j)$ be a spheroid representing
$g_l'=g_l^{-1}p^{i+1}_i(h_l)$, and let $Q_j$ be obtained by gluing up
$f_1,\dots,f_r$ by $(m+2)$-cells.
Then $P_j\to P_i$ extends to a map $Q_j\to P_i$, and the composition
$P_{j+1}\to P_j\i Q_j$ induces an epimorphism on $\pi_{m+1}$.
Setting $i_0=0$ and $i_{n+1}=j(i_n)$, we have $X=\invlim Q_{i_n}$, where
each $Q_{i_n}$ is $m$-connected and each bonding map $Q_{i_{n+1}}\to Q_{i_n}$
is $(m+1)$-connected. \qed
\enddemo

\remark{Remark} By a well-known example depending on a deep result of J. F. Adams,
Proposition 3.5 does not generalize to the case $m=\infty$, i.e.\ there exists
an infinite-dimensional non-cell-like compactum $X$ with $\pi_n(X)=0$
for all $n$ (see \cite{EH; 5.5.10}, \cite{DS1; 10.3.1}).
\endremark
\medskip

The following ``Whitehead theorem in Steenrod homotopy'' answers a question of
Koyama.

\proclaim{Theorem 3.6} Let $X$, $Y$ be connected compacta of dimensions $\le m$
and $\le m+1$ respectively, and let us consider a Steenrod
homotopy class $f\:(X,x)\ssm (Y,y)$.
If $f_*\:\pi_n(X)\to\pi_n(Y)$ is a bijection for $n\le m$ and
a surjection for $n=m+1$, then $f$ is a Steenrod homotopy equivalence.
\endproclaim

Just as with the classical Whitehead Theorem, the proof of Theorem 3.6
works to yield slightly more general assertions (a) and (b) below;
addendum (c) follows from Theorem 3.6 itself along with the last assertion of
addendum (b).

\proclaim{Addenda} Let $X$ and $Y$ be connected compacta and
$f\:(X,x)\ssm (Y,y)$ a Steenrod homotopy class.

(a) If  $\dim X\le m$, $\dim Y\le m+1$, and $i\:Y\emb Z$ is an inclusion into
a connected compactum $Z$ such that the composition
$\pi_n(X)@>f_*>>\pi_n(Y)@>i_*>>\pi_n(Z)$ is a bijection for $n\le m$ and
a surjection for $n=m+1$, then there exists a $g\:Y\ssm X$ such that
$gf=[\id_X]$ and $[i]fg=[i]$.
In particular, $f_*\:\pi_n(X)\to\pi_n(Y)$ is a split injection for all $n$.

(b) If $\dim Y\le m+1$, and $f_*\:\pi_n(X)\to\pi_n(Y)$ is a bijection for
$n\le m$ and a surjection for $n=m+1$, then there exists a $g\:Y\ssm X$ such
that $fg=[\id_Y]$; and if $i\:Z\emb X$ is an inclusion of a connected
compactum $Z$ of dimension $\le m$, then $g$ can be chosen so that additionally
$gf[i]=[i]$.
In particular, $f_*\:\pi_n(X)\to\pi_n(Y)$ is a split surjection for all $n$.

(c) If $\dim X\le m$ and $\dim Y\le m$, and $f_*\:\pi_n(X)\to\pi_n(Y)$ is
a bijection for all $n\le m$, then $f$ is a Steenrod homotopy equivalence.
\endproclaim

A pro-group version of Theorem 3.6 appears in \cite{EH; 5.5.6} with a sketch
of proof; a different proof is given in \cite{DS2} (see also \cite{DS1} and
\cite{AM; \S4}).
Koyama used it to prove Addendum 3.6(c) under two additional hypotheses:
(i) $f_*\:\check\pi_n(X)\to\check\pi_n(Y)$ is a bijection for $n\le m$;
(ii) a technical condition which is slightly weaker than requiring $Y$ to be
an inverse limit of polyhedra with abelian fundamental groups \cite{Koy}.
He then asked whether these conditions (i) and (ii) are superfluous
\cite{Koy; Problem 2}.
A hopelessly erroneous proof of Addendum 3.6(c) appears in \cite{Li}, where
(without any comments) $\pi_1(P,Q)$ is thought to be a group for any polyhedral
pair $(P,Q)$ and then the vanishing of ``$\derlim$'' of an inverse sequence
of such ``groups'' is taken to imply the Mittag-Leffler condition.
The correct part of the arguments found in \cite{Li} amounts to a proof that
condition (i) can be dropped from the hypothesis of Koyama's theorem, but
represents no progress with respect to the question of necessity of (ii).

\demo{Proof of Theorem 3.6} Let
$f_{[0,\infty)}\:P_{[0,\infty)}\to Q_{[0,\infty)}$ represent $f$.
Without loss of generality all the $P_i$ and $Q_i$ are connected.
By Lemma 2.5(b$_0$), after passing to a subsequence of $P_i$'s we may assume
that $f_{[0,\infty)}$ is level-preserving.
Given any $J\i [0,\infty)$, we write $f_J$ for the restriction
$P_J\to Q_J$ of $f$ and $M_J$ for the mapping cylinder of $f_J$.
We may assume that $M_{[0,\infty)}$ is of dimension at most $m+2$.
From the hypothesis and the proof of Theorem 3.1(a) we get that
$\pi_n(M_{[0,\infty)},\,P_{[0,\infty)})=0$ for $1\le n\le m+1$.
By the proof of Theorem 3.1(b) we have short exact sequences
$$0\to\derlim\pi_{n+1}(M_i,\,P_i)\to\pi_n(M_{[0,\infty)},\,P_{[0,\infty)})
\to\invlim\pi_n(M_i,\,P_i)\to 0$$
for $n\ge 1$ and also that $\derlim\pi_1(P_i)\to\derlim\pi_1(M_i)$ is onto.
Then $\pi_n(M_i,\,P_i)$ have trivial inverse limit for $1\le n\le m+1$ and
trivial derived limit for $2\le n\le m+2$.
By Lemma 3.3 and Lemma 3.7(b) below, the $\pi_n(M_i,\,P_i)$ are Mittag-Leffler
for each $n$ with $1\le n\le m+2$.

By Lemma 3.4(a), for each $i$ there exists a $j$ such that
$\pi_n(M_j,\,P_j)\to\pi_n(M_i,\,P_i)$ is trivial for $1\le n\le m+1$.
By passing to a subsequence, we may assume that
$\pi_n(M_{i+1},\,P_{i+1})\to\pi_n(M_i,\,P_i)$ is trivial for $1\le n\le m+1$.
We may also assume that the image of $\pi_{m+2}(M_{i+1},\,P_{i+1})$ in
$\pi_{m+2}(M_i,\,P_i)$ equals that of $\pi_{m+2}(M_j,\,P_j)$ for all $j>i$.
Set $i_k=(2m+3)k$ and let us fix some triangulation of each $(M_{i_k},P_{i_k})$.
By an induction on $j=0,1,\dots,m+1$, the bonding map
$p\:M_{i_k}\to M_{i_k-j}$ is homotopic rel $P_{i_k}$ to a map
$\phi_j$ sending the $j$-skeleton $M_{i_k}^{(j)}$ into $P_{i_k-j}$.
Thus, writing $i_k'=i_k-m-1$, the resulting map $\psi_k:=\phi_{m+1}$ sends
$M_{i_k}$ into $P_{i_k'}$, restricts to the bonding map
$p\:P_{i_k}\to P_{i_k'}$, and is homotopic to the bonding map
$p\:M_{i_k}\to M_{i_k'}$ by a homotopy $\Psi_k$.

Then $\mu\:M_{i_{k+1}}@>p>>M_{i_k}@>\psi_k>>P_{i_k'}$ and
$\nu\:M_{i_{k+1}}@>\psi_{k+1}>>P_{i_{k+1}'}@>p>>P_{i_k'}$
are homotopic with values in $M_{i_k'}$ by a homotopy $h_k$.
Writing $i_k''=i'_k-m-1=i_{k-1}+1$, similarly to the above, the compositions
of $\mu$ and $\nu$ with the bonding map $P_{i_k'}\to P_{i_k''}$ are homotopic
by a homotopy $h'_k\:M_{i_{k+1}}\x I\to M_{i_k''}$ whose restriction to
$M_{i_k}^{(m)}\x I$ has values in $P_{i_k''}$.
For each $(m+1)$-simplex $\Delta$ of $M_{i_{k+1}}$, we have
$h'_k(\partial(\Delta\x I))\i P_{i_k''}$.
Let $\alpha_\Delta\:(\Delta\x I,\partial)\to (M_{i_{k+1}'},\,P_{i_{k+1}'})$
be a relative spheroid whose composition with the bonding map
$M_{i_{k+1}}\to M_{i_{k-1}}$ is homotopic in $(M_{i_{k-1}},\,P_{i_{k-1}})$
to the composition of $h'_k|_{\Delta\x I}$ and the bonding map
$M_{i_k''}\to M_{i_{k-1}}$.
Amending $\psi_{k+1}$ by all $\alpha_\Delta$'s, we obtain that the compositions
of $\mu$ and the amended $\nu$ with the bonding map $P_{i_k'}\to P_{i_{k-1}}$
are homotopic (with values in $P_{i_{k-1}}$) by a homotopy $h''_k$.
Thus the compositions
$g_{k+1}\:M_{i_{k+1}}@>\psi_{k+1}>>P_{i_{k+1}'}@>p>>P_{i_{k-1}}$ are the integer
slices of
a level-preserving map $g_{[0,\infty)}\:M^{is}_{[0,\infty)}\to P^i_{[0,\infty)}$,
where $s(j)=j+2$ and $M^i_j=M_{i_j}$, whose restriction to $P^{is}_{[0,\infty)}$
coincides with the shift
$[p^{i_2}_{i_0},p^{i_\infty}_{i_\infty})\:
P^{is}_{[0,\infty)}\to P^i_{[0,\infty)}$.

By construction%
\footnote{Or, alternatively, by (iii)\imp (i) in Proposition 2.8(b).}%
, $g_{[0,\infty)}$ is properly homotopic to
$[p^{i_2}_{i_0},p^{i_\infty}_{i_\infty})\:
M^{is}_{[0,\infty)}\to M^i_{[0,\infty)}$.
Indeed, the amended $\psi_{k+1}$ is homotopic to the original one by a homotopy
$\Psi'_{k+1}$ with values in $M_{i_{k+1}'}$.
The homotopy $h''_k$ is homotopic to the composition of $h'_k$ and the bonding
map $P_{i_k''}\to P_{i_{k-1}}$ by a $2$-homotopy bounded by the compositions
of $\Psi'_k$ and $\Psi'_{k+1}$ with the bonding maps.
In turn, $h'_k$ is relatively homotopic to the composition of $h_k$ and
the bonding map $P_{i_k'}\to P_{i_k''}$.
Finally, $h_k$ is homotopic to the bonding map $M_{i_{k+1}}\to M_{i_k'}$ by
a $2$-homotopy bounded by the compositions of $\Psi_k$ and $\Psi_{k+1}$ with
the bonding maps. \qed
\enddemo

\proclaim{Lemma 3.7} Let $G_i$ be an inverse sequence of groups and $H_i$
an inverse sequence of their subgroups.

(a) \cite{Ir} If the inverse sequence of the pointed sets $G_i/H_i$ of right cosets
satisfies the Mittag-Leffler condition, then $\derlim H_i\to\derlim G_i$ is onto.

(b) If all the $G_i$ are countable, the converse holds.
\endproclaim

As the proof of Lemma 3.3 already uses both left and right actions, it
seems unlikely that Lemma 3.7(b) can be similarly proved by a mere count of
cardinalities.
We use the Baire Category Theorem (see Corollary 7.11, where it is deduced from
Bourbaki's ``Mittag-Leffler Theorem'').

Part (a) will be used a little later (in the proof of Theorem 3.15(b$'$)).

\demo{Proof. (a)} For variety, we indicate a geometric proof in the finitely
presented case.
Let $\dots\to (P_2,Q_2)\to (P_1,Q_1)$ be an inverse sequence of pairs of
compact connected polyhedra with $\pi_1(Q_i)=H_i$ and $\pi_1(P_i)=G_i$,
the inclusion $H_i\i G_i$ being induced by $Q_i\i P_i$.
Let $(X,A)=\invlim (P_i,Q_i)$ and write $S_i=\pi_1(P_i,Q_i)$.
Assuming that $\im[S_{i+1}\to S_i]=\im[S_j\to S_i]$ for all $i$ and all $j>i$,
let us prove that $\pi_0(A)\to\pi_0(X)$ is surjective.
Represent an element of $\pi_0(X,b)$, where $b\in A$, by a proper ray
$b_{[0,1]}\ell_1b_{[1,2]}\ell_2\dots\:[0,\infty)\to P_{[0,\infty)}$, where
$b_i=p^\infty_i(b)$ and the $\ell_i\:(I,\partial I)\to(P_i,b_i)$ are loops.
Then the path $b_{[0,1]}\ell_1b_{[1,2]}$ is homotopic with values in
$(P_{[0,2]},Q_{[0,2]};\,b_2)$ to a path
$\ell_1'\:(I,\partial I;\,\{0\})\to(P_2,Q_2;\,b_2)$.
Similarly, $\ell_1'\ell_2b_{[2,3]}$ is homotopic with values in
$(P_{[1,3]},Q_{[1,3]};\,b_3)$ to a path
$\ell_2'\:(I,\partial I;\,\{0\})\to(P_3,Q_3;\,b_3)$.
Combining all such homotopies together yields a homotopy between the original
proper ray and some proper ray with values in $Q_{[0,\infty)}$. \qed
\enddemo

\demo{(b)} Suppose that $G_i/H_i$ is not Mittag-Leffler.
Then there exists a $k$ such that the images of
$G_{k+i}/H_{k+i}$ in $G_k/H_k$ do not stabilize.
Note that these are the same as the images of the compositions
$G_{k+i}\to G_k\to G_k/H_k$.
Let $A_i$ be the image of $G_{k+i}$ in $G_k$ and let $B_i=A_i\cap H_k$.
Then the images of $A_i$ in $A_0/B_0$ do not stabilize; nor do the preimages
$A_iB_0$ of these images.
The anti-automorphism $g\mapsto g^{-1}$ sends them to $B_0A_i$, and it follows
that each pointed set $B_0\q A_i/A_{i+1}$ of double cosets is nontrivial.
Each $G_{k+i}\to A_i$ is onto, so $\derlim G_{k+i}\to\derlim A_i$ is onto.
The composition $\derlim H_{k+i}\to\derlim G_{k+i}\to\derlim A_i$ factors
through $\derlim B_i$, therefore $\derlim B_i\to\derlim A_i$ is surjective.
Then the proof of Lemma 3.3 shows that $A:=\invlim A_0/A_i$ is the orbit
of $B:=\invlim B_0/B_i$ under the left action of $A_0$ on $A$.

The inverse limit topology on $A$ is induced by the ultrametric $d$ defined by
$d(x,y)=\frac1n$ whenever $p^\infty_n(x)=p^\infty_n(y)\in A_0/A_n$ but
$p^\infty_{n+1}(x)\ne p^\infty_{n+1}(y)\in A_0/A_{n+1}$.
Clearly $A$ is complete in this metric (in fact it is the completion of
$A_0/\bigcap A_i$).
A ball of radius $\frac1n$ is a point-inverse of $p^\infty_n$.
If such a ball $(p^\infty_n)^{-1}(gA_n)$ intersects $B$, then $gA_n$ contains
some $b\in B_0$.
In this case let $B_0aA_{n+1}$ be a nontrivial element of $B_0\q A_n/A_{n+1}$.
Then $baA_{n+1}$ is contained in $gA_n$ and is disjoint from $B_0$.
Hence $(p^\infty_{n+1})^{-1}(baA_{n+1})$ is contained in
$(p^\infty_n)^{-1}(gA_n)$ and is disjoint from $B$.
We proved that $B$ is nowhere dense in $A$.
Hence by the Baire Category Theorem, $A$ cannot be the union of the countably
many translates $gB$ of $B$ under the left action of $A_0$ on $A$ by isometries.
This is a contradiction. \qed
\enddemo

Pontryagin considered the equivalence relation on inverse sequences of groups,
generated by the operation of passing to a subsequence \cite{P; Ch.~III, \S I}.
Similarly to the proof of Proposition 2.6, (ii)\iff (iii), two inverse
sequences $G_i$, $H_i$ of abelian groups are equivalent in Pontryagin's sense
iff they are related by a {\it pro-isomorphism}, which is a collection of maps
$f_i\:G_i\to H_i$ commuting with the bonding maps and such that there exist
an increasing sequence $k\:\N\to\N$ and homomorphisms
$g_i\:H_{k(i+1)}\to G_{k(i)}$ such that the compositions
$g_if_{k(i+1)}\:G_{k(i+1)}\to G_{k(i)}$ and $f_{k(i)}g_i\:H_{k(i+1)}\to H_{k(i)}$
equal the bonding maps.

\proclaim{Theorem 3.8} A collection of homomorphisms $f_i\:G_i\to H_i$
between inverse sequences of countable groups commuting with the bonding maps
is a pro-isomorphism if and only if $\invlim f_i\:\invlim G_i\to\invlim H_i$ and
$\derlim f_i\:\derlim G_i\to\derlim H_i$ are bijections.
\endproclaim

This answers a question of Koyama \cite{Koy; Problem 1}, who proved
the special case where each $f_i(G_i)$ is normal in $H_i$ \cite{Koy; Lemma 2}.
The finitely generated abelian case was first announced by Keesling
\cite{Ke; Theorem 2.4}.

Modulo Lemma 3.7(b), the proof of Theorem 3.8 is similar to that of
\cite{Koy; Lemma 2}, but since one's reading of that proof is complicated by
a misprint and several further references, for convenience we include
the details.

\demo{Proof} The ``only if'' part follows since for any inverse sequence of
groups $\Gamma_i$, the bonding maps induce bijections
$\invlim\Gamma_{i+1}\to\invlim\Gamma_i$ and
$\derlim\Gamma_{i+1}\to\derlim\Gamma_i$.

Conversely, each $f_i$ yields short exact sequences
$1\to K_i\to G_i\to f_i(G_i)\to 1$ and $1\to f_i(G_i)\to H_i\to L_i\to *$,
where $K_i=\ker f_i$ and $L_i$ is the pointed set $H_i/f(G_i)$ of right cosets.
By the hypothesis and the six-term exact sequence 3.1(d,d$'$),
$\invlim K_i=\derlim K_i=1$ and $\invlim L_i=*$.
By Lemma 3.7(b), the $L_i$ satisfy the Mittag-Leffler condition, and by Lemma 3.3
so do the $K_i$.
Let $k(0)=0$ and assume that $k(i)$ is defined.
Then by Lemma 3.4(a), there exist $j>k(i)$ and $k(i+1)>j$ such that
$L_{k(i+1)}\to L_j$ and $K_j\to K_{k(i)}$ are trivial.
Given an $h\in H_{k(i+1)}$, its image in $H_j$ equals $f_j(g_j)$ for some
$g_j\in G_j$.
Given any $g'_j\in G_j$ with $f_j(g'_j)=f_j(g_j)$, its image $g\in G_i$
equals the image of $g_j$.
Hence $g_i\:H_{k(i+1)}\to G_{k(i)}$, $h\mapsto g$, is well defined.
It is straightforward to verify that it is a homomorphism and that the
$g_i$'s invert the $h_i$'s as required. \qed
\enddemo

We now mention two known corollaries of Theorem 3.6; a third will appear
in Theorem 3.15(e).

The following was proved in \cite{DS2} and also follows from the
$\pi_\infty$-criterion of Siebenmann \cite{S1} as well as from \cite{EH; 5.5.6}.

\proclaim{Theorem 3.9}
Let $X$ and $Y$ be connected finite dimensional compacta.
A Steenrod homotopy class $f\:(X,x)\ssm(Y,y)$ is a Steenrod homotopy
equivalence iff it can be represented by a level-preserving map
$f_{[0,\infty)}\:(P_{[0,\infty)},p_{[0,\infty)})\to
(Q_{[0,\infty)},q_{[0,\infty)})$ that is an isomorphism in pro-homotopy.
\endproclaim

See Proposition 2.6 for the definition of an isomorphism in pro-homotopy.
An interesting alternative proof of Theorem 3.9 is given in \cite{G\"u}.

\demo{Proof} The ``only if'' part follows from Lemma 2.5(b).
Conversely, the induced maps $(f_i)_*\:\pi_n(P_i,p_i)\to\pi_n(Q_i,q_i)$
yield an isomorphism of pro-groups for each $n$.
So by the ``if'' part of Theorem 3.8 and by Theorem 3.1(b),
$f_*\:\pi_n(X)\to\pi_n(Y)$ is an isomorphism for all $n$.
Hence by Theorem 3.6 $f$ is a Steenrod homotopy equivalence. \qed
\enddemo

If $P$ is a non-compact polyhedron, then $\pi_n(P):=[(S^n,pt),(P,pt)]$
is isomorphic to $\spi_n(P)$ by Proposition 2.3(a).

\proclaim{Theorem 3.10} {\rm (Edwards--Geoghegan)} Let $X$ be an $m$-dimensional
connected compactum.
The following are equivalent:

\smallskip
(i) $\pi_n(X)$ is discrete for all $n\le m$;

(ii) $X$ is pointed Steenrod homotopy dominated by a compact polyhedron;

(iii) $X$ admits a map into a polyhedron inducing isomorphisms on $\pi_n$
for all $n\le m$;
\smallskip

(i$'$) $\pi_n(X)$ is discrete for all $n$;

(ii$'$) $X$ is pointed Steenrod homotopy dominated by a compact $m$-polyhedron;

(iii$'$) $X$ admits a map into an $(m+1)$-polyhedron inducing isomorphisms on
$\pi_n$ for all $n$.
\endproclaim

There is a number of different proofs of Theorem 3.10 in the literature.
Though many of them are quite obscure from the geometric viewpoint:
the original proof \cite{EG1}, \cite{EG2} (see also a simplification in
\cite{D3}) involves non-separable polyhedra; Dydak \cite{D2; 7.4} converts
bonding maps to fibrations; Ferry \cite{Fe2; Theorem 4 plus the first
paragraph of its own proof} works with Hilbert cube manifolds.
A closed geometric proof can be extracted from Kodama's papers \cite{Ko2},
\cite{Ko3}, which the present author learned of only when both proofs presented
below had been written down.
The second one may in fact be considered to be a simplified version of
Kodama's argument.

\demo{Proof. (i)\imp(ii$'$)} This will be proved by a combination of
Addendum (a) to Theorem 3.6 with the trick in the proof of
Proposition 3.5(b).
An alternative proof will be given below based on Theorem 3.12.

Suppose that $X$ is the limit of an inverse sequence of compact
$m$-polyhedra $P_i$.
Since $\pi_n(X)$ is Hausdorff for each $n\le m$, by Theorem 3.1(b) and
Lemma 3.3 the inverse sequence $\pi_n(P_i)$ is Mittag-Leffler for each
$n\le m+1$.
By reindexing, we may assume that for each $i$ and each $n\le m+1$, the image
$G^{[n]}_i$ of $\pi_n(P_{i+1})$ in $\pi_n(P_i)$ equals the image of each
$\pi_n(P_j)$ with $j>i$.
Then each $G^{[n]}_{i+1}$ surjects onto $G^{[n]}_i$ and so does $\pi_n(X)$,
for each $n\le m+1$.

Since $\pi_n(X)$ is discrete for $n\le m$, by Lemma 3.4(b) it injects into
$\pi_n(P_k)$ for some $k$, which may be assumed to be the same for all $n\le m$.
Set $q=k+m+1$.
Since $G^{[1]}_q$ surjects onto $G^{[1]}_{q-1}$, each $x\in\pi_1(P_q)$ has
the same image in $\pi_1(P_{q-1})$ with some $y_x\in G^{[1]}_q$.
Attaching a $2$-cell to $P_q$ along some representative loop of $xy_x^{-1}$
for each $x\in\pi_1(P_q)$, we obtain a (possibly noncompact) polyhedron
$P^{[1]}_q$ such that the bonding map $P_q\to P_{q-1}$ factors through
the inclusion $P_q\i P^{[1]}_q$, and $\pi_1(X)$ maps onto $\pi_1(P_q^{[1]})$.
If the cells are attached to $P_q\x[0,\infty)$ instead of to $P_q$, it is easy
to keep the result locally compact.
Next, since $\pi_2(P_{q-1})$ maps into $G^{[2]}_{q-2}$, so does
$\pi_2(P^{[1]}_q)$.
Hence each $x\in\pi_2(P^{[1]}_q)$ has the same image in $\pi_2(P_{q-2})$
as some element $y_x$ of the image of $G^{[2]}_q$ in $\pi_1(P_q)$.
Continuing the process, we will eventually factor the bonding map
$P_q\to P_k$ through a sequence of inclusions
$P_q\i P_q^{[1]}\i\dots\i P_q^{[m+1]}=:Q$ such that $\pi_n(X)$ maps onto
$\pi_n(Q)$ for each $n\le m+1$.

Since $\pi_n(X)$ injects into $\pi_n(P_k)$ for all $n\le m$, it also injects
into $\pi_n(Q)$ for all $n\le m$.
Thus the triple $X\to P_q\to Q$ satisfies the hypothesis of Addendum (a) to
Theorem 3.6, except that $Q$ need not be compact.
Nevertheless, the proof of Addendum 3.6(a) applies to show that $P_q$ pointed
Steenrod homotopy dominates $X$. \qed
\enddemo

\demo{(ii)\imp(iii) and (ii$'$)\imp(iii$'$)}
If $K$ is a compact polyhedron and $d\:(K,b)\ssm (X,x)$ and
$u\:(X,x)\ssm (K,b)$ are such that $du=[\id_X]$, let us consider the doubly
infinite mapping telescope $P:=\Tel(\dots@>q>>K@>q>>K@>q>>\dots)$, where
$q\:(K,b)\to (K,b)$ is a PL map representing the Steenrod homotopy class $ud$
by Proposition 2.3(a).
From the Mather trick (see \cite{FR}) and Proposition 2.3(a) we get
a proper map $X\x\R\to P$ such that the composition $X\x\{0\}\i X\x\R\to P$
induces an isomorphism on all $\pi_n$. \qed
\enddemo

\demo{(i)\imp(i$'$), (ii)\imp(ii$'$), (iii)\imp(i), (iii$'$)\imp(i$'$)}
Obvious. \qed
\enddemo

\remark{Remark} As observed in \cite{EG1}, compacta satisfying the equivalent
conditions of Theorem 3.10 need not have the shape of any compact polyhedron.
Indeed, let $P$ be a polyhedron that is homotopy dominated by a compact
polyhedron $K$, but is not homotopy equivalent to any compact polyhedron
(see \cite{FR}).
If $d\:K\to P$ and $u\:P\to K$ are such that $du\simeq\id_P$, the
limit $X$ of the inverse sequence $\hat K:=(\dots@>ud>>K@>ud>>K)$ is such that
the projection $f\:X\to P$ induces an isomorphism on Steenrod homotopy groups.
Indeed, it is not hard to see that the level-preserving maps
$d_{[0,\infty)}\:\hat K_{[0,\infty)}\to P\x [0,\infty)$ and
$u_{[0,\infty)}\:P\x [0,\infty)\to\hat K_{[0,\infty)}$, defined respectively
by $d$ and $u$ on the integer levels, and extended to all levels using
the homotopies $d\simeq dud$ and $udu\simeq u$, are mutually inverse in
semi-proper homotopy.
Now if there exists a Steenrod homotopy equivalence $g\:L\ssm X$ for
some compact polyhedron $L$, then $[f]g\:L\ssm P$ can be represented by a map
$h\:L\to P$ by Proposition 2.3(a).
Since $h$ induces isomorphisms of homotopy groups, it is a homotopy equivalence
by the classical Whitehead theorem, which is a contradiction.
\endremark

\remark{Remark}
Ferry proved that, allowing the case $m=\infty$, for a connected compactum $X$,
$\pi_n(X)$ is discrete for all $n<m$ iff $X$ has the shape of a compactum
that is LC$_n$ for all $n<m$ \cite{Fe2}.
The case $m=1$ was first proved by Krasinkiewicz (see \cite{DS1});
concerning the ``if'' part see \S6.
In particular, the following assertions can be added to the list of
Theorem 3.10:

\smallskip
(iv) $X$ has the shape of an LC$_m$ compactum;
\smallskip
(iv$'$) $X$ has the shape of an LC$_\infty$ compactum.
\medskip
\endremark

\proclaim{Lemma 3.11} {\rm (Dydak)} Let $A_i\to B_i\to C_i\to D_i$ be an exact
sequence of inverse sequences of groups; or more generally of pointed sets where
each $A_i$ is a group acting on $B_i$ so that the nonempty point-inverses of
$B_i\to C_i$ are precisely the orbits of this action.

If both the $A_i$ and the $C_i$ satisfy the Mittag-Leffler condition and the
$D_i$ satisfy the dual Mittag-Leffler condition, then the $B_i$ satisfy
the Mittag-Leffler condition.
\endproclaim

Lemma 3.11 was proved in \cite{D5}, where it is stated in a weaker form
(insufficient for its applications below).
In \S8 we shall need not only Lemma 3.11 but also its proof, so for convenience
we reproduce it below, aiming at relative readability.

\demo{Proof}
By reindexing, we may assume that $\im(A_{i+1}\to A_i)=\im(A_{i+2}\to A_i)$,
$\im(C_{i+1}\to C_i)=\im(C_{i+2}\to C_i)$ and
$\ker(D_i\to D_{i-1})=\ker(D_i\to D_{i-2})$ for each $i$.
It suffices to prove that $\im(B_{n+2}\to B_n)=\im(B_{n+3}\to B_n)$ for every $n$.

So let $b_n$ be the image of some $b_{n+2}\in B_{n+2}$.
Since $\im(C_{n+2}\to C_{n+1})=\im(C_{n+4}\to C_{n+1})$, the image
of $b_{n+2}$ in $C_{n+1}$ is the image of some $c_{n+4}\in C_{n+4}$.
Then the image of $c_{n+4}$ in $D_{n+1}$ equals the image of $b_{n+2}$ in
$D_{n+1}$, which is trivial due to the exactness of the rows.
Since $\ker(D_{n+4}\to D_{n+1})=\ker(D_{n+4}\to D_{n+3})$, the image of $c_{n+4}$
must be trivial already in $D_{n+3}$.
Hence the image of $c_{n+4}$ in $C_{n+3}$ is the image of some
$b_{n+3}\in B_{n+3}$.
By construction, $b_{n+2}$ and $b_{n+3}$ map to the same element in $C_{n+1}$.
Hence their images in $B_{n+1}$ are related by the action by some
$a_{n+1}\in A_{n+1}$.
Since $\im(A_{n+1}\to A_n)=\im(A_{n+3}\to A_n)$, the image of $a_{n+1}$ in $A_n$
is the image of some $a_{n+3}\in A_{n+3}$.
Finally, the result of the action of $a_{n+3}$ on $b_{n+3}$ maps onto $b_n$. \qed
\enddemo

\remark{Remark} The above argument shows that for a fixed $n$,
(i) if the images of $A_i$'s in $A_n$ stabilize, and $C_i$ are Mittag-Leffler
and $D_i$ are dual Mittag-Leffler, then the images of $B_i$'s in $B_n$
stabilize; (ii) if the images of $C_i$'s in $C_n$ stabilize, each
$A_{i+1}\to A_i$ is an isomorphism, and $D_i$ are dual Mittag-Leffler, then
the images of $B_i$'s in $B_n$ stabilize.
This will be used in \S8.
\endremark
\bigskip

If $P$ is a locally compact polyhedron, we say that $P$ is {\it properly
$m$-connected at infinity} if every compact subset $Q$ of $P$ is contained
in a compact subset $R$ of $P$ such that every proper map $\R^n\to\Cl(P\but R)$
with $n\le m$ extends to a proper map $\R^n\x [0,\infty)\to\Cl(P\but Q)$.
Clearly, this property implies local $n$-connectedness of
the one-point compactification of $P$.
The converse implication does not hold (see Theorem 6.12a below).

\proclaim{Theorem 3.12} {\rm (Kodama)} Let $X$ be the inverse limit of compact
polyhedra $P_i$, and fix some $m\ge 0$.
Then $\pi_n(X)$ is discrete for all $n<m$ if and only if $P_{[0,\infty)}$
is properly $m$-connected at infinity.
\endproclaim

Theorem 3.12 is essentially a restatement of some results of
\cite{Ko2}, \cite{Ko3}.
Kodama's arguments are visibly simplified due to our use of Lemmas 3.11 and 3.13.
The homological analogue of Theorem 3.12 was obtained by Dydak (see
Theorem 6.12b).

\demo{Proof} We only consider the case $m>0$.
Theorem 3.1(a) furnishes the following exact sequence of inverse
sequences:
$$\pi_n(X)\to\pi_n(P_{[k,\infty]})\to\pi_n(P_{[k,\infty]},X)\to\pi_{n-1}(X)\to
\pi_{n-1}(P_{[k,\infty]}).$$
It is easy to see that $P_{[k,\infty]}$ Steenrod deformation retracts onto
$P_k$, so their Steenrod homotopy groups (pointed sets when $n=0$) can be
identified.

If $\pi_n(X)$ is discrete for each $n<m$, then by Theorem 3.1(b) and
Lemma 3.3, $\pi_n(P_k)$ satisfy the Mittag-Leffler condition for each $n\le m$.
By Lemma 3.4(c), $\pi_n(P_k)$ satisfy the dual Mittag-Leffler condition
for each $n<m$.
Hence by Lemma 3.11, $\pi_n(P_{[k,\infty]},X)$ satisfy the Mittag-Leffler
condition for each $n\le m$.
On the other hand, by Theorem 3.1(c),
$\invlim\pi_n((P_{[k,\infty]},X)=\pi_n(X,X)$, which is trivial.
Then by Lemma 3.4(a), for each $i$ there exists a $j>i$ such that
$\pi_n(P_{[j,\infty]},X)\to\pi_n(P_{[i,\infty]},X)$ is trivial for each $n\le m$.

Conversely, suppose that for each $i$ there exists a $j>i$ such that
$\pi_n(P_{[j,\infty]},X)\to\pi_n(P_{[i,\infty]},X)$ is trivial.
Then $\pi_{n-1}(X)\to\pi_{n-1}(P_j)$ is injective, from the exact sequence above.
Hence $\pi_{n-1}(X)$ is discrete.

To complete the proof, we use Lemma 3.13 below.
Clearly, a base-ray preserving proper map $\R^n\to P_{[k,\infty)}$ extends to
a proper map $\R^n\x[0,\infty)\to P_{[k,\infty)}$ if and only if it is
``null-homotopic'' by a base-ray preserving proper homotopy.
Moreover, every proper map $\R^n\to P_{[k,\infty)}$ is properly homotopic
to a base-ray preserving one, as long as every base-ray preserving map
$\R^1\to P_{[k,\infty)}$ is ``null-homotopic'' by a base-ray preserving
proper homotopy. \qed
\enddemo

\proclaim{Lemma 3.13} Let $X$ be the inverse limit of compact
polyhedra $P_i$.
There is a natural pointed bijection between $\pi_n(P_{[k,\infty]},X)$ and
the set of base-ray preserving proper homotopy classes of base-ray preserving
proper maps $\R^n\to P_{[k,\infty)}$.
\endproclaim

The inverse sequence of the pointed sets (groups for $n>1$)
$\pi_n(P_{[k,\infty]},X)$ was considered by Kodama under the name
``the movability pro-group'' \cite{Ko1}; more precisely, Kodama defined
these by means of the statement of Lemma 3.13 and apparently did not
realize that they coincide with the relative Steenrod homotopy sets (which
were familiar to him).

\demo{Proof}
By definition, $\pi_n(P_{[k,\infty]},X)$ is the set of base-ray
preserving proper homotopy classes of base-ray preserving proper maps
$(\R^n\x [0,\infty),\R^n\x\{0\}\but B^n)\to
(P_{[k,\infty)}\x[0,\infty),P_{[k,\infty)}\x\{0\})$,
where $B^n\i\R^n\x\{0\}$ is an $n$-ball (we omit base rays from the notation).
Since $[0,\infty)$ deformation retracts onto $\{0\}$, we may assume without
loss of generality that $B^n$ does map into $P_{[k,\infty)}\x\{0\}$ under all
maps and homotopies in question.
This implies the assertion. \qed
\enddemo

Theorem 3.12 immediately implies

\proclaim{Corollary 3.14} {\rm (Dydak \cite{D4})} If $X$ and $Y$ are
compacta such that $\pi_i(X)$ and $\pi_i(Y)$ are discrete for all $i\le n$
and $\pi_i(X\cap Y)$ are discrete for all $i\le n-1$ then $\pi_i(X\cup Y)$
are discrete for all $i\le n$.
\endproclaim

Dydak's proof uses homology, universal covers and nearly
Steenrod connected compacta (concerning the latter, see Theorem 6.12).
Under the stronger hypothesis that $\pi_i(X\cap Y)$ is discrete for all
$i\le n$, Corollary 3.14 was obtained by Kodama \cite{Ko2}, \cite{Ko3}.
The case $n=0$ was obtained by Krasinkiewicz (see Theorem 8.3 below).

\demo{Alternative proof of Theorem 3.10, (i)\imp (ii$'$)}
In the terminology of \S6, this will be the Steenrod-theoretic version of
the Ferry Construction.
Specifically, we will construct a Steenrod retraction $P_{[k,\infty]}\ssm X$
for some $k$.
It will follow that $P_k$ Steenrod homotopy dominates $X$, for it is easy to
see that $P_{[k,\infty]}$ Steenrod deformation retracts onto $P_k$.

Let us fix a triangulation of $P_{[0,\infty]}$, where the diameters of
the simplices tend to zero as we approach $X$.
Given some $n$ and $k_n$, let $Q^n=(\dots\to Q^n_{k_n+1}\to Q^n_{k_n})$
denote the inverse sequence $\dots\to P_{[k_n,k_n+2]}^{(n)}\to
P_{[k_n,k_n+1]}^{(n)}\to P_{k_n}^{(n)}$, whose limit is the $n$-skeleton
$P_{[k_n,\infty]}^{(n)}$.
Then a Steenrod retraction $P_{[k_n,\infty]}^{(n)}\cup X\ssm X$ is the proper
homotopy class of a proper retraction
$r^n\:Q_{[0,\infty)}^n\cup P_{[0,\infty)}\to P_{[0,\infty)}$.
Since $\pi_n(X)$ is discrete for all $n<0$ (a vacuous condition), by
Theorem 3.12, for each $i$ there exists a $j=j(i)>i$ such that each vertex of
$P_{[j,\infty)}$ is the endpoint of a proper ray $[0,\infty)\to
P_{[i,\infty)}$.
This yields an $r^0$, with $k_0=j(0)$.
Assuming that $r^n$ has been defined, we can define an $r^{n+1}$ since by
Theorem 3.12, for each $i$ there exists a $j=j(i)>i$ such that every proper
map $f\:S^n\x [0,\infty)\to P_{[j,\infty)}$ with
$f(S^n\x\{0\})=\partial\Delta$ for some $(n+1)$-simplex $\Delta$ of
$P_{[j,\infty)}$ extends to a proper map
$\bar f\:B^{n+1}\x [0,\infty)\to P_{[i,\infty)}$ such that
$\bar f(B^{n+1}\x\{0\})=\Delta$.
Finally, $[r^m]\:P_{[k,\infty]}\ssm X$ is the desired Steenrod retraction. \qed
\enddemo

\definition{Steenrod fibrations}
Let $P=(\dots@>p_1>>P_1@>p_0>>P_0)$ and $Q=(\dots@>q_1>>Q_1@>q_0>>Q_0)$ be
inverse sequences of compact polyhedra and PL maps and let $w\:P\to Q$ be
a {\it PL level map}, i.e.\ a sequence of PL maps $w_i\:P_i\to Q_i$ such that
$w_ip_i=q_iw_{i+1}$ for each $i$.
We say that $w$ satisfies the {\it homotopy lifting property} if for each $i$
there exists a $j>i$ such that given a compact polyhedron $K$, a homotopy
$f_t\:K\to Q_j$ and a lift $\bar f_0\:K\to P_j$ of $\bar f_0$ (i.e.\
$w_j\bar f_0=f_0$), there exists a lift $\tl f_t\:K\to P_i$ of $q^j_i f_t$
(i.e.\ $w_i\tl f_t=q^j_i f_t$) such that $\tl f_0=p^j_i\bar f_0$.
Finally, a map $f\:E\to B$ between compacta is called a {\it Steenrod fibration}
if $E=\invlim P$, $B=\invlim Q$ and $f=\invlim w$ for some $P$, $Q$ and $w$ as
above, with $w$ satisfying the homotopy lifting property.

Clearly, the homotopy lifting property is satisfied, with $j(i)=i$, if
\roster
\item"(i)" each $f_i$ is a (Serre) fibration.
\endroster
In particular, this is the case if
\roster
\item"(ii)" each $f_{i+1}$ is the composition of
a fibration $\phi_{i+1}\:P_{i+1}\to R_{i+1}$ and the pullback
$q_i^*(f_i)\:R_{i+1}\to Q_{i+1}$ of $f_i$.
\endroster
In the latter situation (ii), $f$ will be a fibration, as it is identified with
the map $\psi^\infty_0$ corresponding to the decomposition
$E=\invlim(\dots@>\psi_1>>B_1@>\psi_0>>B_0=B)$, where each $\psi_i$ is
a pullback of $\phi_i$.
In the more general situation (i), $f$ does not have to be a fibration, as shown
by the following examples.

The PL map $f\:[0,2]\to [0,1]$, $x\mapsto\min(x,1)$, is a Steenrod
fibration (even with the $f_i$ being fibrations) but not a fibration.
Another such example is the projection of the topologist's sine curve (see
Example 2.10) onto $[0,1]$.
On the other hand, if $B$ is a compactum that contains no nondegenerate paths
(for instance, the pseudo-arc), every map $f\:E\to B$ is a fibration.
As observed in \cite{MR; Example 4}, it does not have to be a Steenrod fibration,
for instance, when $E=B\vee B$ (for any choice of a non-isolated basepoint in $B$)
and $f$ is the folding map.
\enddefinition

\proclaim{Theorem 3.15} Let $f\:E\to B$ be a map between compacta.

\medskip
(a) Suppose that $E=\invlim R$ and $B=\invlim Q$, where
$R=(\dots@>r_1>>R_1@>r_0>>R_0)$ and $Q=(\dots@>q_1>>Q_1@>q_0>>Q_0)$ are inverse
sequences of compact polyhedra and PL maps.
Then there exist an inverse sequence $P=(\dots@>p_1>>P_1@>p_0>>P_0)$ of compact
polyhedra and PL maps with $\invlim P=E$ and PL level maps $u\:P\to R$ and
$w\:P\to Q$ such that $\invlim u=\id_E$ and $\invlim w=f$.

\medskip
(b) \cite{C2; 3.3} Suppose that $f$ is a Steenrod fibration and $b\in B$ is
such that $F:=f^{-1}(b)$ is non-empty; pick some $\tl b\in F$.
Then there are natural actions of $\pi_1(E)$ on $\pi_n(F)$ for $n\ge 1$ and
of $\pi_1(B)$ on $\pi_0(F)$ and a natural exact sequence of pointed sets
$$\dots\to\pi_1(F)\to\pi_1(E)@>f_*>>\pi_1(B)\to\pi_0(F)\to\pi_0(E)@>f_*>>\pi_0(B)$$
whose maps to the left of $\pi_1(B)$ (resp.\ $\pi_2(B)$) are homomorphisms of
groups (resp.\ right $\Z\pi_1(E)$-modules).
In addition, $\pi_2(B)\to\pi_1(F)$ is also $\pi_1(E)$-equivariant,
$\pi_1(F)\to\pi_1(E)$ is a crossed module, $\pi_1(B)\to\pi_0(F)$ is
$\pi_1(B)$-equivariant with respect to the right regular action on $\pi_1(B)$,
and the non-trivial point-inverses of $\pi_1(F)\to\pi_0(E)$ coincide with
the orbits of $\pi_1(B)$.

\medskip
(b$'$) In the notation of (b), if $B$ is connected, $f_*\:\pi_0(E)\to\pi_0(B)$
is surjective.

\medskip
(c) \cite{C2; 3.4} If $f$ is a Steenrod fibration and $b,b'\in B$ represent
the same element of $\pi_0(B)$, then $f^{-1}(b)$ and $f^{-1}(b')$ have
the same shape.

\medskip
(d) If $E$ is finite-dimensional and $f$ is {\rm cell-like}, that is,
$f^{-1}(b)$ is cell-like for each $b\in B$, then $f$ is a Steenrod fibration.

\medskip
(e) \cite{DS2} A cell-like map between connected finite-dimensional compacta is
a Steenrod homotopy equivalence.
\endproclaim

From part (a) it follows that the class of Steenrod fibrations coincides with
the (a priori wider) class of shape fibrations, which were introduced in \cite{MR}.%
\footnote{Let us show more: If $f$ is a shape fibration, and $R$, $Q$ are as
in (a), then there exist $P$, $u$ and $w$ as in (a) and such that $w$ satisfies
the homotopy lifting property.
Indeed, by (a) we may assume that there is a PL level map $v\:R\to Q$ with
$\invlim v=f$.
Then by \cite{MR; Theorem 1}, $v$ satisfies the AHLP of \cite{MR}.
The required $P$, $u$ and $w$ are then furnished by the proof of
\cite{MR; Theorem 2}.}

A weakened version of (c) was initially proved in \cite{MR; Theorem 3}.
Part (d) strengthens \cite{MR; Theorem 5}.
Weakened versions of (d) were initially obtained by Dydak (1975) and
K. Morita (1975) (see \cite{D2; 8.5}), improving earlier results by Bogatyj
\cite{Bog}, Y. Kodama (1974) and Kuperberg \cite{Ku1}.

The assumption of finite dimensionality is necessary in (e)
(see \cite{DS1; 10.3.1}) and in (d) \cite{MR; Example 6}.

\demo{Proof. (a)} Let $\Gamma\i E\x B$ be the graph of $f$.
Let $P_i$ be a polyhedral neighborhood of the image of $\Gamma$ in $R_i\x Q_i$
such that each $P_{i+1}$ maps into $P_i$ and $\invlim P_i=\Gamma$.
Define $u_i\:P_i\to R_i$ and $w_i\:P_i\to Q_i$ to be the restrictions of
the projections onto the factors of $R_i\x Q_i$. \qed
\enddemo

\demo{(b)} Let $P$, $Q$ and $w$ be as in the definition of a Steenrod fibration,
with $j(i)=i+1$ in the homotopy lifting property.

Let us prove a homotopy lifting property for level-preserving representatives of
Steenrod homotopy classes.
Consider a polyhedron $X$ and a level-preserving map
$\phi_{[0,\infty)}\:X\x I\x[0,\infty)\to Q_{[0,\infty)}$ along with a lift
$\bar\phi_{[0,\infty]}\:X\x[0,\infty)\to P_{[0,\infty)}$ of its restriction to
$X\x\{0\}\x[0,\infty)$.
By the hypothesis, every its shifted level $q_{i-1}\phi_i$ lifts to a map
$\tl\phi_{i-1}\:X\x I\to P_i$ restricting to $p_{i-1}\bar\phi_i$ on $X\x\{0\}$.
Since $q^{[i,i+1]}_{[i-2,i-1]}\:Q_{[i,i+1]}\to Q_{[i-2,i-1]}$ factors through
$q_{i-1}\x\id_I\:Q_i\x I\to Q_{i-1}\x I$, we also obtain a lift
$\tl\phi_{[i-2,i-1]}$ of $q^{[i,i+1]}_{[i-2,i-1]}\phi_{[i,i+1]}$, restricting
to $p^{[i,i+1]}_{[i-2,i-1]}\bar\phi_{[i,i+1]}$ on $X\x\{0\}\x[i,i+1]$.
Two consecutive instances of this construction yield two possibly different
lifts $\tl\phi'_{i-2}$ and $\tl\phi''_{i-2}$ of $q^i_{i-2}\phi_i$ restricting
to $p^i_{i-2}\bar\phi_i$ on $X\x\{0\}$.
By the hypothesis, $p_{i-3}\tl\phi_{i-2}'$ is homotopic to
$p_{i-3}\tl\phi_{i-1}''$ through lifts of $q^i_{i-3}\phi_i$ that restrict to
$p^i_{i-3}\bar\phi_i$ on $X\x\{0\}$.
It follows that $p_{[i-3,i-2]}\tl\phi_{[i-2,i-1]}$ is homotopic with support
over $[i-3,i-2.5]$ to a new lift $\tl\phi_{[i-3,i-2]}'$ of
$q^{[i,i+1]}_{[i-3,i-2]}\phi_{[i,i+1]}$ that now agrees with its neighbor
$\tl\phi_{[i-4,i-3]}'$.
Thus we obtain a lift $\tl\phi_{[0,\infty)}\:X\x[0,\infty)\to P_{[0,\infty)}$ of
$q^{[0,\infty)}_{[-3,\infty-3]}\phi_{[0,\infty)}$ restricting to
$p^{[0,\infty)}_{[-3,\infty-3]}\bar\phi_{[0,\infty)}$ on $X\x\{0\}\x[0,\infty)$.

The thus established property easily implies that the inclusion induced map
$\pi_n(CF,F;\,\tl b)\to\pi_n(MC(f),E;\,\tl b)$ is a bijection for $n\ge 1$.%
\footnote{Another possible approach is to observe that
$f_*\:\pi_n(E,F;\,\tl b)\to\pi_n(B;b)$ is a bijection.
This leads to weaker conclusions: the exact sequence becomes shorter by one term;
the $\pi_1(E)$-actions reduce, via the restriction of scalars, to
$\pi_1(F)$-actions; and the $\pi_1(B)$-action reduces to a $\pi_1(E)$-action.}
Since $CF$ is cell-like, being an inverse limit of cones,
$\partial\:\pi_n(CF,F;\,\tl b)\to\pi_{n-1}(F;\,\tl b)$ is a bijection by
Theorem 3.1(a).
Hence the assertion follows from the homotopy exact sequence of $(MC(f),E)$,
given by Theorem 3.1(a). \qed
\enddemo

\demo{(b$'$)} We may assume that $E$ is connected; else we may replace it
with the connected component containing $\tl b$ and verify that the restriction
of $f$ is a Steenrod fibration.
Let $P$, $Q$ and $w$ be as in the definition of a Steenrod fibration,
with $j(i)=i+1$ in the homotopy lifting property.
Then by Theorem 3.1(b), $\pi_0(E)\to\pi_0(B)$ can be identified with
$\derlim (w_i)_*\:\derlim\pi_1(P_i)\to\derlim\pi_1(Q_i)$.
Without loss of generality, all the $P_i$ are connected, so the cokernel of
$(w_i)_*\:\pi_1(P_i)\to\pi_1(Q_i)$ can be identified with
$S_i:=\pi_0(MC(w_i),P_i)$.
By Lemma 3.7(a), it suffices to show that the inverse
sequence $\dots\to S_1\to S_0$ satisfies the Mittag-Leffler condition.
Indeed, let $F_i=w_i^{-1}(b_i)$, where $b_i=p^\infty_i(b)$.
By the homotopy lifting property, each bonding map $S_{i+1}\to S_i$ factors
through $T_i:=\pi_1(MC(w_i|_{F_i}),F_i)$.
Each $T_i$ can be identified with $\pi_0(F_i)$ and so is finite.
Hence $S_i$ satisfy the Mittag-Leffler condition. \qed
\enddemo

\demo{(c)} Let $P$, $Q$ and $w$ be as in the definition of a Steenrod fibration,
with $j(i)=i+1$ in the homotopy lifting property.
Let $R_i=w_i^{-1}(b_i)$ and $R'_i=w_i^{-1}(b'_i)$, where $b_i=q^\infty_i(b)$
and $b_i'=q^\infty_i(b')$, and let $v\:R_{[0,\infty)}\to [0,\infty)$ be
the level-preserving projection.
Pick a level-preserving Steenrod path $\ell_t\:[0,\infty)\to Q_{[0,\infty)}$
connecting $b_{[0,\infty)}$ with $b'_{[0,\infty)}$.
Then by the proof of (b), the proper homotopy
$q^{[0,\infty)}_{(-3,\infty-3)}v\ell_t\:R_{[0,\infty)}\to Q_{[0,\infty)}$ lifts
to a homotopy $h_t\:R_{[0,\infty)}\to P_{[0,\infty)}$ such that $h_0$ is
the inclusion composed with $p^{[0,\infty)}_{(-3,\infty-3)}$.
In particular, we obtain a proper map $h_1\:R_{[0,\infty)}\to R'_{[0,\infty)}$;
similarly, one constructs $h'_1\:R'_{[0,\infty)}\to R_{[0,\infty)}$.
Since the folding map $[-1,1]\to [0,1]$, $x\mapsto |x|$, is null-homotopic
keeping the endpoints fixed, $h_1h_1'$ and $h'_1h_1$ are properly homotopic
to the identity. \qed
\enddemo

\demo{(d)} Let $P$, $Q$ and $w$ be as in (a), with $\dim P_i\le m$ for each $i$.
Let us triangulate each $Q_i$ so that each $q_i$ is simplicial with respect to
the triangulation of $Q_{i+1}$ and some subdivision of the second barycentric
subdivision of the triangulation of $Q_i$.
Pick a point $v\in B$, and for each $i$ let $N_i^v$ be the simplicial
neighborhood of $q^\infty_i(v)$ in $Q_i$.
Then each $q_i(N_{i+1}^v)\i N^v_i$ and $\invlim N^v_i=\{v\}$.
Moreover, the union of all $N_{i+1}^u$'s that intersect $N_{i+1}^v$ is mapped
into $N_i^v$ due to our choice of the triangulations.
Let $M^v_i=w_i^{-1}(N^v_i)$.
Using that $f^{-1}(v)=\invlim M^v_i$ is cell-like and that $\dim P_i\le m$, by
the proof of Proposition 3.5(a) we may assume that each
$p_i|_{M^v_{i+1}}\:M^v_{i+1}\to M^v_i$ is null-homotopic.

Let $Q'_i=\bigcup_{v\in B}N_i^v$, in other words the simplicial neighborhood of
$q^\infty_i(B)$, and let $P'_i=\bigcup_{v\in B}N_i^v\x M_i^v\i Q_i'\x P_i$
(note that the union is finite).
Let $p'_i=(q_i\x p_i)|_{P'_i}$ and $q'_i=q_i|_{Q'_i}$, and let
$w'_i\:P'_i\to Q'_i$ be the restriction of the projection.
Given a linear map $\phi\:\Delta^{n+1}\to Q'_{i+2}$ and a partial lift
$\bar\phi_\partial\:\partial\Delta^{n+1}\to P'_{i+2}$, we have
$\phi(\Delta^{n+1})\i N^v_{i+2}$ for some $v\in X$ and therefore
$\bar\phi_\partial(\partial\Delta^{n+1})$ lies in the union of all products
$N_{i+2}^u\x M^u_{i+2}$ where $N_{i+2}^u$ intersects $N_{i+2}^v$.
Then $p_{i+1}'\bar\phi_\partial(\partial\Delta^{n+1})\i N_{i+1}^v\x M^v_{i+1}$.
Let $\tl\phi_\partial\:\partial\Delta^{n+1}\to M^v_{i+1}$ be the composition
of $p_{i+1}'\bar\phi_\partial$ and the projection.
Then the shift $p_i\tl\phi_\partial\:\partial\Delta^{n+1}\to M^v_i$ extends to
a map $\tl\phi\:\Delta^{n+1}\to M^v_i$.
Let $\bar\phi:=(q^{i+2}_i\phi)\x\tl\phi\:\Delta^{n+1}\to N^v_i\x M^v_i$.
Thus $\bar\phi\:\Delta^{n+1}\to P'_i$ is a lift of $q'_iq'_{i+1}\phi$ extending
$p'_ip'_{i+1}\bar\phi_\partial$.
In particular, $w'\:P'\to Q'$ satisfies the homotopy lifting property. \qed
\enddemo

\demo{(e)} This follows from parts (b) and (d) combined with Theorem 3.6. \qed
\enddemo

\head 4. Homology and cohomology \endhead

\definition{Steenrod homology}
Following \cite{St; \S IV}, the group $H_n(X,A)$ is defined to be the locally
finite homology group $H_{n+1}^\lf(P_{[0,\infty)},\,P_0\cup Q_{[0,\infty)})$,
where $(X,A)$ is the inverse limit of the compact polyhedral pairs $(P_i,Q_i)$.
We recall that the locally finite homology
$H_i^\lf(P_{[0,\infty)},\,P_0\cup Q_{[0,\infty)})$ can be defined as the group
of oriented proper pseudo-bordism classes of proper maps from oriented
$i$-pseudo-manifolds $(M,\partial M)$ into
$(P_{[0,\infty)},\,P_0\cup Q_{[0,\infty)})$.
Since locally finite homology is an invariant of proper homotopy type,
$H_n(X,A)$ is a shape invariant of the pair $(X,A)$ by Corollary 2.2.
(This in particular entails the topological invariance of Steenrod
homology.)
The induced maps $f_*$, where $f$ is a map or a Steenrod homotopy class, and
the boundary maps $\delta_*$ (which are also Steenrod homotopy invariants)
are defined in the obvious ways.
\enddefinition

\proclaim{Theorem 4.1} {\rm (Milnor \cite{Mi})} Steenrod homology
satisfies all Eilenberg--Steenrod axioms, as well as the following:

\smallskip
(i) Map Excision Axiom: $f_*\:H_n(X,A)\to H_n(Y,B)$ is an isomorphism
for any map $f\:(X,A)\to(Y,B)$ that restricts to a homeomorphism
$X\but A\to Y\but B$.

\smallskip
(ii) Cluster Axiom: $H_n((\bigsqcup X_i\but pt)^+,pt)\simeq\prod H_n(X_i,pt)$
naturally for any compacta $X_1,X_2,\dots$, where $Y^+$ denotes the one-point
compactification of $Y$.

\smallskip
(iii) Milnor's Exact Sequence: if
$(X,A)=\invlim(\dots@>q^1>>(Y_1,Z_1)@>q^0>>(X_0,Z_0))$, where $(Y_i,Z_i)$ are
pairs of compacta, there is a natural short exact sequence of groups
$$0\to\derlim H_{n+1}(Y_i,Z_i)\to H_n(X,A)\to\invlim H_n(Y_i,Z_i)\to 0.$$

(iv) Steenrod homology is the only ordinary homology theory (in the sense of
the Eilenberg--Steenrod axioms) on pairs of compacta that satisfies
(ii) for compact polyhedra $X_i$ and (i).
\endproclaim

Note that (ii) is a special case of (iii), using that the {\it cluster}
$(\bigsqcup X_i\but pt)^+$ is homeomorphic to
$\invlim(\dots\to X_1\vee X_2\vee X_3\to X_1\vee X_2\to X_1)$.

A version of (iii) is found already in \cite{St; Theorem 8}.

\definition{Pontryagin cohomology}
The group $H^n(X,A)$ is defined to be the compactly supported cohomology
group $H^{n+1}_\c(P_{[0,\infty)},\,P_0\cup Q_{[0,\infty)})$, where $(X,A)$ is
the inverse limit of compact polyhedral pairs $(P_i,Q_i)$.
Remarks similar to those about $H_n(X,A)$ apply, and $H_n(X,A;\,G)$ can be
similarly defined.
In addition, we have
$$H^{n+1}_\c(P_{[0,\infty)},P_0)\simeq\dirlim
H^{n+1}(P_{[0,\infty)},P_0\cup P_{[k,\infty)})
\simeq\dirlim H^n(P_{[k,\infty)})\simeq\dirlim H^n(P_k)$$
where the second isomorphism uses that $H^{n+1}(P_{[0,\infty)},P_0)=0$ since
$P_{[0,\infty)}$ deformation retracts onto $P_0$; and the third that
$P_{[k,\infty)}$ deformation retracts onto $P_k$.
The composite isomorphism
$$H^n(X)\simeq\dirlim H^n(P_k)$$
constitutes the original definition of $H^n(X)$, which for finite-dimensional
compacta is due to Pontryagin \cite{P; Ch.~III, \S II}.
As a corollary of this isomorphism and Proposition 2.3(b), we obtain
\enddefinition

\proclaim{Proposition 4.2} If $G$ is a countable abelian group, there
is a natural isomorphism $H^n(X;\,G)\simeq [X,\,K(G,n)]$.
\endproclaim

The cohomological analog of Theorem 4.1 holds, with the direct product
in (ii) replaced by the direct sum, and with the exact sequence of (iii)
replaced by a natural isomorphism $H^n(X,A)\simeq\dirlim H^n(Y_i,Z_i)$
\cite{Mi}.

\proclaim{Theorem 4.3} If $X\i\R^m$ is compact, there are natural (with
respect to inclusion) isomorphisms
$H^n(X)\simeq\tl H_{m-n-1}(\R^m\but X)$ and
$H_n(X)\simeq\tl H^{m-n-1}(\R^m\but X)$.
\endproclaim

The second isomorphism is due to Steenrod \cite{St}.
The first isomorphism traces back to the original duality theorem of Alexander
(1922), which was stated for $\bmod 2$ Betti numbers and for polyhedra $X$
(possibly wildly embedded).
These restrictions were gradually removed by Alexandroff, Frankl, Lefschetz and
eventually by Pontryagin (see \cite{P} and references there).

It is well-known that upon adoption of the map excision axiom these
isomorphisms become obvious, modulo the Poincar\'e duality on the non-compact
manifold $S^m\but X$ (see e.g.\ \cite{Fe3}).
Since the verification of the map excision axiom takes more than a few lines,
it seems worthwhile to be also aware of a direct version of this argument:

\demo{Proof} Let $\dots\i U_2\i U_1\i U_0$ be a sequence of polyhedral
neighborhoods of $X$ with $\bigcap U_i=X$.
Then the definition of $H_n(X)$ and an appropriate kind of the Poincar\'e
duality yield $$H_n(X)\simeq H_{n+1}^\lf(U_{[0,\infty)},\,U_0)\simeq
H^{m-n}(U_{[0,\infty)},\,\Fr\, U_{[0,\infty)}).$$
Using excision, the exact sequence of a pair, and the vanishing of
$\tl H^i(\R^m\x[0,\infty))$, the latter group is isomorphic to
$$\tl H^{m-n-1}(\R^m\x[0,\infty)\but U_{[0,\infty)})\simeq
\tl H^{m-n-1}(\R^m\but X).$$
The latter isomorphism holds since the projection
$\R^m\x[0,\infty)\but U_{[0,\infty)}\to\R^m\but X$ is a homotopy equivalence.
Similarly,
$$\multline H^n(X)\simeq H^{n+1}_\c(U_{[0,\infty)},\,U_0)\simeq
H_{m-n}(U_{[0,\infty)},\,\Fr\, U_{[0,\infty)})\\
\simeq\tl H_{m-n-1}(\R^m\x[0,\infty)\but U_{[0,\infty)})\simeq
\tl H_{m-n-1}(\R^m\but X). \qed\endmultline$$
\enddemo

\remark{Remark}
The above argument yields an elementary proof of the Jordan curve theorem,
using only basic facts about (co)homology of polyhedra; the Poincar\'e duality;
and Lemma 2.1.
(The latter is needed since $H^1(S^1)=H^2_\c(P_{[0,\infty)},P_0)$ is
computed twice: with $P_i=U_i$ and with $P_i=S^1$.)
\endremark
\medskip

By analogy with the direct limit interpretation of Pontryagin cohomology, it is
convenient to reserve a special notation for the inverse limits that appear in
Milnor's short exact sequences:

\definition{\v Cech homotopy and homology} If $X$ is the inverse limit of
polyhedra $P_i$, consider the topological group (pointed space for $n=0$)
$\check\pi_n(X,x):=\invlim\pi_n(P_i,p_i)$ and the topological group
$\check H_n(X):=\invlim H_n(P_i)$.
By Proposition 2.6, $\check\pi_n$ (resp\. $\check H_n$) and indeed the entire
short exact sequence 3.1(b) (resp\. 4.1(iii)) is a shape invariant of $(X,x)$
(resp.\ $X$).
The canonical epimorphisms $\pi_n(X)\to\check\pi_n(X)$ and
$H_n(X)\to\check H_n(X)$ will be denoted by $\check\tau$.
\enddefinition

Here is the (absolute) ``Hurewicz theorem in Steenrod homotopy''.

\proclaim{Theorem 4.4} Let $m>1$ and let $X$ be a compactum with $\pi_n(X)=0$
for all $n<m$.
Then the Hurewicz homomorphism $\pi_n(X)\to H_n(X)$ is an isomorphism for
$n\le m$ and an epimorphism for $n=m+1$.
\endproclaim

The first assertion of Theorem 4.4 was conjectured in \cite{Ch; p.\ 300} and
proved 40 years later by Lisitsa \cite{Li}.
His proof of the second assertion, however, contains an error: in the case
$m=2$ it depends on the claim that if $X$ is the inverse limit of polyhedra
$P_i$, then already $\derlim\pi_{m+2}(P_i)\to\derlim H_{m+2}(P_i)$ is
allegedly an epimorphism.
But this claim is false for $m=2$.
Indeed, let $P_i$ be the cone of the composition $S^3@>p^i>>S^3@>h>>S^2$
of a degree $p^i$ map and the Hopf map.
The identity on $S^2$ extends to a map $f_i\:P_{i+1}\to P_i$; let
$X=\invlim(\dots\to P_2\to P_1)$.
Clearly each $(f_i)_*\:H_4(P_{i+1})\to H_4(P_i)$ is the inclusion $\Z@>p>>\Z$
onto the index $p$ subgroup, so $\derlim H_4(P_i)\ne 0$.
Since $P_i$ is simply-connected, every map $S^4\to P_i$ is homotopic to
one where the cone vertex has precisely one point-inverse; it follows that
$\pi_4(hp^i)$ surjects onto $\pi_4(P_i)$.
From the exact sequence $$\pi_4(S^2)@>>>\pi_4(hp^i)@>>>
\pi_3(S^3)@>hp^i>>\pi_3(S^2),$$
$\pi_4(hp^i)$ is a quotient of $\pi_4(S^2)\simeq\Z/2$.
So each $\pi_4(P_i)$ is finite and $\derlim\pi_4(P_i)=0$.
Thus $\derlim\pi_4(P_i)\to\derlim H_4(P_i)$ cannot be surjective.

\remark{Remark}
The problem can be fixed in this example by using the surjectivity of
the connecting homomorphism $\invlim\pi_3(P_i)\to\derlim H_4(P_i)$ for
the short exact sequence of inverse sequences
$$0\to H_4(P_i)\to\Gamma_3(P_i)\to\pi_3(P_i)\to 0$$
obtained from Lemma 4.5 below.
Indeed, each $\Gamma_3(P_{i+1})\to\Gamma_3(P_i)$ is an isomorphism (between
copies of $\Z$).
Hence $\derlim\Gamma_3(P_i)=0$ and so the connecting homomorphism must be
surjective.
\endremark

\proclaim{Lemma 4.5} {\rm (Whitehead \cite{Wh})} Let $K$ be a simply-connected
polyhedron with a fixed triangulation, and let $\Gamma_n(K)$ be the image of
$\pi_n(K^{(n-1)})\to\pi_n(K^{(n)})$.
The Hurewicz homomorphisms fit into a long exact sequence
$$\dots\to\Gamma_n(K)\to\pi_n(K)\to H_n(K)\to\Gamma_{n-1}(K)\to\dots$$
\endproclaim

The use of geometric language enables us to simplify Whitehead's argument.

\demo{Proof} The assertion of the lemma becomes obvious if $\Gamma_n$ is replaced
with the group $\Gamma_n'$ of $\partial$-spherical singular oriented
pseudo-bordism classes of $\partial$-spherical singular oriented
$(n+1)$-pseudo-manifolds in $K$.
(An $n$-pseudo-manifold $M$ with boundary is {\it $\partial$-spherical}
if $\partial M$ is PL homeomorphic to $S^{n-1}$.
A pseudo-bordism $W$ between such $M_0$, $M_1$ is {\it $\partial$-spherical}
if the closure of $\partial W\but(\partial M_0\cup\partial M_1)$ is PL
homeomorphic to $S^{n-1}\x I$.
The group operation on pseudo-bordism classes is given by boundary connected
sum.)
Since $K$ is simply-connected, $\Gamma_n'$ is isomorphic to the group
$\Gamma_n''$ of simply-connected $\partial$-spherical singular oriented
pseudo-bordism classes of simply-connected $\partial$-spherical singular
oriented $(n+1)$-pseudo-manifolds in $K$.

If $M$ is a $\partial$-spherical oriented $(n+1)$-pseudo-manifold, it
collapses onto a union $L$ of $M^{(n-1)}$ with some $n$-simplices of $M$.
If additionally $M$ is simply-connected, the Hurewicz homomorphism
$\pi_n(L,\,L^{(n-1)})\to H_n(L,\,L^{(n-1)})$ is an isomorphism.
Next, $[\partial M]\in H_n(M)$ is trivial, in particular it maps trivially
to $H_n(L,\,L^{(n-1)})$.
It follows that the inclusion $\partial M\i M$ is homotopic to a map
$\partial M\to L^{(n-1)}=M^{(n-1)}$.
Similarly if $W$ is a simply-connected $\partial$-spherical oriented
pseudo-bordism between $M_0$ and $M_1$, any maps $\partial M_0\to M_0^{(n-1)}$
and $\partial M_1\to M_1^{(n-1)}$ homotopic to the inclusions
$\partial M_0\i M_0$ and $\partial M_1\i M_1$ are homotopic with values in
$W^{(n)}$. \qed
\enddemo

\demo{Proof of Theorem 4.4} By Proposition 3.5(a), $X$ is an inverse limit of
$(m-1)$-connected polyhedra $P_i$.
Obviously, the Hurewicz homomorphisms fit into the long exact sequence
$$\dots\to\Gamma_n'(X)\to\pi_n(X)\to H_n(X)\to\Gamma_{n-1}'(X)\to\dots$$
where $\Gamma_n'(X)$ is the group of the classes of proper singular oriented
pseudo-bordism, restricting to proper homotopy on the boundary, of proper
singular oriented $(n+2)$-pseudo-manifolds with boundary $S^n\x [0,\infty)$
in $P_{[0,\infty)}$.
Similarly to the proof of Theorem 3.1(b), there is a short exact sequence
$$0\to\derlim\Gamma_{n+1}'(P_i)\to\Gamma_n'(X)\to\invlim\Gamma_n'(P_i)\to 0$$
where $\Gamma_n'(P_i)$ is defined in the proof of Lemma 4.5 and
coincides with $\Gamma_n(P_i)$ as shown in that proof.
Since each $P_i$ is $(m-1)$-connected, $\Gamma_n(P_i)=0$ for $n\le m$.
From the hypothesis and Theorem 3.1(b),
$\derlim\pi_m(P_i)=0$, so by Lemma 3.3, the $\pi_m(P_i)$ are Mittag-Leffler.
Then by Lemma 4.6 below, the $\Gamma_{m+1}(P_i)$ are Mittag-Leffler and therefore
$\derlim\Gamma_{m+1}(P_i)=0$.
Thus $\Gamma_n'(X)=0$ for $n\le m$, and the assertion follows. \qed
\enddemo

\proclaim{Lemma 4.6} Let $P_i$ be an inverse sequence of $(m-1)$-connected
polyhedra, where $m>1$.
If the $\pi_m(P_i)$ are Mittag-Leffler, then so are the $\Gamma_{m+1}(P_i)$.
\endproclaim

\demo{Proof} Suppose that $\pi_m(P_k)$ maps onto the image of $\pi_m(P_j)$ in
$\pi_m(P_i)$, $k>j$.
By cellular approximation, $\pi_m(P_k^{(m)})$ maps onto the image of
$\pi_m(P_j^{(m)})$ in $\pi_m(P_i^{(m+1)})$.
Pick a homotopy equivalence $f\:\bigvee S^m\to P_j^{(m)}$, and let
$\alpha_1,\dots,\alpha_r\in\pi_m(\bigvee S^m)$ be the basis given by
the factors of the wedge.
Then there exist $\beta_1,\dots,\beta_r\in\pi_m(P_k^{(m)})$ mapping onto
the images of $\alpha_1,\dots,\alpha_r$ in $\pi_m(P_i^{(m+1)})$.
Combining together their representatives yields a map
$g\:\bigvee S^m\to P_k^{(m)}$ such that the two compositions
$h_1\:\bigvee S^m@>g>>P_k^{(m)}@>>>P_i^{(m+1)}$ and
$h_2\:\bigvee S^m@>f>>P_j^{(m)}@>>>P_i^{(m+1)}$ induce the same map on $\pi_m$.
Then $h_1$ and $h_2$ are homotopic.
Since $f$ is invertible up to homotopy, $P_j^{(m)}\to P_i^{(m+1)}$
factors up to homotopy through $P_k^{(m)}\to P_i^{(m+1)}$.
Hence $\pi_{m+1}(P_k^{(m)})$ maps onto the image of
$\pi_{m+1}(P_j^{(m)})$ in $\pi_{m+1}(P_i^{(m+1)})$.
Therefore $\Gamma_{m+1}(P_k)$ maps onto the image of $\Gamma_{m+1}(P_j)$ in
$\Gamma_{m+1}(P_i)$. \qed
\enddemo

\remark{Remark} The following versions of Theorem 4.4 are found in
the literature.

Christie \cite{Ch; 5.12} proved that if $X$ is UV$_{m-1}$, where $m>1$,
the Hurewicz homomorphism $\check\pi_m(X)\to\check H_m(X)$ is an isomorphism.
This was rediscovered by K. Kuperberg (1972) and partially by T. Porter (1973).
Artin and Mazur \cite{AM; 4.5} prove a pro-group version of this result,
which additionally implies that under the same hypothesis,
$\pi_{m-1}(X)\to H_{m-1}(X)$ is also an isomorphism; an easier proof is found
in \cite{DS1}.
Both these results follow immediately from Proposition 3.5(a) and the classical
Hurewicz Theorem.

Kodama and Koyama \cite{KK}, \cite{Ko4} prove that if $X$ is UV$_{m-1}$,
where $m>2$, the Hurewicz homomorphism $\pi_m(X)\to H_m(X)$ is an isomorphism.
This result follows by the proof of Theorem 4.4, if instead of Lemma 4.6 one uses
that each $\Gamma_{m+1}(P_i)$ is finite for $m>2$.
(The latter holds since $\pi_{m+1}(\vee S^m)$ is finite for $m>2$, which in turn is
easy to see using the Pontryagin--Thom construction.)
As noted in \cite{KK}, this result does not extend to the case $m=2$, by
considering the suspension of the $p$-adic solenoid for some $p$.

Finally, as noted by Kuperberg \cite{Ku2}, $\check\pi_n(X)=0$ for all $n<3$
does not imply $\check\pi_3(X)\simeq\check H_3(X)$.
Indeed, consider the double suspension of a $p$-adic solenoid with $p$ odd.
The assertion follows since an odd degree self-map of $S^3$ induces
the identity on $\pi_4(S^3)\simeq\Z/2$, which in turn is easy to see using
the Pontryagin--Thom construction.
\endremark

\proclaim{Proposition 4.7}
If $\pi_0(X)=0$, the Hurewicz homomorphism $\check\pi_1(X)\to\check H_1(X)$
is surjective with kernel equal to the closure of the commutator subgroup.
\endproclaim

The first assertion was proved by Marde\v si\'c--Ungar, while stated by Dydak
\cite{D1}.
It does not generalize to connected compacta $X$ \cite{D1}.
If $E$ is the Hawaiian earring (see Example 5.6), the commutator subgroup of
$\check\pi_1(E)$ is not closed in $\check\pi_1(E)$ \cite{D1}.
An erroneous version of Proposition 4.7 which did not address these
two counterexamples to it, was stated without proof in \cite{Ch}.

\demo{Proof} Suppose that $X$ is the inverse limit of polyhedra $P_i$.
Let us write $\pi_1(P_i)=G_i$.
Each $G_i\to H_1(P_i)$ is onto with kernel $G_i'$.
By Theorem 3.1(b) and Lemma 3.3, the $G_i$ are Mittag-Leffler.
Then the commutator subgroup of the stable image of $G_j$, $j>i$, in $G_i$
is the stable image of $G_j'$, $j>i$, hence the $G_i'$ are also Mittag-Leffler.
Hence $f\:\invlim G_i\to\invlim H_1(P_i)$ is onto.
Since $f$ is a continuous homomorphism to an abelian group, its kernel
contains the closure of the commutator subgroup of $\invlim G_i$.
Conversely, if a thread $(g_0,g_1,\dots)\in\ker f$, then each $g_i\in G_i'$,
moreover $g_i$ lies in the stable image of $G_j$, $j>i$.
Then $g_i=[h_{i1},h_{i1'}]\dots[h_{ir},h_{ir'}]$, where each $h_{ik}$ lies
in the stable image of $G_j$.
Let $\hat g_i=[h_1,h_1']\dots[h_r,h_r']$ for some preimages $h_k$ of $h_{ik}$
in $\invlim G_j$, then $\hat g_0,\hat g_1,\dots$ converges to the thread
$(g_0,g_1,\dots)$. \qed
\enddemo

As a final remark, if $X$ is a connected compactum, $\pi_0(X)\to\tl H_0(X)$ is
an epimorphism by the classical Hurewicz Theorem and Theorem 3.1(b).

\head 5. Some examples \endhead

The purpose of this section is twofold: to illustrate the considerations
of the preceding sections and to motivate those of the subsequent ones.

\example{Example 5.1 (null-sequence)} Let $\N^+$ be the one-point
compactification of the countable discrete space $\N$.
Theorem 4.1(iii) implies
$H_0(\N^+)\simeq\invlim H_0(\{1,\dots,i\})\simeq\prod\Z$ (countable product).
On the other hand, since a singular $0$-cycle is supported by only finitely many
points, $\stau\:\sH_0(\N^+)\to H_0(\N^+)$ is an injection onto the subgroup
$\left<(1,1,\dots)\right>\oplus(\bigoplus\Z)$ (countable sum).
At the same time, it is easy to see that both $\pi_0(\N^+)$ and $\spi_0(\N^+)$
are homeomorphic to $\N^+$.
\endexample

\example{Example 5.2 (Cantor set)} Let $\Z_p=\invlim(\dots\to\Z/p^2\to\Z/p)$ be
the topological group of $p$-adic integers, where each $\Z/p^n$ is endowed
with the discrete topology.
We have $\pi_0(\Z_p)\cong\invlim\pi_0(\Z/p^i)\cong\Z_p$ (as spaces).
Clearly $\spi_0(\Z_p)@>\sstau>>\pi_0(\Z_p)$ is a homeomorphism.
On the other hand, it is easy to see that $\sH_0(\Z_p)$ is isomorphic, as
a topological group, to the additive group of the group ring $\Z[\Z_p]$.
At the same time, a single element of $H_0(\Z_p)\simeq\invlim H_0(\Z/p^i)
\simeq\invlim\Z[\Z/p^i]$ may well ``involve'' all points of $\Z_p$.
For instance, the elements $\sum_{i=0}^{n-1} m_it^i\in\Z[\Z/2^n]$ where
each $m_{\alpha_0+2\alpha_1+\dots+2^{n-1}\alpha_{n-1}}=(-1)^{\alpha_0}+
(-1)^{\alpha_1}\cdot 2+\dots+(-1)^{\alpha_n}\cdot 2^{n-1}\ne 0$ (here each
$\alpha_i=0$ or $1$) form a thread in $\invlim\Z[\Z/2^n]\simeq H_0(\Z_2)$.
\endexample

\example{Example 5.3 ($p$-adic solenoid)}
Let $\Sigma_p$ be the $p$-adic solenoid, i.e.\ the mapping torus of
the homeomorphism $h\:\Z_p\to\Z_p$, $h(a)=a+1$.
From the exact sequence of the pair $(\Sigma_p,\Z_p)$ we get
$0\to\Z@>\partial>>\Z_p\to\spi_0(\Sigma_p)\to 0$, whence
$\spi_0(\Sigma_p)\cong\Z_p/\Z$ (as spaces).
Since the same argument also applies to $\pi_0(\Sigma_p)$ (see Theorem 3.1(a)),
$\spi_0(\Sigma_p)@>\sstau>>\pi_0(\Sigma_p)$ is a bijection.
The Steenrod group $\pi_0(\Sigma_p)$ can be alternatively computed using
the Milnor exact sequence of Theorem 3.1(b): $\invlim\pi_0(S^1)=pt$, so
$\pi_0(\Sigma_p)\simeq\derlim\pi_1(S^1)\simeq\Z_p/\Z$ (see Example 3.2).

As for the homology, $\sH_0(\Sigma_p)@>\sstau>>H_0(\Sigma_p)$ has
a nontrivial kernel, contrary to an assertion in \cite{Fe3; p.\ 152}.
We have $\tl H_0(\Sigma_p)=\derlim(\dots@>p>>\Z@>p>>\Z)\simeq\Z_p/\Z$.
Since the composition $\check H_0(pt)\simeq H_0(pt)\to H_0(\Sigma_p)\to
\check H_0(\Sigma_p)\simeq\check H_0(pt)$ is the identity,
$\check H_0(\Sigma_p)$ splits off as a direct summand in $H_0(\Sigma_p)$.
Hence $H_0(\Sigma_p)\simeq\Z\oplus(\Z_p/\Z)$.
On the other hand, from the exact sequence of the pair $(\Sigma_p,\Z_p)$
we obtain $\sH_0(\Sigma_p)\simeq\Z[\Z_p/\Z]$.
(A computation of $H_0(\Sigma_p)$ from the pair $(\Sigma_p,\Z_p)$ is also
possible, but far less trivial; see \cite{M2; Appendix}.)
\endexample

\remark{Remark}
In fact, it is well-known that if $X$ is an inverse limit of fibrations (or
even Steenrod fibrations) between polyhedra, $\spi_n(X)@>\sstau>>\pi_n(X)$
is a bijection for all $n$ \cite{Vo}, \cite{BK; IX.3.1}, \cite{Fe2; 5.5(i)}.
\endremark

\example{Example 5.4 (cluster of solenoids)}
Let $X$ be the cluster of solenoids, i.e.\ the one-point compactification of
the countable disjoint union of copies of $\Sigma_p\but pt$, with compactifying
point $b$.
We claim that $\stau\:\sH_0(X)\to H_0(X)$ is not surjective.
By the Cluster Axiom (see Theorem 4.1), $H_0(X,b)$ is naturally isomorphic to
the countable product of copies of $\tl H_0(\Sigma_p)\simeq\Z_p/\Z$.
A singular $0$-cycle has support in some finite sub-wedge
$\Sigma_p\vee\dots\vee\Sigma_p$ of the cluster.
Hence $\sH_0(X,b)$ maps into (in fact, onto) the countable sum
$\bigoplus\Z_p/\Z\subset\prod\Z_p/\Z$.
\endexample

\example{Example 5.5 (wedge of solenoids)} We claim that
$\stau\:\spi_0(X)\to\pi_0(X)$ is not surjective, where $X=\Sigma_p\vee\Sigma_p$.
Since the Hurewicz homomorphism $h\:\pi_0(X)\to\tl H_0(X)$ is surjective
(see the concluding remark in \S4), it sends some element of $\pi_0(X)$ onto
$(\alpha+\Z,\,\alpha+\Z)\in(\Z_p/\Z)\oplus(\Z_p/\Z)$ for some $\alpha\in\Z_p\but\Z$.
On the other hand, a representative of any element of $\spi_0(X,b)$, where $b$
is the basepoint of the wedge, has support in one of the two factors of the wedge.
Hence $\spi_0(X,b)$ maps into (in fact, onto) the subset
$[(\Z_p/\Z)\oplus 0]\cup [0\oplus(\Z_p/\Z)]$ of $H_0(X,b)$.
\endexample

\remark{Remark} Every $1$-dimensional Steenrod connected compactum has the shape
of either the Hawaiian earring or a finite wedge of circles or a point
(A. Trybulec; see \cite{DS1; 7.3.3}).
\endremark

\example{Example 5.6 (Hawaiian earring)}
Let $E$ be the Hawaiian earring, that is the one-point compactification of
$\R\x\N$.
Then $\pi_1(E)\simeq\invlim\pi_1(\bigvee_{i=1}^r S^1)\simeq\invlim F_r$, where
$F_r$ denotes the free group $\left<x_1,\dots,x_r\mid\,\right>$.
On the other hand, if $f\:(S^1,pt)\to (E,x)$ is a map, where $x$ is the
compactifying point, we may assume that $f$ is transverse to each
$x_n:=(0,n)\in\R^n\x\N$.
The preimage of all these points is a framed $0$-submanifold of $S^1\but x$
with components indexed by the letters $x_1^{\pm1},x_2^{\pm1}\dots$ so that
each $x_i$ occurs only finitely many times.
Similarly to the usual Pontryagin--Thom construction, it is easily verified that
$\spi_1(E)$ is isomorphic to the bordism group of such framed manifolds.
It follows that $\spi_1(E)$ is isomorphic to the subgroup $L$ of
$\invlim F_r$ consisting of all threads $(\dots,g_2,g_1)$ such that when
$g_i$'s are written as reduced words $w_i$ in their free groups, each $x_j$
occurs only boundedly many times (i.e.\ at most $n(j)$ times for some function
$n(j)$) in each $w_i$.
Thus $\spi_1(E)@>\sstau>>\pi_1(E)$ has nontrivial cokernel.

In contrast, the abelianization $\sH_1(E)$ of $L$ clearly maps onto
$H_1(E)\simeq\prod_{i=1}^\infty\Z$.
It turns out that the kernel of this surjection splits off as a direct summand
in $\sH_1(E)$ and is isomorphic to $\bigoplus_{\frak c}\Q\oplus\bigoplus_p A_p$,
where $\frak c$ is the cardinality of the continuum, $p$ runs over all primes,
and $A_p$ is the $p$-adic completion of $\bigoplus_{\frak c}\Z_p$ \cite{EK2}.
This is found by a direct computation of the abelianization of $L$.
\endexample

\example{Example 5.7 (Hawaiian snail)}
Let $E^n$ be the $n$-dimensional Hawaiian earring, now viewed as
the one-point compactification of $\R^n\x\Z$.
Let $\sigma_n\:E^n\to E^n$ be the homeomorphism given by the shift $\Z\to\Z$,
$i\mapsto i+1$.
Let $X$ be the mapping torus of $\sigma_n$, where $n\ge 2$.
From the Milnor exact sequence it is easy to see that $\pi_n(X)\simeq H_n(X)=0$.
On the other hand, if $\tl X$ is the universal cover of $X$, it is easy to
see that $\spi_n(\tl X)$ contains $\Z[\Z]$.
Hence $\spi_n(X)$ does not inject into $\pi_n(X)$.

To compute $\sH_n(X)$, we first note that by the Hurewicz Theorem (which is
proved in the required generality in \cite{Sp}, for instance),
$\sH_n(E^n)\simeq\spi_n(E^n)$ for $n>1$;
and by the Pontryagin--Thom construction (see Example 5.6 or the proof of
Theorem 1.1; compare \cite{EK3}), $\spi_n(E^n)\simeq\prod_{i=1}^\infty\Z$ for
$n>1$.
Then the exact sequence of the pair $(X,E^n)$ yields
$\prod_{i=1}^\infty\Z@>\partial>>\prod_{i=1}^\infty\Z\to\sH_n(X)\to 0$,
where the boundary map $\partial$ is given by
$1-t\:\Z[[t^{\pm 1}]]\to\Z[[t^{\pm 1}]]$.
But this is an isomorphism (its inverse is given by $1+t+t^2+\dots$), so
$\sH_n(X)=0$.
\endexample

\example{Example 5.8 (Hawaiian snail with partitions)}
Let $X$ be the mapping torus of $\sigma_{n-1}$ from the previous example,
where $n\ge 3$.
Consider an LC$_{n-1}$ compactum $Y$ obtained from $X$ by attaching
small $n$-disks $D_1,D_2,\dots$.
We will show that $\check\tau\stau\:\spi_n(Y)\to\check\pi_n(Y)$ is not
surjective.
(Note that Theorem 6.5 below implies that
$\check\tau\stau\:\spi_n(Y^+)\to\check\pi_n(Y^+)$ is surjective, where $Y^+$
is obtained by attaching a $2$-disk to $Y$ along the orbit of
the compactifying point.)

Clearly $H_n(X)\simeq\check H_n(X)\simeq\Z$ and
$\pi_n(X)\simeq\check\pi_n(X)\simeq\Z[\Z]$.
Pick some $a\in\pi_n(X)$ whose image in $H_n(X)$ is non-trivial, and let
$b\in\pi_n(Y)$ be the image of $a$.
Suppose that $b=\stau([f])$ for some spheroid $f\:S^n\to Y$.
Then $f$ lifts to the universal cover $Y_\infty$.
Since $S^n$ is compact, there exists a finite number $r$ such that the image of
the lift of $f$ is contained in a union of $r$ translates of a fundamental
domain $F_0$ of the action of $\Z$ on $Y_\infty$.
Let $f_{r+1}\:S^n\to Y_{r+1}$ be a lift of $f$ into the $(r+1)$-fold cyclic
cover of $Y$.
Then $f_{r+1}(S^n)$ is disjoint from some translate $F_i$ of
the fundamental domain of $\Z/(r+1)$.

By Theorem 3.15(b), $\pi_n(X_{r+1})\to\pi_n(X)$ is an isomorphism, where $X_{r+1}$
is the $(r+1)$-fold cyclic cover of $X$.
So $a$ is the image of some $a_{r+1}\in\pi_n(X_{r+1})$.
The image $b_{r+1}\in\pi_n(Y_{r+1})$ of $a_{r+1}$ is represented by a translate
of $f_{r+1}$, so $b_{r+1}$ maps trivially to
$\pi_n(Y_{r+1},\,\Cl(Y_{r+1}\but F_j))$ for some copy $F_j$ of the fundamental
domain.
In particular, it maps trivially to $H_n(Y_{r+1},\,\Cl(Y_{r+1}\but F_j))$.
Let $F'_j=F_j\cap X_{r+1}$ be the corresponding fundamental domain of
the covering map $X_{r+1}\to X$.
Since $Y_{r+1}$ is $n$-dimensional,
$H_n(X_{r+1},\,\Cl(X_{r+1}\but F'_j))\to H_n(Y_{r+1},\,\Cl(Y_{r+1}\but F_j))$
is injective.
Since $\Cl(X_{r+1}\but F'_j)$ deformation retracts onto a copy of the
$(n-1)$-dimensional compactum $E^{n-1}$,
$H_n(X_{r+1})\to H_n(X_{r+1},\,\Cl(X_{r+1}\but F'_j))$ is injective.
Thus $a_{r+1}$ maps trivially to $H_n(X_{r+1})$.
But then $a$ maps trivially to $H_n(X)$, which contradicts our choice of $a$.
\endexample

\remark{Remark} Another LC$_{n-1}$ compactum $Y$ such that
$\check\tau\stau\:\spi_n(Y)\to\check\pi_n(Y)$ is not surjective is
$Y=S^1\vee E^n$, where the basepoint of the Hawaiian earring $E^n$ is taken
at its compactifying point.
This was observed by Dydak \cite{D2; Example 8.10} and Zdravkovska \cite{Zd},
of which the present author was unaware when writing up Example 5.8.
Zdravkovska also shows that $\spi_n(X)\to\pi_n(X)$ is not surjective for any
compactum $X$, Steenrod homotopy equivalent to $Y$.
In contrast, Ferry proved that every UV$_1$ compactum is Steenrod homotopy
equivalent to a compactum $X$ such that
$\stau\:[Z,X]^\tr\to [Z,X]$ is a bijection for every finite-dimensional
compactum $Z$ \cite{Fe2}.
His construction always produces an infinite dimensional $X$.
\endremark

\head 6. Comparison of theories \endhead

We recall that a compactum $X$ is {\it locally $n$-connected (L$\sitC_n$)}, if
every neighborhood $U$ of every $x\in X$ contains a neighborhood $V$ of $x$ such
that $\spi_i(V)\to\spi_i(U)$ is trivial for all $i\le n$.
We say that $X$ is {\it semi-L$\sitC_n$} if it is L$\sC_{n-1}$ and every $x\in X$
has a neighborhood $V$ such that $\spi_n(V)\to\spi_n(X)$ is trivial.

The singular group $\spi_n(X)$ can be topologized similarly to
$\pi_n(X)$, that is by declaring a subset to be open if it is the preimage
of an open subset of $\check\pi_n(X)$; or equivalently if it is
a point-inverse of the homomorphism induced by some map of $X$ into
a polyhedron.
Note that $\spi_n(X)$ and $\pi_n(X)$ are in general non-Hausdorff, whereas
$\check\pi_n(X)$ is metrizable.

\proclaim{Theorem 6.1} Fix some $n\ge 0$ and let $X$ be a compactum.

\medskip
(a) The following conditions are equivalent:

\smallskip
\noindent
{\rm L$\sC_n$:} $X$ is locally $n$-connected with respect
to $\spi_i(\cdot)$;

\noindent
{\rm LC$_n$:} $X$ is locally $n$-connected with respect to
$\pi_i(\cdot)$;

\noindent
{\rm L\v C$_n$:} $X$ is locally $n$-connected with respect to
$\check\pi_i(\cdot)$.

\medskip
(b) If $X$ is LC$_1$ or $n=0$, the following conditions are also equivalent to
the preceding ones:

\smallskip
\noindent
{\rm HL$\sC_n$:} $X$ is locally $n$-connected with respect
to $\sH_i(\cdot)$;

\noindent
{\rm HLC$_n$:} $X$ is locally $n$-connected with respect to
$H_i(\cdot)$;

\noindent
{\rm HL\v C$_n$:} $X$ is locally $n$-connected with respect to
$\check H_i(\cdot)$.

\smallskip
(c) If $X$ is L$\sitC_n$, then $\spi_n(X)@>\sstau>>\pi_n(X)$
and $\pi_n(X)@>\check\tau>>\check\pi_n(X)$ are bijective.

\smallskip
(d) If $X$ is semi-L$\sitC_n$, then $\spi_n(X)@>\check\tau\sstau>>\check\pi_n(X)$
is bijective.

\smallskip
(e) If $X$ is L$\sitC_{n-1}$, then $\spi_n(X)@>\check\tau\sstau>>\check\pi_n(X)$
is onto a dense subset.

\bigskip
(f) $X$ is semi-L$\sitC_n$ iff $X$ is L$\sitC_{n-1}$ and $\spi_n(X)$
is discrete.

(g) $X$ is semi-LC$_n$ and $\spi_n(X)@>\sstau>>\pi_n(X)$ is onto
iff $X$ is LC$_{n-1}$ and $\pi_n(X)$ is discrete.

\medskip
(h) $X$ is semi-L\v C$_n$ iff $X$ is L\v C$_{n-1}$ and
$\check\pi_n(X)$ is discrete.
\endproclaim

Most of Theorem 6.1 is contained in the literature.
Hurewicz proved that $X$ is L$\sC_1$ and HL\v C$_n$ if and only if it is
L$\sC_n$ \cite{Hu}.
Borsuk proved an equivalent of the second assertion of (c) \cite{Bo}.
Part (d) was proved by Kuperberg \cite{Ku1} and under the L$\sC_n$ hypothesis also
by Bogatyi \cite{Bog}.
Dydak proved part (e) and that if $X$ is L$\sC_n$, then $\check\pi_n(X)$ is
discrete \cite{D2; Theorem 8.7}.
He obtains these along with both assertions of (c) by a modification of
the proof of his pro-group version of the Smale theorem (see Theorem 3.15(e));
this local-global approach, similar in spirit to the Zeeman spectral sequence,
provides a substantial alternative to the methods used in the present section.

For clarity we suppress all basepoints in the proof and leave it
to the reader to verify that this does not cause confusion.
The proofs of all parts depend on

\demo{Ferry's Construction \cite{Fe2; p.\ 381}}
We represent $X$ as the limit of an inverse sequence $P$ of polyhedra
and PL maps, and use the notation introduced in the beginning of \S2.
Let us triangulate $P_{[0,\infty)}$ so that the diameters of simplices tend to
zero as they approach $X$.
Assuming that $X$ is L$\sC_n$, it follows by induction on $n$ that for each
$\eps>0$ there exists a $k$ such that
$Y:=X\cup P_{[k,\infty)}^{(n+1)}$ retracts onto $X$ by a retraction
$r_k^{(n+1)}$ that is $\eps$-close to the identity.
Since $r_k^{(n+1)}$ is close to the identity and $P_{[0,\infty]}$ is
L$\sC_\infty$, for each $l$ there exists a $k$ such that $r_k^{(n+1)}$ is
homotopic to the identity by a homotopy $Y\x I\to P_{[l,\infty]}$ keeping
$X$ fixed.
If $X$ is only semi-L$\sC_n$, the restriction of $r_k^{(n)}$ to
$P_{[k,\infty)}^{(n)}$ still extends to a map $P_{[k,\infty)}^{(n+1)}\to X$,
which however is no longer close to the identity on top-dimensional simplices.
\enddemo

\demo{Proof of Theorem 6.1. (e)} Let us represent the given
$\alpha\in\check\pi_n(X)$ by a Steenrod spheroid
$\Phi\:S^n\x [0,\infty)\to P_{[0,\infty)}$.
We may assume that each $\phi_i=\Phi|_{S^n\x\{i\}}\:S^n\to P_i$ goes into
the $n$-skeleton.
If $X$ is L$\sC_{n-1}$, for each $l$ there exists a $k$ such that $r_k^{(n)}$ is
homotopic to the identity within $P_{[l,\infty]}$, and therefore
$\psi_k:=r_k^{(n)}\phi_k\:S^n\to X$ is homotopic to $\phi_k$ within
$P_{[l,\infty]}$.
The strong deformation retraction of $P_{[l,\infty]}$ onto $P_l$ carries this
homotopy onto a homotopy between $p^\infty_l\psi_k$ and $p^k_l\phi_k$.
It also carries $\Phi|_{S^n\x[l,k]}$ onto a homotopy between $p^k_l\phi_k$
and $\phi_l$.
Hence $[\phi_l]\in\pi_n(P_l)$ is the image of $[\psi_k]\in\spi_n(X)$.
It follows that $\check\tau\stau([\psi_k])\in\check\pi_n(X)$ converge to
$\alpha$. \qed
\enddemo

\demo{(d): Surjectivity} We use the notation in the proof of (e).
If $X$ is semi-L$\sC_n$, each $\psi_k$ is homotopic
to $\psi_{k+1}$ via $r_k^{(n+1)}\Phi|_{S^n\x[k,k+1]}$.
So the $\check\tau\stau([\psi_k])$ are all equal.
Since $\check\pi_n(X)$ is Hausdorff, they also equal their limit $\alpha$. \qed
\enddemo

\demo{(d): Injectivity \cite{Fe2; p.\ 381}} We will prove that
$\spi_n(X)@>\check\tau\sstau>>\check\pi_n(X)@>p^\infty_k>>\pi_n(P_k)$
is injective.
Given a singular spheroid $\phi\:S^n\to X$ in the kernel of this composition,
it bounds a null-homotopy $\hat\phi\:B^{n+1}\to P_{[l,\infty]}$.
By the cellular approximation theorem, we may assume that its image is
in the $(n+1)$-skeleton (apart from $\partial B^{n+1}$, which is mapped
into $X$).
Then $r_k^{(n+1)}\hat\phi$ is a null-homotopy of $\phi$. \qed
\enddemo

\demo{(c)} This follows either from (d), Theorem 3.1(b) and part (1) of
the following lemma, or from the injectivity in (d), Theorem 3.1(b) and
part (2) of the following lemma.
\enddemo

\proclaim{Lemma 6.2} Fix some $n\ge 0$ and let $X$ be a compactum.

(1) {\rm (Borsuk \cite{Bo})} If $X$ is L$\sitC_n$, then
$\pi_n(X)@>\check\tau>>\check\pi_n(X)$ is injective.

(2) If $X$ is L$\sitC_n$, then $\spi_n(X)@>\sstau>>\pi_n(X)$
is surjective.

(3) {\rm (Hurewicz \cite{Hu})} If $X$ is L$\sitC_{n-1}$
and L\v C$_n$, then $\spi_n(X)@>\check\tau\sstau>>\check\pi_n(X)$ is
injective.
\endproclaim

\demo{Proof. (1) \cite{Fe2; p.\ 381}} Given a map $\phi\:S^{n+1}\to P_k$, its
image may be assumed to lie in the $(n+1)$-skeleton, and then
$r_k^{(n+1)}\phi\:S^{n+1}\to X$ is homotopic to $\phi$ with values in
$P_{[l,\infty)}$.
Thus $\spi_{n+1}(X)$ maps onto the image of $\pi_{n+1}(P_k)$ in
$\pi_{n+1}(P_l)$.
It follows that the $\pi_{n+1}(P_i)$ are Mittag-Leffler, and the assertion follows
from Theorem 3.1(b). \qed
\enddemo

\demo{(2)} Given a Steenrod spheroid $\chi\:S^n\x[0,\infty)\to P_{[0,\infty)}$
with image in the $(n+1)$-skeleton, each $t\in [k,\infty)$ yields a singular
spheroid $\phi_t=r_k^{(n+1)}\chi|_{S^n\x\{t\}}\:S^n\to X$.
The image of $[\phi_k]$ in $\pi_n(X)$ is represented by
$\Phi_k^+\:S^n\x[0,\infty)\to P_{[0,\infty)}$ defined by
$\Phi_k^+(x,t)=\Pi_t(\phi_k(x))$ (the homotopy $\Pi_t$ is defined in the
beginning of \S2).
The family of maps $\Phi_t\:S^n\x [t,\infty)\to P_{[t,\infty)}$, defined by
$\Phi_t(x,s)=\Pi_s(\phi_t(x))$, yields a proper homotopy between
$\Phi_k$ and the map $\Psi\:S^n\x [k,\infty)\to P_{[k,\infty)}$ defined by
$\Psi(x,t)=\Phi_t(x,t)$.
On the other hand, the same family $\Phi_t$ yields a homotopy between $\Psi$
and $r_k^{(n+1)}\chi|_{S^n\x [k,\infty)}$.
The latter is homotopic to $\chi|_{S^n\x[k,\infty)}$ since $r_k^{(n+1)}$ is
homotopic to the identity.
The resulting homotopy $h\:S^n\x[k,\infty)\x I\to P_{[0,\infty]}$
between $\chi|_{S^n\x[k,\infty)}$ and $\Psi$ is such that $h^{-1}(C)$ is
compact for every compact $C\i P_{[0,\infty)}$.
Therefore it can be pushed in so as to become a proper homotopy
within $P_{[k,\infty)}$.
The resulting proper homotopy between $\chi|_{S^n\x[k,\infty)}$ and $\Phi_k$
can be extended to a proper homotopy between $\chi$ and $\Phi_k^+$ using
Borsuk's lemma. \qed
\enddemo

\demo{(3) \cite{DS3; \S3}}
Let $\phi\:S^n\to X$ be a singular spheroid that is \v Cech trivial.
Arguing as in the proof of the injectivity in 6.1(d), we obtain that it bounds
a null-homotopy $\Phi\:B^{n+1}\to P_{[k,\infty]}$ for any $k$ given
in advance.
In particular, given an $\eps>0$, we may choose $k$ so that $r_k^{(n)}$ is
$\eps$-close to the identity and the simplices of the triangulation of
$P_{[k,\infty)}$ have diameters $<\eps$.
By cellular approximation, we may assume that $\Phi$ sends
the $n$-skeleton $T^{(n)}$ of some triangulation $T$ of the interior of
$B^{n+1}$ into the $n$-skeleton of $P_{[k,\infty)}$, and moreover that
$\Phi(\partial\sigma)$ has diameter $<3\eps$ for every $(n+1)$-simplex
$\sigma$ of $T$.
Then the image of the restriction $\phi_\sigma\:\partial\sigma\to X$ of
$r_k^{(n)}\Phi\:T^{(n)}\to X$ has diameter $<5\eps$.
Moreover, since $P_{[k,\infty]}$ is L$\sC_\infty$, given any $\delta>0$, we may
choose $\eps$ so that each $\phi_\sigma$ is $\delta$-null-homotopic in
$P_{[k,\infty]}$.
Since $X$ is L\v C$_n$, $\phi_\sigma$ bounds a $\delta$-null-homotopy
$\Phi_\sigma\:\sigma\to P_{[k_1,\infty]}$ for any $k_1$ given in advance.

We may now repeat the entire process, with a smaller $\eps$, separately for each
$\phi_\sigma\:S^n\to X$, where $\sigma$ runs over all $(n+1)$-simplices of $T$.
The construction converges to a map $f$ into $X$ defined on the union $U$ of
the $n$-skeleta $T_i^{(n)}$ of a sequence of successive subdivisions $T_i$ of
$T$ such that each $(n+1)$-simplex $\Delta$ of $T_i$ has diameter $<2^{-i}$
and $f(\partial\Delta)$ has diameter $<2^{-i}$.
Since $U$ is a dense subset of $B^{n+1}$, this also defines $f$ on the
whole of $B^{n+1}$, and this extension is easily seen to be well-defined and
continuous. \qed
\enddemo

\demo{Proof of Theorem 6.1, continued. (a)} If $X$ is L$\sC_n$, then by
a local version%
\footnote{By ``a local version of $\frak X$'', we mean ``an appropriate
straightforward generalization of $\frak X$, applied to closed neighborhoods
of points''.
For instance, right now we are dealing with the following generalization of
Lemma 6.2(2): if $Z\i Y\i X$ are compacta, where $X$ is L$\ssC_n$ and $Y$ is
a neighborhood of $Z$ in $X$, then the image of $\stau\:\spi_n(Y)\to\pi_n(Y)$
contains the image of $i_*\:\pi_n(Z)\to\pi_n(Y)$.
In some cases, $\frak X$ can be applied directly by virtue of Theorem 6.11.}
of Lemma 6.2 (2), $X$ is LC$_n$.
(Note that an alternative proof of (2) can be obtained from (d) and (1).)
If $X$ is LC$_n$, then it is L\v C$_n$ since $\check\tau$ is
always an epimorphism.
Finally, to prove that L\v C$_n$ implies L$\sC_n$, we may assume, arguing
by induction, that $X$ is L$\sC_{n-1}$.
Then the assertion follows from a local version of Lemma 6.2 (3).
\qed
\enddemo

\demo{(h), \imp} If an element of $\pi_n(X)$ maps trivially to
$\pi_n(P_k)$, it can be represented by a level-preserving Steenrod spheroid
$\chi\:S^n\x[0,\infty)\to P_{[0,\infty)}$ whose restriction to $S^n\x\{k\}$
bounds a null-homotopy $\phi\:D^{n+1}\to P_k$.
Let $\psi\:\R^{n+1}\to P_{[k,\infty)}$ be the proper map obtained by combining
$\phi$ with the restriction of $\chi$ to $S^n\x [k,\infty)$.
We may assume that $\psi$ is cellular with respect to some triangulation of
$\R^{n+1}$.
Since $X$ is L$\sC_{n-1}$, the retraction $r^{(n)}_k\:Y\to X$ is homotopic to
the identity by a homotopy $h\:Y\x I\to P_{[l,\infty]}$ keeping $X$ fixed.
If $\Delta$ is an $(n+1)$-simplex of $\R^{n+1}$, let
$\eta\:(\partial\Delta)\x I\to P_{[l,\infty]}$ be the restriction of $h$ to
$\psi(\partial\Delta)\x I$.
By using the homotopy $\Pi_t$ from the beginning of \S2, we may assume that
$\eta$ restricts to a proper map $(\partial\Delta)\x[0,1)\to P_{[l,\infty)}$.
Since $X$ is semi-L\v C$_n$, this Steenrod spheroid is \v Cech trivial
as long as $l$ is large enough.
Since $\Delta$ is arbitrary, it follows that $\chi$ is also \v Cech trivial.
\qed
\enddemo

\demo{(h), \when}
By Lemma 3.4(b), the hypothesis implies that $\check\pi_n(X)$ injects into
$\pi_n(P_k)$ for some $k$.
On the other hand, if $U\i X$ is so small that its image in $P_k$ lies in
the star of some vertex of $P_k$, the composition
$\check\pi_n(U)\to\check\pi_n(X)\to\pi_n(P_k)$ is trivial.
Then also $\check\pi_n(U)\to\check\pi_n(X)$ must be trivial. \qed
\enddemo

\demo{(g), \imp} For convenience of notation, we assume that $n>0$;
the case $n=0$ is similar.
By the hypothesis, every element of $\pi_n(X)$ can be
represented by a spheroid $\phi\:S^n\to X$.
If the composition $S^n@>\phi>>X@>p^\infty_k>>P_k$ is null-homotopic for
some $k$, then $\phi$ bounds a disk $\psi\:D^{n+1}\to P_{[k,\infty)}$.
Since $X$ is L$\sC_{n-1}$, given any $\eps>0$, the first stage of Hurewicz's
Construction (from the proof Lemma 6.1(1)) represents $[\psi]\in\spi_n(X)$
as a finite sum of the classes of compositions
$(S^n,pt)@>p>>(S^n\vee I,0)@>f_i>>(X,x)$ where $p$ collapses a hemisphere of
$S^n$ onto the whisker $I$, sending the basepoint to the end of the whisker,
and each $f_i(S^n)$ is of diameter at most $\eps$.
Since $X$ is semi-LC$_n$, each of these is Steenrod null-homotopic.
Thus $\stau([\phi])=0$.
This proves that the composition
$\pi_n(X)@>\check\tau>>\check\pi_n(X)@>p^\infty_k>>\pi_n(P_k)$ is injective
for some $k$, which implies the assertion. \qed
\enddemo

\demo{(g), \when}
By the hypothesis, $\pi_n(X)@>\check\tau>>\check\pi_n(X)$ is a bijection
and $\check\pi_n(X)$ is discrete.
Hence by (e), $\spi_n(X)@>\sstau>>\pi_n(X)$ is onto, and from (h), $X$ is
semi-LC$_n$. \qed
\enddemo

\demo{(f)} If $X$ is semi-L$\sC_n$, $\check\pi_n(X)$ is discrete by the proof
of the injectivity in (d).
(Alternatively, by a local version of (e), taking into account
that $\check\pi_n(X)$ is Hausdorff, it follows that $X$ is semi-L\v C$_n$,
and the assertion follows from (h).)

Conversely, if $\spi_n(X)$ is discrete, then by the definition of its topology
it injects into $\check\pi_n(X)$; let $G_n$ denote the image of this injection.
If $n=0$, clearly the underlying set of $\check\pi_n(X)$ can be identified
with the set of quasi-components of $X$ (which coincides with the set of
components since $X$ is compact).
Hence $G_0=\check\pi_0(X)$.
If $n>0$, the topological group $\check\pi_n(X)$ can be endowed with
the left-invariant metric defined as the inverse limit of the metrics on
$\pi_n(P_i)$ where the distance between every two points equals $1$.
Since the induced metric on $G_n$ is also left-invariant and its induced
topology is discrete, it follows that its induced uniform structure is
discrete.
So every Cauchy sequence in $G_n$ is eventually constant, thus $G_n$ is
closed in $\check\pi_n(X)$.
Hence by (e), $G_n=\check\pi_n(X)$, and we conclude that $\check\pi_n(X)$
is discrete in both cases ($n=0$ and $n>0$).
Then by (h), $X$ is semi-L\v C$_n$ and consequently also semi-L$\sC_n$. \qed
\enddemo

\demo{(b)} If $X$ is L$\sC_n$, then $\sH_n(X)@>\sstau>>H_n(X)$ and
$H_n(X)@>\check\tau>>\check H_n(X)$ are isomorphisms similarly to the proof
of (c), using pseudo-manifolds instead of $S^n\x [0,\infty)$.
If $X$ is L$\sC_{n-1}$ and HL$\sC_n$, this argument still works: the retraction
$r_k^{(n+1)}$ becomes defined ``up to attaching handles to the top-dimensional
simplices'', that is, it is now defined on $Y':=X\cup Q_k$, where $Q_k$ is
a polyhedron obtained by replacing each top-dimensional simplex $\Delta^{n+1}$
of $P_{[k,\infty)}^{(n+1)}$ by some $(n+1)$-pseudo-manifold with boundary
$\partial\Delta^{n+1}$.
This $r_k^{(n+1)}\:Y'\to X$ is homotopic to the projection $Y'\to Y$ (instead
of $\id_Y$) by a homotopy $Y'\x I\to P_{[l,\infty)}$ keeping $X$ fixed.

Obviously L$\sC_0$ is equivalent to HL$\sC_0$; and L\v C$_0$ is equivalent to
HL\v C$_0$.
By the above, HL$\sC_0$ implies HLC$_0$ which in turn implies
HL\v C$_0$.
Since singular $1$-cycles are all spherical, L$\sC_1$ implies HL$\sC_1$.
Then by the above it also implies HLC$_1$ and HL\v C$_1$.
It remains to consider the case $n>1$.
Arguing by induction, we may assume that $X$ is L$\sC_{n-1}$.
Then by the above, HL$\sC_n$ implies HLC$_n$, which in turn implies
HL\v C$_n$.
On the other hand, L$\sC_n$ is equivalent to HL$\sC_n$ by Lemma 6.3(a) below,
and L\v C$_n$ is equivalent to HL\v C$_n$ by Lemma 6.3(b) below
along with Theorem 6.1(d). \qed
\enddemo

\proclaim{Lemma 6.3} Let $f\:X\to Y$ be a map between compacta such that
for any map $g$ of an $(n-1)$-polyhedron into $X$, the composition $fg$ is
null-homotopic.
Then

\medskip
(a) $\im[\sH_n(X)\to\sH_n(Y)]$ is contained in $\im[\spi_n(Y)\to\sH_n(Y)]$,
and if $n>1$, $\ker[\spi_n(X)\to\sH_n(X)]$ is contained in
$\ker[\spi_n(X)\to\spi_n(Y)]$.

\medskip
(b) If $X$ is L$\sitC_{n-2}$ and $n>1$,
$\im[\check H_n(X)\to\check H_n(Y)]\i\im[\check\pi_n(Y)\to\check H_n(Y)]$ and
$\ker[\check\pi_n(X)\to\check H_n(X)]\i\ker[\check\pi_n(X)\to\check\pi_n(Y)]$.
\endproclaim
\medskip

The second assertion of (a) is essentially the well-known ``Eventual
Hurewicz Theorem'' of Ferry \cite{Fe1; 3.1} (the name seems to originate from
F. Quinn's ``Ends of maps --- I'' and has become standard).
Indeed, the hypothesis of Lemma 6.3 is satisfied if $f$ is a composition
$f_{n-1}\dots f_0$ where each $f_i$ is trivial on $\spi_i$.

It will follow from Theorem 6.7 that the first assertion of (b) also holds
for $n=1$ if $X$ is assumed to be L$\sC_0$.

\demo{Proof. (a)} Given a singular cycle $\phi\:Z\to X$, where $Z$ is
an $n$-pseudo-manifold with a fixed triangulation, by the hypothesis $f\phi$
is homotopic to a $\psi$ that sends $Z^{(n-1)}$ to the basepoint.
Thus $f_*[\phi]=\sum[\psi_i]$, where $\psi_i\:S^n\to Y$ are given by
restricting $\psi$ to the top-dimensional simplices.

Given a singular $(n+1)$-chain $\phi\:Z\to X$ with $\partial Z=S^n$,
by the hypothesis $f\phi$ is homotopic to a $\psi$ that sends $Z^{(n-1)}$
to the basepoint.
This $\psi$ factors though the quotient $Z'=Z/Z^{(n-1)}$, where
$\pi_n(Z')\to H_n(Z')$ is an isomorphism by the classical Hurewicz Theorem.
Thus $S^n=\partial Z$ is null-homotopic in $Z'$ and therefore in $Y$. \qed
\enddemo

\demo{(b)} Suppose that $X$ and $Y$ are the inverse limits of compact
polyhedra $P_i$ and $Q_i$.
For each $l$ there exists a $k>l$ such that
$r^{(n-1)}_k\:X\cup P_{[k,\infty)}^{(n-1)}\to X$ is defined and homotopic with
values in $P_{[l,\infty]}$ to the identity (see Ferry's Construction).
Then the inclusion $\phi\:K\to P_k$ of the $(n-1)$-skeleton $K$ of some
triangulation of $P_k$ is homotopic with values in $P_{[l,\infty]}$ to some
$\psi\:K\to X$.
By the hypothesis, $f\psi\:K\to Y$ is null-homotopic.
By Lemmas 2.1 and 2.5(b$_0$) we may assume (after dropping some of the $P_i$)
that $f$ extends to a level-preserving map $F\:P_{[0,\infty]}\to Q_{[0,\infty]}$.
Then $F\phi$ is null-homotopic with values in $Q_{[l,\infty]}$.
Hence $F|_{P_k}\:P_k\to Q_{[l,\infty]}$ factors up to homotopy through
the $(n-1)$-connected polyhedron $L:=P_k/K$.
By the Hurewicz Theorem, $\pi_n(L)\simeq H_n(L)$.

Consider a Steenrod spheroid $\phi\:S^n\x [0,\infty)\to P_{[0,\infty)}$
such that each $\phi_k=\phi|_{S^n\x k}\:S^n\to P_k$ represents the trivial
element of $H_n(P_k)$.
Then $F\phi_k$ represents the trivial element of $\pi_n(Q_{[l,\infty]})$.
Since $\phi_k$ is homotopic with values in $P_{[l,k]}$ to $\phi_l$,
the homotopy class $[F\phi_l]$ is also trivial in $\pi_n(Q_{[l,\infty]})$.
Since $Q_{[l,\infty]}$ deformation retracts onto $Q_l$, the class $[F\phi_l]$
is trivial in $\pi_n(Q_l)$.

Let $\psi\:M\to P_{[0,\infty)}$ be a Steenrod cycle, and consider its
level $\psi_k=\psi|_{\psi^{-1}(P_k)}$.
Then $F\psi_k$ is homologous with values in $Q_{[l,\infty]}$ to a spherical
cycle $\phi_l\:S^n\to Q_{[l,\infty]}\to Q_l$.
Since each $\psi_k$ is homologous to $\psi_{k+1}$ with values in $P_{[k,\,k+1]}$,
each $\phi_l$ is homotopic to $\phi_{l+1}$ with values in $Q_{[m+1,\infty]}$,
where $l=l(m+1)$ is chosen similarly to $k=k(l)$.
Then $\chi_m:=\Pi_m\phi_l$ and $\chi_{m+1}$ are level-preserving homotopic
with values in $Q_{[m,\,m+1]}$ (the homotopy $\Pi_t$ is defined in \S2).
Thus $\chi_m$ are the levels of a Steenrod spheroid
$S^n\x [0,\infty)\to Q_{[0,\infty)}$.
On the other hand, since $\psi_m$ is homologous to $\psi_k$ with values in
$P_{[m,k]}$, and $Q_{[m,\infty]}$ deformation retracts onto $Q_m$,
the composition $F\psi_m$ is homologous with values in $Q_m$ to $\chi_m$. \qed
\enddemo

\definition{Spheroids with trunks} An {\it $n$-spheroid with trunks}
in a compactum $X$ is a uniformly continuous map
$\phi\:(\partial R,pt)\to (X,x_0)$, where the uniform space $R$ is a regular
$\eps$-neighborhood of a properly embedded infinite binary tree in $\R^{n+1}$,
with respect to some proper function $\eps\:\R^{n+1}\to (0,1]$.
(The word ``trunk'' is used here in the meaning ``elephant's trunk''.)
Note that $\partial R$ is non-uniformly homeomorphic to the complement of
a tame Cantor set in $S^n$.

If $X$ is the limit of an inverse sequence $P$ of polyhedra and PL maps,
similarly to the proof of Lemma 2.1, every spheroid with trunks
$\phi\:\partial R\to X$ extends to a uniformly continuous map
$\bar\Phi\:R\to P_{[0,\infty]}$ that restricts to a proper map
$\Phi\:\Int R\to P_{[0,\infty)}$, and any two such extensions are homotopic
through such extensions.
We shall say that $\phi$ {\it represents} $[\Phi]\in\pi_n(X)$.
\enddefinition

\proclaim{Proposition 6.4} If $X$ is an LC$_n$ compactum and $n\ge 0$, every
element of $\pi_{n+1}(X)$ is representable by an $(n+1)$-spheroid with trunks.
\endproclaim

\remark{Remark} Proposition 6.4 can be generalized as follows.
If $X$ is an LC$_n$ compactum, every element of $\pi_{n+d}(X)$ can be represented
by a uniformly continuous map $\partial R\to X$,
where $R$ is a regular $\eps$-neighborhood in $\R^{n+d+1}$ of a properly
embedded mapping telescope of an inverse sequence of $(d-1)$-polyhedra $P_i$
(starting with $P_0=pt$), with respect to some proper function
$\eps\:\R^{n+d+1}\to (0,1]$.
In fact, $P_i$ can be more specifically described as the dual skeleta of
the successive barycentric subdivisions of some triangulation of $S^{n+d}$,
with the obvious bonding maps.
In this case, if $2d-1\le n+d$, then $\invlim P_i$ is homeomorphic to
the Menger cube $\mu_{d-1}$ by Bestvina's characterization of $\mu_{d-1}$
(see \cite{GHW}).
Thus $\partial R$ is non-uniformly homeomorphic to the complement of a tame
(by construction) copy of $\mu_{d-1}$ in $S^{n+d}$.
\endremark

\demo{Proof} Let $\Phi\:S^{n+1}\x [0,\infty)\to P_{[0,\infty)}$ be
a Steenrod spheroid.
Let $K^{[1]},K^{[2]},\dots$ be the successive barycentric subdivisions of some
triangulation of $S^{n+1}$.
(Strictly speaking, to obtain the binary rather than the $(n+1)!$-ary tree in
the end, we should have decomposed each operation of barycentric subdivision into
a sequence of dichotomic subdivisions.)
Let $L_i$ be the $n$-skeleton of $K^{[i]}$ and let $N_i$ be the third derived
neighborhood of $L_i$ in $K^{[i]}$.
Note that $S^{n+1}\but\bigcup\Int N_i$ is a tame Cantor set $C$,
and there is a continuous map $\phi\:S^{n+1}\to S^{n+1}$ with
$\phi(\bigcup\Int N_i)=\bigcup L_i$ and $\phi(C)=S^{n+1}$.
(It is constructed analogously to the map $[0,1]\to [0,1]$ that converts
a ternary expansion $A=0.\alpha_1\alpha_2$ into a binary expansion
$B=0.\beta_1\beta_2\dots$ by replacing all twos with ones and halting
after the first occurrence of a one in $A$.)
Let $N=\bigcup N_i\x [i,\infty)\subset S^{n+1}\x [0,\infty)$.
On the other hand, let $N'$ be the third derived neighborhood of
$\bigcup L_i\x [i,\infty)$ in $S^{n+1}\x [0,\infty)$.
Let $E$ and $E'$ be the exteriors of $N$ and $N'$ respectively in
$S^{n+1}\x[0,\infty)$.

Let $f\:N'\x[0,1]\to P_{[l,\infty]}$ be the homotopy between
$\Phi|_{N'}$ and $r_k^{(n+1)}\Phi|_{N'}$.
Using the homotopy $\Pi_t$ from the beginning of \S2, we may assume that
$f^{-1}(X)=N'\x\{1\}$.
Let $h\:S^{n+1}\x[0,\infty)\to E'\cup\partial N'\x [0,1)$ be a homeomorphism
with $h(E)=E'$.
Let $\Phi'$ be defined by $\Phi h$ on $E$ and by $fh$ on $N$.
Then $\Phi'$ is properly homotopic to $\Phi$ (using $f$) and extends by
continuity to a uniformly continuous map
$\bar\Phi'\:S^{n+1}\x [0,\infty]\but C\x\{\infty\}\to P_{[0,\infty]}$. \qed
\enddemo

\proclaim{Theorem 6.5} Let $X$ be an LC$_{n-1}$ compactum, $n\ge 1$.

(a) If $\spi_1(X)=1$, then $\spi_n(X)@>\sstau>>\pi_n(X)$ is an epimorphism.

(b) If $\spi_1(X)$ acts trivially on $\spi_n(X)$, then
$\spi_n(X)@>\check\tau\sstau>>\check\pi_n(X)$ is an epimorphism.
\endproclaim

Note that in the case $n\ge 2$, the hypotheses in (a) and (b) can be weakened
using Theorem 3.15(b) to ``$\spi_1(X)$ is finite'' and ``a finite index
subgroup of $\spi_1(X)$ acts trivially on $\spi_n(X)$'', respectively.

The necessity of simple connectedness in Theorem 6.5 is shown by Example 5.8.

\demo{Proof}
Let $\phi\:(\partial R,pt)\to (X,x)$ be an $n$-spheroid with trunks realizing
the given element of $\pi_n(X)$ by Proposition 6.4.
We may assume that the infinite binary tree $T$ properly embedded in $\R^{n+1}$
is a full subcomplex of some triangulation of $\R^{n+1}$, which we fix from
now on, and that $R$ is the derived neighborhood of $T$.
Since $X$ is L$\sC_{n-1}$, $\phi$ extends to a map $\partial R\cup K\to X$,
also denoted $\phi$, where $K$ is the union of the dual cells to all but
finitely many of the edges of $T$.
Let $B_0,B_1,\dots$ be the closures of the components of $R\but K$.
Thus all but finitely many of $B_i$'s are the dual cells of vertices of $T$,
and the remaining ones are the derived neighborhoods of finite subtrees of $T$.
Without loss of generality, $B_0$ is the derived neighborhood of a finite
subtree $T_1\i T$, and $B_1,B_2,\dots$ are the dual cells of vertices of $T$.
In particular, all $B_i$'s are PL balls.

Let $U_i=\{U_{i\alpha}\}$ be a sequence of covers of $X$ such that each
$U_{i\alpha}$ has diameter at most $2^{-i}$ and each $U_{i+1,\alpha}$ is
contained in some $U_{i\beta}$.
Without loss of generality, the diameter of $X$ is at most $1$, in which
case we may assume that $U_0=\{X\}$.
Then $T$ is the union of an increasing chain of finite (simplicial) subtrees
$T_1\i T_2\i T_3\i\dots$ such that if $B_j$ is the dual cell of a vertex
of $T\but T_i$, then $\phi(\partial B_j)$ is contained in some $U_{i+1,\beta}$.
Let $R_i$ be the derived neighborhood of $T_i$, and let $B_{i1},\dots,B_{ir_i}$
be the dual cells of the vertices of $T_{i+1}\but T_i$ (where $T_0=\emptyset$).
By reindexing and duplicating the elements of each $U_i$ we may assume that
each $\phi(\partial B_{ij})\i U_{ij}$.
Let $\phi_0\:S^n\to x\in X$ be the constant map.
Assuming that $\phi_i\:(S^n,pt)\to (X,x)$ is defined, let us define
$\phi_{i+1}$ by attaching each $\phi|_{\partial B_{ij}}$ to $\phi_i$ along
some path in $U_{kl}$ from a point in $\phi(\partial B_{ij})$ to a point in
$\phi_i(S^n)$, where $U_{kl}$ contains $U_{ij}$ and if $k<m<i$ then any
$U_{mp}$ containing $U_{ij}$ is disjoint from $\phi_i(S^n)$.

Since $X$ is compact, $k=k(i)\to\infty$ as $i\to\infty$.
Since each $\phi_{i+1}$ is $2^{-k(i)}$-close to $\phi_i$, the sequence
$\phi_1,\phi_2,\dots$ uniformly converges to a continuous map
$\phi_\infty\:(S^n,pt)\to (X,x)$.
Let $P_i$ be the nerve of the cover $U_i$ (with duplicate elements omitted),
let us fix simplicial bonding maps $P_{i+1}\to P_i$, arising due to
the hypothesis that $U_{i+1}$ refines $U_i$, and let $p^\infty_i\:X\to P_i$
be the projection.
Each $p^\infty_i\phi_\infty\:(S^n,pt)\to (P_i,p_i)$ is homotopic to
$p^\infty_i\phi_i$ by the ``rectilinear'' homotopy.
Consecutive pairs of these homotopies yield homotopies $h^{(1)}_i$
between $p^\infty_i\phi_{i+1}$ and $p^\infty_i\phi_i$.
A representative $\Phi\:\Int R\to P_{[0,\infty)}$ of $[\phi]$ extending to
a uniformly continuous map $\bar\Phi\:R\to P_{[0,\infty]}$ agreeing with
$\phi$ on $\partial R$ can be obtained by combining some homeomorphisms
$f_i\:(S^n,pt)\to(\partial R_i,pt)$ with the ``rectilinear'' homotopies
$h^{(2)}_i\:(S^n,pt)\x I\to (P_{[i,i+1]},p_{[i,i+1]})$ between the maps
$p^\infty_i\phi f_{i+1}$ and $p^\infty_i\phi f_i$.

By construction, each $\phi_i=(\psi_{11}+\dots+\psi_{1r_1})+\dots+
(\psi_{i1}+\dots+\psi_{ir_i})$, where each $\psi_{ij}\:(S^n,pt)\to (X,x)$ is
freely (i.e.\ unpointed) homotopic to $\phi|_{\partial B_{ij}}$.
Similarly, each $\phi f_i=(\chi_{11}+\dots+\chi_{1r_1})+\dots+
(\chi_{i1}+\dots+\chi_{ir_i})$, where each $\chi_{ij}\:(S^n,pt)\to (X,x)$ is
freely homotopic to $\phi|_{\partial B_{ij}}$.
Since $\spi_1(X)$ acts trivially on $\spi_{n+1}(X)$,
the forgetful map $\spi_n(X)\to[S^n,\,X]^\tr$ into the set of
free homotopy classes is a bijection.
So each $\psi_{ij}$ is (pointed) homotopic to $\chi_{ij}$ by a homotopy $h_{ij}$,
and consequently each $\phi_i$ is homotopic to $\phi f_i$ by a homotopy $h_i$.
Moreover, if $\spi_1(X)=1$, the homotopies $h_{ij}$ can be chosen to go
through maps $(S^n,pt)\to (Y,0)\to (X,x)$, where $Y$ is the one-point
compactification of $\R^n\sqcup [0,\infty)$, such that the image of $\R^n$
is contained in $U_{ij}$.
Hence each map $S^{n+1}\to P_i$ given by ``rectilinear'' null-homotopies
of $p^\infty_i\psi_{ij}$ and $p^\infty_i\chi_{ij}$ along with the homotopy
$p^\infty_i h_{ij}$ is null-homotopic (by a ``rectilinear'' homotopy).
It follows that each $p^\infty_ih_{i+1}$ is homotopic to $p^\infty_ih_i$ by a
``rectilinear'' homotopy extending $h^{(1)}_i$ and $h^{(2)}_i$. \qed
\enddemo

Theorem 6.5 and Theorem 6.1(g) have

\proclaim{Corollary 6.6} Let $X$ be a compactum with $\spi_1(X)=1$,
and let $n\ge 1$.
Then $X$ is semi-LC$_n$ iff $X$ is LC$_{n-1}$ and $\pi_n(X)$
is discrete.
\endproclaim

A $0$-dimensional counterpart of Corollary 6.6 (without the assumption of
simple connectedness) is obtained in \S8 (see Corollary 8.8).
There will follow a sketch of an alternative approach to proving an assertion,
close to Corollary 6.6; the author hopes that this approach will help to
shed light on the necessity of simple connectedness in the hypothesis of
Corollary 6.6.
\bigskip

Let us now turn to homological versions of the above results.

\proclaim{Theorem 6.7 (Eda--Kawmura, Shchepin)}
Let $X$ be an LC$_{n-1}$ compactum, $n\ge 1$.
Then $\sH_n(X)@>\sstau>>H_n(X)$ is an epimorphism.
\endproclaim

The case $n=1$ was established by Shchepin \cite{MRS; Theorem 4.1} using
Bestvina's thesis and other non-trivial structure results in topology of
compacta.
The author learned from Shchepin about five years ago that he could also prove
the general case of Theorem 6.7; this proof has never been written up.
It should be noted that some key parts of the argument in \cite{MRS} do not
seem to generalize for arbitrary $n$, and the author does not know how these
difficulties have been tackled by Shchepin.
Eda and Kawamura proved the following corollary of Theorem 6.7: if $X$ is
L$\sC_{n-1}$ and $n\ge 1$, $\sH_n(X)@>\check\tau\sstau>>\check H_n(X)$ is
an epimorphism \cite{EK1}.

The restriction $n\ge 1$ is necessary in Theorem 6.7 by Example 5.4.

\demo{Proof} Similarly to Proposition 6.4 one shows that if $X$ is L$\sC_{n-1}$
(or just HL$\sC_{n-1}$), every element of $H_n(X)$, $n\ge 1$, can be represented
by a uniformly continuous map $M\to X$ from an oriented $n$-pseudo-manifold
that admits a proper $\eps$-map $M\to T$ (i.e.\ a map whose point-inverses have
diameters $<\eps$) onto the infinite binary tree, with respect to some proper
function $\eps\:T\to (0,1]$.
Since $X$ is L$\sC_{n-1}$ (or just HL$\sC_{n-1}$), ``each trunk'' in this ``cycle
with trunks'' can be partitioned into an infinite sum of small singular cycles.
Then by Lemma 6.3(a), the L$\sC_{n-1}$ condition allows to replace $M$ with
the sum of a singular compact pseudo-manifold $M_0\to X$ and a finite number of
spheroids with trunks $\partial R_i\to X$.
The latter can be converted within their Steenrod homology classes (not within
their Steenrod homotopy classes!) into genuine spheroids $S^n\to X$ by
the proof of Theorem 6.5 (which is simplified as one does not have to keep
track of basepoints now). \qed
\enddemo

\proclaim{Theorem 6.8} Fix some integer $n\ge 0$ and let $X$ be a compactum.

\smallskip
(a) {\rm (Shchepin)} HL$\sitC_n$ \imp\ HLC$_n$ \iff\ HL\v C$_n$.

\smallskip
(b) {\rm (Jussila \cite{Ju}; see also \cite{Br; VI.10.6})} If $X$ is
semi-HLC$_n$, then $H_n(X)@>\check\tau>>\check H_n(X)$ is an isomorphism.

\smallskip
(c) If $X$ is HL$\sitC_n$, then $\sH_n(X)@>\sstau>>H_n(X)$ is an isomorphism.

\smallskip
(d) {\rm (Marde\v si\'c \cite{Ma}; see also \cite{Br; VI.12.6})} If $X$ is
semi-HL$\sitC_n$, then the composition $\sH_n(X)@>\check\tau\sstau>>\check H_n(X)$
is an isomorphism.

\smallskip
(e) If $X$ is HL$\sitC_{n-1}$, then $\sH_n(X)@>\check\tau\sstau>>\check H_n(X)$ is
onto a dense subset.

\bigskip
(f) $X$ is semi-HL$\sitC_n$ iff $X$ is HL$\sitC_{n-1}$ and $\sH_n(X)$ is discrete.

\smallskip
(g)  $X$ is semi-HLC$_n$ iff $X$ is HLC$_{n-1}$ and $H_n(X)$ is discrete.

\smallskip
(h) $X$ is semi-HL\v C$_n$ iff $X$ is HL\v C$_{n-1}$ and
$\check H_n(X)$ is discrete.
\endproclaim

Part (g) follows immediately from (a), (b) and (h).
Part (c) follows from (a), (b) and (d) (compare \cite{Br; V.12.15}), but
we also give a direct proof below.

Part (a) was proved by Shchepin \cite{MRS; Theorem 4.10} in the case $n=1$.
The author learned from Shchepin about five years ago that he could also
prove the general case of (a); this proof has never been written up.
While our proof of (a) is based on the same idea as the argument of \cite{MRS},
namely, the Hurewicz Construction (see the proof of Lemma 6.1(3)), our
realization of this idea in higher dimensions depends on a more delicate
technique of the geometric approach to cohomology, developed by Buoncristiano,
Rourke and Sanderson \cite{BRS}.

The implication HL$\sC_\infty$ \imp HLC$_\infty$ cannot be reversed
for the cluster of copies of an arbitrary non-simply-connected acyclic
$2$-polyhedron \cite{EKR}.

\definition{Pseudo-comanifolds} We call a proper PL map $f\:W\to P$ between
polyhedra a {\it (co-connected) $k$-pseudo-comanifold} if there exists
a triangulation of $P$ such that the preimage $f^{-1}(\Delta^i)$ of every
its $i$-simplex is a (pseudo-connected) $(i-k)$-pseudo-manifold with boundary
equal to the preimage of the boundary $f^{-1}(\partial\Delta^i)$.
Such a triangulation of $P$ will be called {\it transverse} to $f$.
Co-manifolds originate in \cite{BRS}, where they are called ``mock bundles''.

An embedded $k$-pseudo-comanifold $f\:W\emb P$ is {\it co-orientable} if
$H_k(P,\,P\but W)$ contains no torsion.
In this case $H^k(P,\,P\but W)$ is free abelian, and a choice of a set of its
generators, representable by cocycles with disjoint supports, is called
a {\it co-orientation} of $f$.
(For reference, a $k$-pseudo-manifold $M$ without boundary is orientable iff
the compactly supported cohomology $H^k_ñ(M)$ is torsion-free.
In this case the locally finite homology $H_k^\lf(M)$ is free abelian, and
a choice of a set of its generators, representable by cycles with disjoint
supports, corresponds to an orientation of $M$.)
If $f\:W\to P$ is an arbitrary pseudo-comanifold, where $P$ is compact,
$f$ is the projection of an embedded pseudo-comanifold $\bar f\:W\to P\x\R^N$
for some $N$, and the {\it co-orientability (co-orientation)} of $f$ can
be well-defined as that of $\bar f$.

It is easy to see that $H^k(P;\,\Z/2)$ (respectively $H^k(P)$) is isomorphic to
the group of (oriented) pseudo-cobordism classes of (oriented)
$k$-pseudo-comanifolds in $P$, cf\. \cite{BRS}.

Let $f\:W\to P$ be an oriented $k$-pseudo-comanifold.
Note that if $P$ is an oriented $i$-pseudo-manifold with respect to
a triangulation of $P$, transverse to $f$, then $W$ will be an oriented
$(i-k)$-pseudo-manifold, cf\. \cite{BRS; II.1.2}.
On the other hand, if $\phi\:M\to P$ is simplicial with respect to
a triangulation of $P$, transverse to $f$, the pullback $\phi^*(f)\:N\to M$ is
an oriented $k$-pseudo-comanifold, cf\. \cite{BRS; bottom of p.\ 23}.
Now if $\phi$ happens to be an oriented singular $i$-pseudo-manifold, we obtain
an oriented singular $(i-k)$-pseudo-manifold $f^!(\phi)\:N\to W$.
Using the techniques of \cite{BRS}, it is easy to show that
$f^!\:H_i(P)\to H_{i-k}(W)$ is well-defined.
Furthermore, $f_*f^![\phi]=[f]\smallfrown[\phi]$, cf\. \cite{BRS; p.\ 29}.

$$\CD
N@>f^!(\phi)>>W\\
@V\phi^*(f)VV@VfVV\\
M@>\phi>>P
\endCD$$

If $f\:W\to P$ is a co-connected oriented $0$-pseudo-comanifold, it is easy to
see that $[f]=\pm[\id_P]\in H^0(P)$.
Let us call such an $f$ an {\it elaboration} of the polyhedron $P$ if
specifically $[f]=[\id_P]\in H^0(P)$.
Then, in particular, $f_*f^!=\id\:H_i(P)\to H_i(P)$ for each $i$.
Hence $H_i(P)$ can be identified with a direct summand of $H_i(W)$.
In addition, using once again the co-connectedness of $f$, for each oriented
singular $i$-pseudo-manifold $\phi\:M\to P$, simplicial with respect to
a triangulation of $P$, transverse to $f$, the elaboration $\phi^*(f)\:N\to M$
induces an isomorphism $H_i(N)\to H_i(M)$ due to the pseudo-connectedness of $N$.
\enddefinition

\demo{Elaborated Ferry Construction}
Assuming that $X$ is HL$\sC_n$, we replace the retraction
$r_k^{(n+1)}\:P_{[k,\infty)}^{(n+1)}\cup X\to X$ and its homotopy to
the identity with values in $P_{[l,\infty]}$ keeping $X$ fixed (see
``Ferry's Construction'' in the beginning of this section) by
the following data:

\medskip
(i) an elaboration $q^{(n+1)}\:Q^{n+1}\to P_{[k,\infty)}^{(n+1)}$;

\smallskip
(ii) a retraction $\bar r_k^{(n+1)}\:Q^{n+1}\cup X\to X$;

\smallskip
(iii) a pseudo-cobordism $W^{n+2}\to P_{[k,\infty)}^{(n+1)}$ between
$q^{(n+1)}$ and $\id_{P_{[k,\infty)}^{(n+1)}}$;

\smallskip
(iv) a map $W^{n+2}\cup X\to P_{[l,\infty]}$ restricting to $\bar r_k^{(n+1)}$
on $Q^{n+1}\cup X$ and to the inclusion on $P_{[k,\infty)}^{(n+1)}$.
\medskip

These data are constructed similarly to the Ferry Construction, item (iii)
being granted by the definition of an elaboration.

If $X$ is semi-HL$\sC_n$, then $\bar r_k^{(n)}\:Q^n\cup X\to X$ restricted
to $Q^n$ extends to a continuous map $Q^{n+1}\to X$.
\enddemo

\demo{Proof of Theorem 6.8. (c), (d), (e)}
Using the elaborated Ferry Construction, these parts can be now be proved
similarly to the corresponding parts of Theorem 6.1. \qed
\enddemo

\definition{Fractalized pseudo-manifolds}
An inverse sequence $\dots@>f_1>>P_1@>f_0>>P_0$ of polyhedra and PL maps will be
called {\it fractalizing} if each $f_i$ is an elaboration, transverse to
a triangulation $T_i$ of $P_i$ such that $f_i$ is simplicial with respect to
$T_{i+1}$ and some subdivision of the barycentric subdivision of $T_i$.
Such a sequence of triangulations $(T_0,T_1,\dots)$ will also be called {\it
fractalizing}.
The limit $F$ of such an inverse sequence is a {\it fractalized} $P_0$
and $p^\infty_0\:F\to P_0$ is a {\it fractalization} of $P_0$.

If $Q$ is an $n$-pseudo-manifold and $f\:M\to Q$ its fractalization, then by
the above $f_*\:H_n(M)\to H_n(Q)$ is an isomorphism, and in particular $M$ has
a well-defined fundamental class $[M]$ in Steenrod homology.
If $X$ is a compactum and $\phi\:M\to X$ is a map, we say that it {\it represents}
$\phi_*([M])\in H_n(X)$.
\enddefinition

\proclaim{Lemma 6.9} Let $X$ be a HL\v C$_n$ compactum, and let $\phi\:M\to X$
be a singular fractalized oriented $n$-pseudo-manifold.
If $\check\tau([\phi])=0\in\check H_n(X)$, then $\phi$ bounds a singular
fractalized oriented null-pseudo-bordism.
\endproclaim

The idea of the proof is based on the Hurewicz Construction (see the proof of
Lemma 6.1(3)).
The case $n=1$ was obtained by Shchepin \cite{MRS; Theorem 4.4}.
This case is special in that $2$-dimensional fractal pseudo-manifolds
reduce to $2$-dimensional fractal manifolds, which are all homeomorphic to
each other, see \cite{MRS}.

\demo{Proof} We may assume that the assertion holds in dimensions $<n$.
Suppose that $M$ is the limit of a fractalizing inverse sequence
$\dots@>f_1>>Q_1@>f_0>>Q_0$ with $Q_0$ an oriented $n$-pseudo-manifold.
Suppose that $X$ is the inverse limit of polyhedra $P_i$ and PL maps $p_i$.
Without loss of generality, $\phi$ extends to a level-preserving map
$\phi_{[0,\infty]}\:Q_{[0,\infty]}\to P_{[0,\infty]}$ (see \S2).
By an abuse of notation, let us redefine $Q_{[0,\infty]}$ by replacing each
mapping cylinder $MC(f_i)$ with a pseudo-cobordism between the elaboration
$f_i\:Q_{i+1}\to Q_i$ and the identity of $Q_i$.
Then $\phi_{[0,\infty)}$ is a Steenrod cycle representing $[\phi]$, thus by
the hypothesis each $\phi_i$ bounds a singular oriented pseudo-manifold
$\hat\phi_i\:\hat Q_i\to P_i$.
Given an $\eps>0$, there exists a $k$ such that there is a triangulation
of the pair $(\hat Q_k,Q_k)$, extending a triangulation of $Q_k$ in
a fractalizing sequence of triangulations, and such that the image of
every its simplex under $\hat\phi_k$ is $\eps$-close to $X$ and has
diameter at most $\eps$.
Note that for each simplex $\Delta$ of this triangulation of $Q_k$, the inverse
sequence $\dots\to (f^{k+2}_k)^{-1}(\Delta)\to (f^{k+1}_k)^{-1}(\Delta)\to\Delta$
is fractalizing.
Then similarly to Ferry's Construction (see the proof of Theorem 6.1),
the inductive hypothesis implies that for each $\delta>0$, the number $\eps$
can be chosen so that $\phi_{[k+1,\infty]}$ extends to a level-preserving map
$\psi_{[k+1,\infty]}\:R_{[k+1,\infty]}\to P_{[k+1,\infty]}$, where the inverse
sequence $\dots@>g_{k+1}>>R_{k+1}@>g_k>>R_k:=Q_k\cup\hat Q_k^{(n)}$
is fractalizing, each $g_i^{-1}(Q_i)\i R_{i+1}$ is identified with $Q_{i+1}$ so
that $g_i|_{Q_{i+1}}$ is identified with $f_i|_{Q_{i+1}}$, and
$\psi_{[k+1,\infty]}$ is $\delta$-close to the composition of the projection
onto $R_k$ and the restriction of $\hat\phi_k$.
(Here $R_{[k+1,\infty]}$ consists of pseudo-cobordisms rather than mapping
cylinders.)

Let $\Delta$ be an $(n+1)$-simplex of $\hat Q_k$, and let
$S_i=(g^i_k)^{-1}(\partial\Delta)\i R_i$ for each $i>k$.
Let $S_\infty$ be the limit of the fractalizing inverse sequence
$\dots\to S_{k+1}\to S_k$.
Then the restriction $\chi_{[k+1,\infty)}$ of $\psi_{[k+1,\infty)}$ to
$S_{[k+1,\infty)}$ is a Steenrod cycle representing the image of
the fundamental class of $S_\infty$ in $H_n(X)$.
Recall that the image of $S_\infty$ in $X$ has diameter at most $\delta$.
Since $X$ is HL\v C$_n$, for each $\gamma>0$ this $\delta$ can
therefore be chosen so that $\chi_{[k+1,\infty)}$ be \v Cech trivial with
support in the $\gamma$-neighborhood in $X$ of the image of $S_\infty$.

Then each $\chi_i$ bounds a singular connected oriented pseudo-manifold
$\hat\chi_i\:\hat S_i\to P_i$ with image in the $\gamma$-neighborhood of
$\chi_i(S_i)$.
Given any specific value of $i$, we may extend the restriction
$S_i\to\partial\Delta$ of $g^i_k$ to a map $\hat S_i\to\Delta$, simplicial with
respect to some subdivision of the barycentric subdivision of $\Delta$.
Repeating the same for each $(n+1)$-simplex $\Delta$ of $\hat Q_k$, we obtain
an extension of $g^i_k\:R_i\to Q_k\cup\hat Q_k^{(n)}$ to an elaboration
$\hat g^i_k\:\hat R_i\to\hat Q_k$ such that the various $\hat\chi_i$ combine
into a map $\hat\psi_i\:\hat R_i\to P_j$.
Since $\hat\psi_i$ and $\hat\phi_k\hat g^i_k$ are $\gamma$-close, given an $l$,
we may choose $\gamma$ so that they extend to a map
$\hat\psi_{[k,i]}\:\hat R_{[k,i]}\to P_{[l,i]}$ of a pseudo-cobordism
$\hat g^{[i,k]}_k\:\hat R_{[k,i]}\to\hat Q_k$ between $\hat g^i_k$ and
the identity of $\hat Q_k$, which restricts to
$\psi_{[k,i]}\:R_{[k,i]}\to P_{[k,i]}$.

Now given any $l'$, we have a similar string of dependencies
$l'\mapsto\gamma'\mapsto\delta'\mapsto\eps'\mapsto k'$, and the above
procedure applies to $\hat\chi_i$ with $i=k'$ in place of $\hat\psi_k$.
Proceeding in this fashion, we eventually obtain a fractalizing inverse
sequence $(\dots\to T_1\to T_0)=(\dots\to\hat R_i\to\hat Q_k)$ and a proper map
$T_{[0,\infty)}\to P_{[l,\infty)}$, which by construction extends to
a continuous map of the inverse limit $T_\infty\to X$ that yields the required
fractalized oriented null-pseudo-bordism of $\phi$. \qed
\enddemo

\demo{Fractalized Ferry Construction}
Assuming that $X$ is a HL\v C$_n$ compactum, one can apply Lemma 6.9 to
construct:

\medskip
(i) a fractalizing inverse sequence
$\dots@>q_1>>Q_1@>q_0>>Q_0=P_{[k,\infty)}^{(n+1)}$
with inverse limit $Q_\infty$;

\smallskip
(ii) a retraction $\tl r_k^{(n+1)}\:Q_\infty\cup X\to X$;

\smallskip
(iii) a `telescopic pseudo-cobordism'
$q^{[0,\infty]}_0\:W_{[0,\infty]}\to P_{[k,\infty)}^{(n+1)}$,
combining pseudo-cobordisms $q^{[i,i+1]}_i\:W_{[i,i+1]}\to Q_i$
between $q_i$ and $\id_{Q_i}$ (thus $W_{\N\cup\infty}=Q_{\N\cup\infty}$);

\smallskip
(iv) a `telescopic homotopy'
$h^{(n+2)}_k\:W_{[0,\infty]}\cup X\to P_{[l,\infty]}$ restricting to
$\tl r_k^{(n+1)}$ on $Q_\infty\cup X$ and to the inclusion
on $P_{[k,\infty)}^{(n+1)}$.
\enddemo

\proclaim{Corollary 6.10} {\rm (Shchepin)} If $X$ is a HL\v C$_n$ compactum,
every element of $H_n(X)$ is representable by a singular fractalized oriented
$n$-pseudo-manifold.
\endproclaim

This follows from the fractalized Ferry Construction similarly to the proof of
Lemma 6.2(2).
The author learned from E. V. Shchepin about five years ago that he could
prove a version of Corollary 6.10 using another (possibly different) notion of
a fractalized pseudo-manifold, defined by a direct inductive construction.

\demo{Proof of Theorem 6.8, continued. (a)} The implication HLC$_n$ \imp\
HL\v C$_n$ is trivial.
The converse implication follows from local versions of Corollary 6.10 and
Lemma 6.9.
The implication HL$\sC_n$ \imp\ HL\v C$_n$ follows from (d), or
alternatively the implication HL$\sC_n$ \imp\ HLC$_n$ follows from (c).
\qed
\enddemo

\demo{(b)} This is similar to the proof of Corollary 6.6, but easier.

Let us represent an element of $H_n(X)$ by a proper singular oriented
$(n+1)$-pseudo-manifold $\phi_{[0,\infty)}\:Q_{[0,\infty)}\to P_{[0,\infty)}$.
(Here each $Q_{[i,i+1]}$ is an arbitrary oriented pseudo-bordism, rather than
a mapping cylinder.)
If the image of $[\phi_{[0,\infty)}]$ under the composition
$H_n(X)@>\check\tau>>\check H_n(X)@>p^\infty_k>>H_n(P_k)$
is trivial, we may assume that $\phi_k\:Q_k\to P_k$ bounds a singular oriented
$(n+1)$-pseudo-manifold $\hat\phi_k\:\hat Q_k\to P_k$.
Without loss of generality, the Steenrod cycle
$\phi_{[k,\infty)}\cup\hat\phi_k\:Q_{[k,\infty)}\cup\hat Q_k\to P_{[0,\infty)}$
maps simplicially into $P_{[k,\infty)}^{(n+1)}$.
Let us write $Q=Q_{[k,\infty)}\cup\hat Q_k$ and
$\pi=\phi_{[k,\infty)}\cup\hat\phi_k\:Q\to P_{[0,\infty)}$.
Given an $\eps>0$, we may choose $k$ so that the fractalized Ferry Construction
applies, with $\tl r_k^{(n)}$ sending $(q^\infty_0)^{-1}(\partial\Delta)$, for
each $(n+1)$-simplex $\Delta$ of $P_{[0,\infty)}$, into a subset of $X$ of
diameter $<\eps$.
Given an $(n+1)$-simplex $\Delta$ of $Q$, let $\psi_\Delta$ denote the Steenrod
cycle $\phi|_\Delta\cup h^{(n+1)}_k|_{R_{[0,\infty)}}\:
\Delta\cup R_{[0,\infty)}\to P_{[0,\infty)}$, where
$R_i=(q^i_0)^{-1}(\phi(\partial\Delta))$ form a fractalizing inverse sequence
$\dots\to R_1\to R_0$.
Then $[\phi]=[\bigsqcup\psi_\Delta]$, where $\Delta$ runs over all
$(n+1)$-simplices of $Q$, already on the level of locally finite simplicial
chains of $P_{[0,\infty)}$.

Since $X$ is compact, it can be covered by a finite collection of sets
$U_1,\dots,U_r$ of diameters at most $2\eps$ such that each
$\tl r_k^{(n)}((q^\infty_0)^{-1}(\partial\Delta))$ is
contained in $U_i$ for some $i=i(\Delta)$.
For each $i$, consider the Steenrod cycle
$\psi_i:=\bigsqcup_{i=i(\Delta)}\psi_\Delta$.
We have $[\phi]=[\psi_1]+\dots+[\psi_r]$ already on the level of locally finite
chains.
On the other hand, each $[\psi_i]$ lies in the image of $H_n(U_i)$, and
therefore is trivial by the hypothesis.
This proves that the composition
$H_n(X)@>\check\tau>>\check H_n(X)@>p^\infty_k>>H_n(P_k)$ is injective
for some $k$, which implies the assertion. \qed
\enddemo

\demo{(f),(h)} These are proved similarly to the corresponding parts of Theorem
6.1. \qed
\enddemo

We record one byproduct of the preceding discussion.

\proclaim{Theorem 6.10$'$} If $n>0$, Corollary 6.10 holds under the weaker
hypothesis of HLC$_{n-1}$.
\endproclaim

\demo{Proof} This is similar to the proof of Theorem 6.5, but easier, since one
does not need to worry about basepoints, and Proposition 6.4 can be replaced by
the construction in the proof of Theorem 6.8(b). \qed
\enddemo

\proclaim{Theorem 6.11} Let $X$ be a compactum, $n\ge 0$.

\medskip
(a) The following are equivalent:

\medskip
\noindent
(i) $X$ is LC$_n$;

\smallskip
\noindent
(ii) for any closed $Y,Z\i X$ with $Y\i\Int Z$, the inclusion $Y\emb Z$ factors
through an LC$_n$ compactum $W$;

\smallskip
\noindent
(iii) for any closed $Y,Z\i X$ with $Y\i\Int Z$, $\im[\pi_i(Y)\to\pi_i(Z)]$ is
countable for all $i\le n$.

\bigskip
(b) The following are equivalent:

\medskip
\noindent
(i) $X$ is HLC$_n$;

\smallskip
\noindent
(ii) for any for any closed $Y,Z\i X$ with $Y\i\Int Z$, $Y$ is contained in an
LC$_n$ compactum $W$ such that the inclusion $Y\emb Z$ factors through
an inclusion of $Y$ into a compactum $\hat W$ such that $\hat W\but Y$ is
a fractalization of $W\but Y$;

\smallskip
\noindent
(iii) for any closed $Y,Z\i X$ with $Y\i\Int Z$, $\im[H_i(Y)\to H_i(Z)]$ is
finitely generated for all $i\le n$.
\endproclaim

The implication (i)\imp(iii) in (b) is due to Borel and Moore \cite{BoM}, and its
converse to Dydak \cite{D6}.
The earliest result of this kind is the cohomological analogue of (i)\iff(iii),
which is due to Wilder (see \cite{Br}).

\demo{Proof. (a)} Assuming (i), let $Y$ and $Z$ be as in (ii), and suppose that
$(X,Y)=\invlim(P_i,Q_i)$.
By Ferry's Construction (see the proof of Theorem 6.1),
$P_{[k,\infty)}^{(n+1)}\cup X$ retracts onto $X$ for some $k$.
Pick an $l>k$ such that this retraction sends $W:=Q_{[l,\infty)}^{(n+1)}\cup Y$
into $Z$.
Clearly, $W$ is L$\sC_n$.

Assuming (ii), by parts (h) and (c) of Theorem 6.1, $\pi_i(W)$ is discrete for
all $i\le n$.
Hence by Lemma 3.4(b), each $\pi_i(W)$ is countable, which implies (iii).

Assuming (iii), let $\dots\i N_2\i N_1$ be a fundamental sequence of closed
neighborhoods of an arbitrary given point $x\in X$.
Then by Theorem 3.1(c), $\invlim\pi_k(N_i)=0=\derlim\pi_k(N_i)$ for all $k$.
If $G_i=\im[\pi_k(N_{i+1})\to\pi_k(N_i)]$, then we still have
$\invlim G_i=0=\derlim G_i$ for all $k$.
Since $G_i$ are countable for all $k\le n$, by Lemmas 3.3 and 3.4(a), for
each $i$ there exists a $j>i$ such that $G_j\to G_i$ is trivial for all
$k\le n$.
This implies (i). \qed
\enddemo

\demo{(b)} The proof of (iii)\imp(i) is similar to that in (a).
The proof of (i)\imp (ii) is similar to that in (a), using the fractalized Ferry
Construction.

Assuming (ii), similarly to the proof of (ii)\imp(iii) in (a), we get that
$H_i(W)$ injects into a finitely generated abelian group.
Hence it is finitely generated, and the assertion follows from the
commutativity of the diagram
$$\CD
H_i(\hat W)@>f_*>>H_i(W)\\
@Aincl_*AA@Vf^!VV\\
H_i(Y)@>incl_*>>H_i(\hat W)@>(\tl r_k^{(n+1)})_*>>H_i(Z). \qed
\endCD$$
\enddemo

An inverse sequence of groups $G_i$ is said to be {\it nearly
Mittag-Leffler} if for each $i$ there exists a $j>i$ such that for each $k>j$
the image of $G_j\to G_i$ is contained in the normal closure of the image
of $G_k\to G_i$.
If $X$ is the limit of an inverse sequence of compact connected polyhedra
$P_i$, the property that the inverse sequence $\pi_1(P_i)$ is nearly
Mittag-Leffler clearly does not depend on the choice of the basepoint of $X$
and is an invariant of the equivalence relation in Proposition 2.6(iii).
Then we may call $X$ {\it nearly Steenrod connected} if this property is
satisfied.
Thus the property of being nearly Steenrod connected is an (unpointed) shape
invariant.
It was introduced by McMillan \cite{Mc} (under a different name); it is proved
in \cite{Mc} that a continuous image of a nearly Steenrod connected compactum
is nearly Steenrod connected.
Clearly (see Theorems 3.1(b), 4.1(iii) and Lemma 3.3),

\bigskip
\centerline{Steenrod connected (i.e.\ $\pi_0(X)=0$)}
\centerline{$\Downarrow$}
\centerline{nearly Steenrod connected}
\centerline{$\Downarrow$}
\centerline{homologically Steenrod connected (i.e.\ $H_0(X)=0$).}
\bigskip

\proclaim{Theorem 6.12} Let $X$ be the limit of an inverse sequence of
compact connected polyhedra $P_i$.
Let us consider the one-point compactification $P_{[0,\infty]}/X$ of
$P_{[0,\infty)}$.

(a) {\rm (Shrikhande \cite{Sh})} $P_{[0,\infty]}/X$ is LC$_1$ iff $X$ is nearly
Steenrod connected;

(b) {\rm (Dydak \cite{D6})} $P_{[0,\infty]}/X$ is HLC$_n$
iff $H_i(X)$ is discrete for all $i<n$.
\endproclaim

We note that by Theorem 6.1(b), $P_{[0,\infty]}/X$ is LC$_n$ if and only if
it is both LC$_1$ and HLC$_n$.

\demo{Proof of (b)} The proof of Theorem 3.12 works to show that $H_k(X)$ is
discrete for all $k<n$ iff for each $i$ there exists a $j>i$ such that
$H_k(P_{[j,\infty]},X)$ maps trivially to $H_k(P_{[i,\infty]},X)$ for each $k<n$.
By the Map Excision Axiom 4.1(i), the latter is equivalent to
$P_{[0,\infty]}/X$ being HLC$_n$ at the compactifying point $X/X$. \qed
\enddemo

The proof of (a) can be carried out along the lines of the proof of
Theorem 3.12; we leave the details to the interested reader.

\head 7. Covering theory \endhead

\proclaim{Theorem 7.1} Let $X$ be a connected locally connected
compactum and let $d\in\{1,2,\dots,\infty\}$.
The monodromy map yields a bijection between $d$-fold covering maps over $X$
(up to a fiberwise homeomorphism) and representations of the topological group

(a) $\spi_1(X)$;

(b) $\check\pi_1(X)$ (or $\pi_1(X)$)

\noindent
into the discretely topologized symmetric group $S_d$ (up to an inner
automorphism of $S_d$).
Connected covering spaces correspond to transitive representations.
\endproclaim

If $X$ is a semi-L$\sC_1$ compactum, $\spi_1(X)$ and $\check\pi_1(X)$
coincide and are discrete by parts (d), (f), (h) of Theorem 6.1, so the usual
covering theory is a special case of each of the assertions (a) and (b) of
Theorem 7.1.

Fox mentions Theorem 7.1(a) \cite{F1; p.\ 2} with the following comment:
``I did this in lectures at the University of Mexico in the summer of 1951.
It has since been independently discovered by others, and appears for example
in [Spanier's textbook], p.\ 82''.
Concerning those others, see in particular \cite{AM}.

Using Theorem 6.1(e), one can deduce Theorem 7.1(b) directly from 7.1(a):

\demo{Proof of (b)} Corollary 2.5.3, Lemma 2.5.11 and Theorem 2.5.13 of
\cite{Sp} imply that fiberwise homeomorphism classes of connected $d$-fold
covering maps $p\:\tl X\to X$ correspond bijectively to conjugacy classes of
index $d$ subgroups of the singular fundamental group $\spi_1(X)$ that contain
the preimage of some neighborhood of $1$ in $\check\pi_1(X)$ under
the composition $\spi_1(X)@>\sstau>>\pi_1(X)@>\check\tau>>\check\pi_1(X)$.
Since $\check\pi_1(X)$ is zero-dimensional, a base of neighborhoods of $1$ is
given by all clopen sets.
Therefore the subgroups in question are precisely those that contain the kernel
of $\check\tau\stau$ composed with some continuous representation of
$\check\pi_1(X)$ into a discrete group.
On the other hand, by Theorem 6.1(e) the image of $\check\tau\stau$
is dense in $\check\pi_1(X)$.
Thus the fiberwise homeomorphism classes of connected $d$-fold covering spaces
of $X$ correspond bijectively to conjugacy classes of index $d$ subgroups of
$\check\pi_1(X)$ that contain the kernel of a continuous representation into
a discrete group.
If $H$ is such an index $d$ subgroup, the intersection $K$ of all its conjugates
contains the kernel into the discrete group, and so $\check\pi_1(X)/K$ is
discrete.
The action of $\check\pi_1(X)/K$ on the right cosets of its index $d$ subgroup
$H/K$ yields a transitive representation $\check\pi_1(X)/K\to S_d$.
This proves the transitive case of (b), which implies the general case
using that a locally connected compactum has open connected components. \qed
\enddemo

In order to generalize Theorem 7.1 to non-locally-connected spaces, Fox
``corrected'' the notion of a covering \cite{F1}, \cite{F2}.
A minor inaccuracy in Fox's theory was in turn corrected in \cite{Mo}.
See also \cite{MaM} and references there.
We recall that a {\it covering map} is a map $p\:\tl X\to X$ such that there
exists a cover $\{U_\alpha\}$ of $X$ satisfying

\medskip
(i) each $p^{-1}(U_\alpha)=\bigsqcup_\lambda U_\alpha^\lambda$, where each
$p|_{U_\alpha^\lambda}$ is a homeomorphism onto $U_\alpha$.
\medskip

\noindent
In this case we shall say that the cover $\{U_\alpha^\lambda\}$ of $\tl X$
{\it lies over} $\{U_\alpha\}$.
An {\it overlay structure} on $p$ is a cover $\{U_\alpha^\lambda\}$ of
$\tl X$ lying over some cover of $X$ and additionally satisfying

\medskip
(ii) if $U_\alpha^\lambda\cap U_\beta^\mu$ and
$U_\alpha^\lambda\cap U_\beta^\nu$ are both nonempty, then $\mu=\nu$.
\medskip

Two overlay structures $\{U_\alpha^\lambda\}$ and $\{V_\alpha^\lambda\}$ on $p$
are {\it equivalent} if there exists an overlay structure
$\{W_\alpha^\lambda\}$ on $p$ refining each of them (as covers of $\tl X$).
An {\it overlaying} is an equivalence class of overlay structures on
a covering map.
An {\it isomorphism} of overlayings $[p\:\tl X\to X;\, \{U_\alpha^\lambda\}]$
and $[p'\:\tl X'\to X;\, \{V_\beta^\lambda\}]$ is a homeomorphism
$h\:\tl X\to\tl X'$ that is fiberwise (i.e.\ $ph=p'$) and such that
$\{h(U_\alpha^\lambda)\}$ and $\{V_\beta^\lambda\}$ are equivalent as
overlay structures on $p'$.

\proclaim{Proposition 7.2} \cite{F2}, \cite{Mo} If $X$ is locally connected,
every covering over $X$ admits a unique overlay structure.
\endproclaim

Moore's argument is helpful for comprehending the above definitions, so we
reproduce it.

\demo{Proof \cite{Mo}} If $\{V_\beta\}$ is a cover of $X$ by connected sets, any
two overlay structures on the given covering map $p$ that lie over $\{V_\beta\}$
will clearly coincide up to a renumbering.
To construct such a $\{V_\beta\}$, notice that every $x\in X$ has a connected
neighborhood, since by the hypothesis the connected component of $X$
containing $x$ contains a neighborhood of $x$.
Hence every cover $\{U_\alpha\}$ of $X$ has a refinement $\{V_\beta\}$
where each $V_\beta$ is connected.
If $\{U_\alpha^\lambda\}$ is an overlay structure lying over $\{U_\alpha\}$,
condition (ii) ensures that it has a refinement $\{V_\beta^\lambda\}$
lying over $\{V_\beta\}$.
Thus any two overlay structures on a given covering map over $X$ are equivalent.

On the other hand, suppose that $\{U_\alpha^\lambda\}$ lies over
$\{U_\alpha\}$, and $\{V_\beta\}$ is a refinement of $\{U_\alpha\}$.
If each $V_\beta$ is connected, there exists a refinement
$\{V_\beta^\lambda\}$ of $\{U_\alpha^\lambda\}$ lying over $\{V_\beta\}$.
If additionally $V_\beta\cap V_\gamma\ne\emptyset$ implies
$V_\beta\cup V_\gamma\i U_\alpha$ for some $\alpha$, then
$\{V_\beta^\lambda\}$ will be an overlay structure.
To construct such a $\{V_\beta\}$, let $\{W_\gamma\}$ be any cover of $X$ by
sets of diameters at most $\frac\lambda2$, where $\lambda$ is a Lebesgue number
for $\{U_\alpha\}$ (that is, every subset of $X$ of diameter at most $\lambda$ is
contained in at least one $U_\alpha$).
Then $W_\beta\cap W_\gamma=\emptyset$ implies $W_\beta\cup W_\gamma\i U_\alpha$
for some $\alpha$.
Finally $\{V_\beta\}$ is defined to be any refinement of $\{W_\gamma\}$
with connected elements. \qed
\enddemo

There exists a covering over $\N\x I\cup [0,\infty]\x\partial I$ that admits
no overlay structure \cite{F2; Figs.\ 1, 2 and p.\ 78}.
This is especially easy to see using Lemma 7.3(a), cf\. \cite{Mo; Example 1}.

If $X$ is the one-dimensional Hawaiian snail in Example 5.7, that is the
mapping torus of the self-embedding of $\{0\}\cup\{\frac 1{2^n}\mid n\in\N\}$
defined by $x\mapsto\frac x2$, its universal covering admits uncountably many
pairwise inequivalent overlay structures.
However, they are all isomorphic to each other via covering transformations
which are not self-isomorphisms of any overlay structure (compare
\cite{F1; Example 3}).
Using this and Lemma 7.3(b), it is easy to construct a covering map over
$S^1\vee X$ (where the basepoint of $X$ is not in the limit circle)
that admits uncountably many pairwise non-isomorphic overlay structures
\cite{Mo; Example 3}.
(Beware that \cite{Mo; Example 2} is erroneous, as can be seen using either
Lemma 7.3(a) or Corollary 7.5 or Theorem 7.6.)

The following lemma simplifies \cite{F1; Extension Theorem 5.2},
\cite{Mo; Theorem 1.5}.

\proclaim{Lemma 7.3} (a) A covering map over a compactum admits an overlay
structure iff it is induced from a covering map over a polyhedron.

(b) Overlayings $f^*(p)$, $g^*(q)$ over a compactum $X$ induced from
the covering maps $p$, $q$ over polyhedra $P$, $Q$ are isomorphic iff
$f$ and $g$ factor up to homotopy as $f\:X@>h>>R@>f_1>>P$ and
$g\:X@>h>>R@>g_1>>Q$, where $R$ is a polyhedron and $f_1^*(p)=g_2^*(q)$.
\endproclaim

\demo{(a)} Clearly, every covering map induced from an overlaying comes with
the induced overlay structure.
Conversely, suppose we are given a covering $p$ over a compactum $X$ with
an overlay structure $\{U_\alpha^\beta\}$ lying over $\{U_\alpha\}$.
If $U_\alpha^\lambda\cap U_{\beta_1}^{\mu_1}\cap\dots\cap U_{\beta_r}^{\mu_r}$
and $U_\alpha^\lambda\cap U_{\beta_1}^{\nu_1}\cap\dots\cap U_{\beta_r}^{\nu_r}$
are both nonempty, each $\mu_i=\nu_i$.
Therefore the obvious map $q$ of the nerve of $\{U_\alpha^\lambda\}$ onto
the nerve $N$ of $\{U_\alpha\}$ is a covering map.
Thus $p$ is induced from the covering map $q$ by the projection $X\to N$. \qed
\enddemo

\demo{(b)} Let $\{U_\alpha\}$ and $\{V_\beta\}$ be the covers of $P$ and
$Q$ by the open stars of the vertices of some triangulations.
Each $U_\alpha\cap U_\gamma$ and each $V_\beta\cap V_\delta$ are connected,
so $p$ and $q$ admit overlay structures lying over $\{U_\alpha\}$ and
$\{V_\beta\}$.
Let $\{W_\gamma^\lambda\}$ be an overlay structure on the covering map
$f^*(p)=g^*(q)$ that is a common refinement of the induced overlay structures.
If $\{W_\gamma^\lambda\}$ lies over $\{W_\gamma\}$, then $f$ and $g$ factor
up to homotopy through the nerve $N$ of $\{W_\gamma\}$, and the maps
$N\to P$ and $N\to Q$ induce the same covering map from $p$ and $q$. \qed
\enddemo

The following result represents a simplification of Fox's formulation of
the classification of overlayings \cite{F1}.

\proclaim{Theorem 7.4}
Let $X$ be a connected compactum and $d\in\{1,2,\dots,\infty\}$.
There exists a natural bijection between $[X,\,BS_d]$ and the set of isomorphism
classes of $d$-fold overlayings over $X$.
\endproclaim

Note that when $X$ is a polyhedron, this is a mere restatement of the polyhedral
case of Theorem 7.1.

\demo{Proof} By Lemma 7.3(a), every overlaying over $X$ is induced from
a covering over a polyhedron $P$ by a map $f\:X\to P$.
This yields a map $X\to P\to BS_d$, whose homotopy class is well-defined
by Lemma 7.3(b).
Conversely, if $X$ is an inverse limit of polyhedra $P_i$, by Lemma 2.1(a)
every map $f\:X\to BS_d$ extends to $P_{[i,\infty]}$ for some $i$.
Hence it factors up to homotopy through $X\to P_i$ and so gives rise to
an overlaying.
If $f$ is homotopic to a $g\:X\to BS_d$, which extends to $P_{[j,\infty]}$,
by Lemma 2.1(b)
the two extensions are homotopic over $P_{[k,\infty]}$ for some $k$ and
so determine the same overlaying.
The assertion now follows from Proposition 2.3(a). \qed
\enddemo

By Proposition 2.3(b), $[X,\,BS_d]=\dirlim [P_i,\,BS_d]$ where the compactum
$X$ is the inverse limit of the polyhedra $P_i$.
On the other hand, $\dirlim [P_i,\,BS_d]$ can be identified with the orbit
set of the action of $S_d$ on $\dirlim\Hom(\pi_1(P_i),\,S_d)$ generated by
inner automorphisms of $S_d$.
It is well-known (and not hard to see) that Mittag-Leffler inverse sequences
of groups reduce to their topological inverse limits \cite{AS; \S2}.
This implies the following generalization of Theorem 7.1.

\proclaim{Corollary 7.5} \cite{He} Let $X$ be a Steenrod connected compactum,
and let $d\in\{1,2,\dots,\infty\}$.
The monodromy map yields a bijection between the isomorphism classes of
$d$-fold overlayings over $X$ and the conjugacy classes of representations of
the topological group $\check\pi_1(X)$ (or $\pi_1(X)$) into the discretely
topologized symmetric group $S_d$.
\endproclaim

Corollary 7.5 fails for arbitrary compacta: there is a free action of $\Z/3$
on the $2$-adic solenoid $\Sigma_2$ with orbit space again homeomorphic to
$\Sigma_2$; however $\check\pi_1(\Sigma_2)=0$.

Note that the addendum on transitive representations, included in Theorem 7.1,
is dropped in Corollary 7.5: in fact, it is false in this generality, by
considering the one-dimensional Hawaiian snail (see  \cite{F1; Example 3},
compare \cite{F2}).
It can be shown to hold, however, if ``connected'' is understood in the sense
of uniform spaces rather than topological spaces:

\proclaim{Theorem 7.6} Let $X$ be a compactum.
There is a natural bijection between the set of equivalence classes of overlay
structures on a covering map $p\:\tl X\to X$ and the set of uniformities on
$\tl X$ agreeing with $p$.
\endproclaim

We say that a uniformity (=uniform structure) on $\tl X$ {\it agrees} with
a covering map $p\:\tl X\to X$, where $X$ is a compactum, if every $x\in X$ has
a closed neighborhood $U$ such that

\medskip
(iii) $p^{-1}(U)=\bigsqcup_\lambda U^\lambda$ as uniform spaces, where
each $p|_{U^\lambda}$ is a homeomorphism.
\medskip

Thus Fox's conditions (i) and (ii) are effectively replaced by the single
condition (iii).
An equivalent way of formulating (iii) is:
\medskip

(iii$'$) $p^{-1}(U)$ is uniformly homeomorphic to $U\x F$,
where $F$ is discrete, by a homeomorphism $h$ such that
$ph^{-1}\:U\x F\to U$ is the projection.
\medskip

We emphasize that $F$ is assumed to be discrete as a uniform space, which
is a stronger condition than the discreteness of its underlying topology.
For instance, $\{\log 1,\log 2,\dots\}\i\R$ is discrete as a topological subspace
of $\R$, but not discrete as a uniform subspace of $\R$ since
$(\log (n+1)-\log n)\to 0$ as $n\to\infty$.
Our reference for uniform spaces is Isbell's book \cite{Is}.

The above definition is equivalent to I. M. James' notion of covering in
the uniform category (see \cite{BDLM} where James' definition is slightly
corrected).
Theorem 7.6 may be viewed as an improvement of a result of \cite{Mo}, but
the author was led to it through the analysis of Example 5.8.

\demo{Proof} Suppose that $X$ is the inverse limit of compact polyhedra $P_i$.
We may assume that they are triangulated so that each bonding map is
simplicial as a map $P_{i+1}\to P_i'$ into some subdivision of the barycentric
subdivision of $P_i$.
By Lemma 7.3(a) (see also the proof of Theorem 7.4) we may assume that
the given overlaying $p$ is induced from a covering map $q_k$ over some $P_k$.
For each $j>k$ let $q_j$ be the covering map over $P_j$ induced from $q_k$, and
let $E_j$ be its total space.
Then $\tl X$ is the inverse limit of the uniformly continuous maps
$E_{j+1}\to E_j$, where each $E_j$ is endowed with the standard metric
corresponding to a triangulation of $E_j$ making $q_j$ simplicial.
This defines a uniformity on $\tl X$ agreeing with $p$, and by Lemma 7.3(a)
(see also the proof of Theorem 7.4) it does not depend on the choice of $q_k$.

Given a uniformity on $\tl X$ agreeing with $p$, let $\{U_\alpha\}$
be a cover of $X$ with each $U_\alpha$ satisfying (iii$'$).
Without loss of generality it is a finite cover.
Fix some metric on $X$, and let $d$ be the diameter of $X$.
Consider the metric on $F$ where the distance between any two points is $2d$.
Let us represent the uniformity on $\tl X$ by some metric.
Let $\delta>0$ be such that $\delta$-close points of $\tl X$ remain $d$-close
under each $h_\alpha\:p^{-1}(U_\alpha)\to U_\alpha\x F$, and let $\eps>0$
be such that $\eps$-close points of each $U_\alpha\x F$ remain $\delta$-close
under $h_\alpha^{-1}$.
Pick a refinement $\{V_\beta\}$ of $\{U_\alpha\}$ such that each
$\{V_\beta\}$ has diameter at most $\eps$.
If $V_\beta$ is contained in both $U_\alpha$ and $U_\gamma$, the restrictions
$V_\beta\x F@>h_\alpha^{-1}>>p^{-1}(V_\beta)@>h_\gamma>>V_\beta\x F$
have to be of the form $\id_{V_\beta}\x\sigma$ for some permutation
$\sigma\in S_F$.
Thus if $V_\beta^\lambda=h_\alpha(V_\beta\x\{\lambda\})$, the cover
$\{V_\beta^\lambda\}$ of $\tl X$ is well-defined.
If $\{V_\beta\}$ is chosen so that additionally
$V_\beta\cap V_\gamma\ne\emptyset$ implies $V_\beta\cup V_\gamma\i U_\alpha$
for some $\alpha$ (see the proof of Proposition 7.2), then $\{V_\beta^\lambda\}$
is an overlay structure.
It is clear that any two overlay structures defined in this way
are equivalent. \qed
\enddemo

\proclaim{Corollary 7.7}
A fiberwise homeomorphism of overlayings is an isomorphism iff it is a uniform
homeomorphism of the corresponding uniformities.
\endproclaim

\proclaim{Corollary 7.8} \cite{F2}, \cite{Mo}
Every finite-sheeted covering map over a compactum $X$ admits a unique overlay
structure.
\endproclaim

The remainder of this section is devoted to an outline of how overlaying spaces
of compacta can be included in the Steenrod homotopy category.
These constructions will be used sporadically in the next section (most
importantly in the proof of Theorem 8.9), and in turn assume the reader's
familiarity with metrizable uniform spaces (see \cite{Is}).

\definition{Uniform compactified mapping telescope}
If $X$ is a metrizable uniform space, the {\it cone} $CX:=X\x I/X\x \{0\}$
is given the quotient uniformity, which is metrizable (see \cite{Vi}) --- in
contrast to the quotient topology.
The {\it mapping cylinder} $MC(f)$ of a uniformly continuous map $f\:X\to Y$
between metrizable uniform spaces is the image of $\Gamma_f\x I\i X\x Y\x I$ in
the metrizable uniform space $CX\x Y$.
If $X$ is complete, then so is $CX$, and if additionally $Y$ is complete,
so is $MC(f)$.

Next, if $A$ is a closed subset of a metrizable uniform space $X$, the
quotient uniform space $X/A$ is metrizable (see \cite{Vi}).
It follows, by standard techniques of uniform absolute retracts (see \cite{Is},
where they are called ``injective spaces'') that every bounded metric on
$A$, inducing its subspace uniformity, extends to a bounded metric on $X$,
inducing its original uniformity.
(The similar extension lemma for pseudo-metrics is known \cite{Is; III.16}.)

The {\it gluing} $X\cup_h Y$ of metrizable uniform spaces along a uniform
homeomorphism $h$ between a closed subset of $X$ and a closed subset of $Y$ is
given the quotient uniformity $(X\sqcup Y)/h$, which is metrizable
(see \cite{BH; I.5.24}), and is well-defined by the above-mentioned extension
lemma.
If $X$ and $Y$ are complete, then so is $X\cup_h Y$.

Consider an inverse sequence $X=(\dots@>f_1>>X_1@>f_0>>X_0)$ of uniformly
continuous maps between metrizable uniform spaces.
Let $X_{[0,\infty]}$ be the inverse limit of the finite mapping telescopes
$X_{[0,n]}=MC(f_0)\cup_{X_1}MC(f_1)\cup_{X_2}\dots\cup_{X_{n-1}}MC(f_{n-1})$
and the obvious retractions $X_{[0,n+1]}\to X_{[0,n]}$.
Given a $J\i [0,\infty]$, by $X_J$ we denote the preimage of $J$ under the obvious
projection $X_{[0,\infty]}\to [0,\infty]$.
According to the above, if each $X_i$ is compete, then so is
$X_{[0,\infty]}$.
In that case it is the completion of the infinite mapping telescope
$X_{[0,\infty)}$.
\enddefinition

\definition{Convergent and Cauchy inverse sequences}
In the notation of the preceding paragraph, let us call the inverse sequence
$X$ {\it convergent} if every uniform
neighborhood of $\invlim X$ in $X_{\N\cup\infty}$ (or equivalently in
$X_{[0,\infty]}$) contains all but finitely many of $X_i$'s.
We say that $X$ is {\it Cauchy} if for every $\eps>0$ there exists a $k$
such that for every $j>k$, the $\eps$-neighborhood of $X_j$ in $X_{\N}$
(or equivalently in $X_{[0,\infty)}$) contains $X_k$.
Clearly, these notions depend only on the underlying uniform structures.
Here is an example of a Cauchy inverse sequence that is divergent:
$\dots\i(0,\frac14]\i(0,\frac12]\i(0,1]$.
The inverse sequence $\dots\i[2,\infty)\i[1,\infty)\i[0,\infty)$ fails to be
Cauchy.%
\footnote{The reader who is acquainted with Marde\v si\' c's ``resolutions''
can easily check that a convergent inverse sequence of uniformly continuous maps
between metrizable complete uniform spaces is a Marde\v si\'c resolution.}

The close analogy with the definition of a convergent/Cauchy sequence of points
in a metrizable uniform space can in fact be formalized.
If $M$ is a metric space, the space $2^M$ of all nonempty closed subsets of $M$
is endowed with the Hausdorff metric $d_H(A_1,A_2)=\min(\max(d_1,d_2),1)$, where
$d_i=\sup\{d(a,A_i)\mid a\in A_j,\,j\ne i\}$.
If $Y$ is a metrizable uniform space, the induced uniform structure of $2^Y$
is well-defined, and if $Y$ is complete or compact, so is $2^Y$
\cite{Is; II.48, II.49}.
We obtain that the inverse sequence $X$ is convergent (Cauchy) iff the
sequence of the closed subsets $X_i$ of $X_{\N\cup\infty}$ is convergent (Cauchy)
as a sequence of points in $2^{X_{\N\cup\infty}}$, or equivalently in
$2^{X_{[0,\infty]}}$.
This implies parts (a) and (b) of the following lemma.
\enddefinition

\proclaim{Lemma 7.9} Let $X=(\dots@>f_1>>X_1@>f_0>>X_0)$ be an inverse sequence
of uniformly continuous maps between metrizable uniform spaces.

(a) If each $X_i$ is compact, $X$ converges.

(b) If $X$ converges, it is Cauchy; the converse holds when each $X_i$ is
complete.

(c) $X$ is convergent if and only if for each $i$, every uniform neighborhood of
$f^\infty_i(\invlim X)$ in $X_i$ contains all but finitely many of $f^j_i(X_j)$'s.

(d) $X$ is Cauchy if and only if for each $i$ and every $\eps>0$ there exists
a $k$ such that for every $j>k$, the $\eps$-neighborhood of $f^j_i(X_j)$
in $X_i$ contains $f^k_i(X_k)$.

(e) If each $X_i$ is uniformly discrete, $X$ converges if and only if it
satisfies the Mittag-Leffler condition.

(f) If $X$ converges and each $X_i$ is non-empty (resp\. uniformly connected),
then $\invlim X$ is non-empty (resp.\ uniformly connected).

(g) If $Y=(\dots@>f_1>>X_1@>f_0>>X_0)$ is another inverse sequence of uniformly
continuous maps between metric spaces, $f_i\:X_i\to Y_i$ are surjections
commuting with the bonding maps and $X$ is convergent (Cauchy), then $Y$ is
convergent (resp.\ Cauchy).
\endproclaim

\demo{Proof. (a)} We have already deduced (a) from general facts, but it
may be of interest to have an elementary proof.
If $X$ diverges, there exists an open neighborhood $U$
of $X$ in $X_{\N\cup\infty}$ whose complement has a nonempty intersection $C_i$
with $X_{\{i,i+1,\dots,\infty\}}$ for each $i$.
By Cantor's theorem, the compacta $C_i$ have a non-empty intersection.
It is contained in $\bigcap X_{\{i,i+1,\dots,\infty\}}=X$ and at the same time
in the complement of $U$, which is a contradiction. \qed
\enddemo

\demo{(c)} A cover of $\invlim X_{\{0,\dots,i\}}$ is uniform iff it can
be refined by the preimage of a uniform cover of some
$X_{\{0,\dots,i\}}$ (see \cite{Is}).
Hence a neighborhood of $X_\infty$ in $X_{\N\cup\{\infty\}}$ is uniform
iff it contains the preimage $V$ of some uniform neighborhood $U$ of
$f^\infty_i(X_\infty)$ in $X_{\{0,\dots,i\}}$ for some $i$.
But this $U$ contains $f^j_i(X_j)$ (which are subsets of $X_i\i X_{\{0,\dots,i\}}$)
for almost all $j>i$ iff $V$ contains $X_j$ for almost all $j>i$. \qed
\enddemo

\demo{(d)-(g)} Part (d) is proved similarly to (c).
Parts (e) and (f) follow from (c).
Part (g) follows from (c) and (d). \qed
\enddemo

\proclaim{Corollary 7.10 (Bourbaki's Mittag-Leffler Theorem)} Let $L$ be the limit
of an inverse sequence $X=(\dots@>f_1>>X_1@>f_0>>X_0)$ of uniformly continuous
maps between complete metrizable uniform spaces.
If each $f_i(X_{i+1})$ is dense in $X_i$, then $f^\infty_0(L)$ is dense in $X_0$.
\endproclaim

\demo{Proof} The hypothesis implies that each $f^j_i(X_j)$ is dense in $X_i$.
Hence by Lemma 7.9(d), $X$ is Cauchy.
Then by 7.9(b), $X$ converges.
If $U$ is a uniform neighborhood of $f^\infty_0(L)$, by 7.9(c) it contains
$f^j_0(X_j)$ for some $j$.
If $U$ is closed, this implies $U=X_0$.
Since $U$ is arbitrary, $f^\infty_0(L)$ is dense in $X_0$. \qed
\enddemo

Corollary 7.10 is known to generalize for inverse spectra rather than inverse
sequences (see Bourbaki's ``General Topology'') and, in another direction,
for continuous, rather than uniformly continuous bonding maps \cite{Ru}.
By a slight modification of an argument from \cite{Ru} we also obtain

\proclaim{Corollary 7.11 (Baire's Category Theorem)}
The intersection of countably many dense open sets in a complete metrizable
uniform space is dense.
\endproclaim

\demo{Proof} Without loss of generality, the dense open sets $U_i$ are ordered
by inclusion: $U_1\supset U_2\supset\dots$.
If $d$ is a metric on $U_0:=X$, let $d_0=d$, and define a metric on $U_{i+1}$
by $d_{i+1}(x,y)=d_i(x,y)+
|\frac1{d_i(x,\,U_i\but U_{i+1})}-\frac1{d_i(y,\,U_i\but U_{i+1})}|$.
Then $d_{i+1}$ is a complete metric on $U_{i+1}$ such that
$\id\:(U_{i+1},d_{i+1})\to(U_{i+1},d_i)$ is a uniformly continuous homeomorphism.

The assertion follows from the Mittag-Leffler theorem applied to the inverse
sequence $\dots\to(U_1,d_1)\to(U_0,d_0)$. \qed
\enddemo

\definition{Steenrod homotopy category extended}
We call a uniform space $X$ {\it residually finite-dimensional} if every uniform
cover of $X$ has a refinement of finite multiplicity (that is, with
a finite-dimensional nerve); such spaces are called ``finitistic'' in
\cite{SSG}.
A metrizable uniform structure is defined on a finite-dimensional uniform
polyhedron in \cite{Is}; it is complete.
Residually finite-dimensional metrizable complete uniform spaces are precisely
the limits of convergent inverse sequences of finite-dimensional uniform
polyhedra and uniformly continuous bonding maps \cite{Is; V.33}.
In particular, these include overlaying spaces over compacta.

Let $X$ and $Y$ be the limits of convergent inverse sequences
$\dots@>p_1>>P_1@>p_0>>P_0$ and $\dots@>q_1>>Q_1@>q_0>>Q_0$
of finite-dimensional uniform polyhedra and uniformly continuous bonding maps.
Without loss of generality, $P_0=pt$ and $Q_0=pt$.
Then similarly to the proof of Lemma 2.1, every uniformly continuous map
$f\:X\to Y$ extends to a uniformly continuous map
$f_{[0,\infty]}\:P_{[0,\infty]}\to Q_{[0,\infty]}$ sending $P_{[0,\infty)}$ into
$Q_{[0,\infty)}$.
Moreover, every two such extensions are homotopic by a uniformly continuous
(as a map $P_{[0,\infty]}\x I\to Q_{[0,\infty]}$)
homotopy $\rel\, X$, sending $P_{[0,\infty)}\x I$ into $Q_{[0,\infty)}$.
We recall that we agreed to call a map $F\:P_{[0,\infty)}\to Q_{[0,\infty)}$
{\it semi-proper} if for each $k$ there exists an $l$ such that
$f^{-1}(Q_{[0,k]})\i P_{[0,l]}$.

We define a {\it Steenrod homotopy class} $X\ssm Y$ to be the semi-proper
homotopy class of a semi-proper map $f\:P_{[0,\infty)}\to Q_{[0,\infty)}$.
By the virtue of the semi-proper homotopy equivalences
$(\id_X)_{[0,\infty)}$ and $(\id_Y)_{[0,\infty)}$ this definition does not
depend on the choice of the inverse sequences $P$ and $Q$.
Steenrod homotopy groups of a residually finite-dimensional metrizable complete
uniform space are now defined similarly to \S3 (a special case was considered
in the Introduction).
Also, similarly to \S4 one can define Steenrod homology and Pontryagin
cohomology of such spaces, but we shall not need them for the purposes of
this paper.
\enddefinition

\definition{Universal generalized overlaying}
Let $X$ be a Steenrod connected compactum (the general case of a connected
compactum will be a recurring theme throughout the next section).
If $X$ is the limit of an inverse sequence of compact connected polyhedra $P_i$,
let $\tl X=\invlim(\dots\to\tl P_1\to\tl P_0)$, where $\tl P_i$ are their
universal covering spaces.
They come endowed with uniform structures (inherited from $P_i$'s via
the covering maps).
Since $P_i$ converge, so do $\tl P_i$.
It follows that the Steenrod homotopy type of the residually finite-dimensional
metrizable complete uniform space $\tl X$ does not depend on
the choice of $P_i$'s.
With some care, it can be shown that even the fiberwise uniform homeomorphism
class of $\tl X$ does not depend on this choice \cite{LaB} (an alternative
proof should appear as a byproduct in a subsequent paper by the present author).
By definition, $p\:\tl X\to X$ is a Steenrod fibration; it is also a Serre
fibration (see condition (ii) in the discussion preceding the statement of
Theorem 3.15).
Each fiber $p^{-1}(pt)$ is uniformly $0$-dimensional, i.e.\ is an inverse limit
of (countable) discrete uniform spaces.
Obviously, $\pi_0(\tl X)=pt$ and $\pi_1(\tl X)=1$.
If $\check\pi_1(X)$ happens to be discrete, then $p$ is an overlaying, which
may be called the {\it universal overlaying} of $X$.
\enddefinition

\example{Example 7.12 (the Nottingham compactum)} If $R$ is a commutative
associative ring with a unit, the set of formal power series
$f(x)=x+a_2x^2+a_3x^3+\dots$ with coefficients $a_i\in R$ forms a group
under composition: $fg(x)=f(g(x))$.
This group $N(R)$ acts effectively and $R$-linearly on the formal power series
ring $R[[x]]$.
In the case $R=\Z_p$ (the $p$-adic integers) $N(R)$ is known to algebraists
as the {\it Nottingham group}, and in the case $R=\Z$ it is related to
the algebra of stable cohomology operations in complex cobordism
\cite{BSh}, \cite{Bu}.
If $R$ is a topological ring, by endowing $N(R)$ with the topology of
the product $\prod_{i=2}^\infty R$ we make it into a topological group.
Let $N_n(R)$ be the subgroup of $N(R)$ consisting of all power series with
$a_2=\dots=a_n=0$, and let $N^n(R)$ denote the quotient $N(R)/N_n(R)$.
Clearly, $N(R)$ is the inverse limit of the truncating maps
$\dots\to N^2(R)\to N^1(R)=1$.

Now let us consider the topological quotient groups $X=N(\R)/N(\Z)$ and
$P_n=N^n(\R)/N^n(\Z)$.
Since the bonding maps in $\dots\to N^2(\Z)\to N^1(\Z)$ are surjective,
$X$ is the inverse limit of $\dots\to P_2\to P_1$.

The universal covers of $P_i$ are given by $N^i(\R)$, hence $N(\R)$ is a
universal generalized overlay of $X$.
Since $p\:N(\R)\to X$ is a fibration and $N(\R)$ is obviously contractible,
it follows that $\spi_1(X)\simeq N(\Z)$ (as topological groups) and
$\spi_i(X)=0$ for $i>1$.
(This was originally observed in \cite{BB}.)
Since $p\:N(\R)\to X$ is also a Steenrod fibration, we similarly conclude
(compare Theorem 3.15(b)) that the Steenrod homotopy groups of $X$ are
the same; moreover, $\stau\:\spi_i(X)\to\pi_i(X)$ is an isomorphism
for all $i$.

From the algebraic viewpoint, this unexpected coincidence of singular
and Steenrod homotopy is only made possible by the fact that
the group $N(Z)=\pi_1(X)$ is residually nilpotent.%
\footnote{and by the asphericity of $X$ --- see Example 8.6 and the subsequent
remark.}

Indeed, let us consider the following embedding of the Hawaiian earring
$E$ (see \S5) into $X$.
Since each truncating map $p_n\:P_{n+1}\to P_n$ is a bundle with fiber
$S^1$, we can embed the $(n+1)$-fold wedge of circles $Q_{n+1}$ into $P_{n+1}$
by taking the union of some cross-section of $p_n$ over the $n$-fold wedge
$Q_n$ and the fiber of $p_n$ over the basepoint of that wedge.
This yields an embedding of the inverse limit $E$ of $Q_n$'s into $X$.

Recall from Example 5.6, that the image of $\stau\:\spi_1(E)\to\pi_1(E)$
does not contain, for instance, the class of the following Steenrod loop
$\ell$:
$$(\dots,[a_1,a_2][a_1,a_3][a_1,a_4],\ [a_1,a_2][a_1,a_3],\ [a_1,a_2])\in
\invlim(\dots\to F_4\to F_3\to F_2).$$
This class is not representable by a map $S^1\to E$ simply because such a
map cannot wind infinitely many times around the first circle of $E$,
whereas the letter $a_1$ is read off by the Pontryagin--Thom construction
precisely on each pass of the first circle in the positive direction.

What is the class of $\ell$ in $\pi_1(X)$?
It is not hard to check that (in the group $N(R)$ with arbitrary $R$) if
$\alpha=x+ax^n+\dots$ and $\beta=x+bx^m+\dots$, then
$\alpha^{-1}\beta^{-1}\alpha\beta=x+ab(m-n)x^{m+n-1}+\dots$,
where the dots represent terms of higher degrees, cf.\ \cite{Jo}.
\footnote{This equality immediately implies that each $N^n(R)$ is
nilpotent.}
Thus, due to the relations in $N(\Z)$, the commutators $[a_1,a_i]$ find
themselves in increasingly deep subgroups $N_i(\Z)$ of this group and
hence are represented by increasingly short (exponentially shortening) loops.
So their product is representable by a loop, in agreement with our expectations.
\endexample

\remark{Remark}
Let us note incidentally that the universal generalized overlay $\tl E$ of
the Hawaiian Earring is not path-connected.
Indeed, let $F$ be the fiber of the Steenrod fibration $\tl E\to E$.
In the notation of the preceding example, the image of $[\ell]$ in $\pi_0(F)$
(see Theorem 3.15(b)) is representable by a map $t\:S^0\to F$ due to the
$0$-dimensionality of $F$.
Now the image of $[t]$ in $\spi_0(\tl E)$ is nontrivial by virtue of
our choice of $\ell$ and using that $\tl E$ is Steenrod simply-connected.
\endremark

\head 8. Zero-dimensional homotopy \endhead

In the zero-dimensional case, Theorem 3.12 admits the following
modification:

\proclaim{Proposition 8.1} {\rm (Krasinkiewicz \cite{K2; 1.1})} Let $X$ be
the limit of an inverse sequence of compact polyhedra $\dots\to P_2\to P_1$.
Then $X$ is Steenrod connected iff for each $k$ there exists a $j>k$
such that for every $i>j$, every path $(I,\partial I)\to (P_{[j,\infty]},X)$ is
homotopic $\rel\,\partial I$ and with values in $P_{[k,\infty]}$ to a path in
$P_{[i,\infty]}$.
\endproclaim

Krasinkiewicz's original proof was geometric.

\demo{Proof} By Theorem 3.1 and Lemma 3.3, $X$ is Steenrod connected iff
$\pi_1(P_i)$ satisfy the Mittag-Leffler condition.
By Lemma 3.11, this is the case iff $\spi_1(P_{[i,\infty]},X)$ satisfy
the Mittag-Leffler condition.
A homotopy $h$ of a path can be converted to one that keeps the second endpoint
fixed by using the restriction of $h$ to that endpoint. \qed
\enddemo

\proclaim{Theorem 8.2} {\rm (McMillan \cite{Mc}; Krasinkiewicz \cite{DS1; 7.2.4},
\cite{K2})}
A continuous image of a Steenrod connected compactum is Steenrod connected.
\endproclaim

At the same time, a continuous surjection between compacta need not induce
an epimorphism on $\pi_0$, by considering the projection
$\Sigma_p\sqcup\Sigma_p\to\Sigma_p\vee\Sigma_p$ (see Examples 5.3 and 5.5).

The proof below simplifies McMillan's approach by using Proposition 8.1.
It is also much easier than either of the two proofs by Krasinkiewicz.
Yet another proof, involving non-metrizable compacta, is given in the
third paper of the series \cite{KO}.

\demo{Proof} Let $f\:X\to Y$ be the given surjection.
By the proof of Theorem 3.15(a), we may assume that $X=\invlim P_i$,
$Y=\invlim Q_i$ and $f=\invlim f_i$ for some compact polyhedra $P_i$ and $Q_i$
and PL maps $f_i\:P_i\to Q_i$ commuting with the bonding maps.
Pick a point $v\in Y$ and let $v_k$ be its image in $Q_k$.
Choose polyhedral neighborhoods $N_k(v)$ of $v_k$ in $Q_k$ so that each $N_k(v)$
contains the image of $N_{k+1}(v)$ and so that $\invlim N_i(v)=\{v\}$.
Let $M_k(v)=f_k^{-1}(N_k(v))$, then $\invlim M_i(v)=f^{-1}(v)$; in particular,
we may assume that each $M_i(v)$ is nonempty for sufficiently large $i=i(v)$.

By the Ferry Construction for L$\sC_{-1}$ compacta (see the beginning of \S6),
for each $k$ there exist a $j$ and an $n$ such that a path
$\ell\:(I,\partial I)\to (Q_{[j,\infty]},Y)$ is homotopic $\rel\partial$ with
values in $Q_{[k,\infty]}$ to an $\ell'$ such that $y_i:=\ell'(\frac in)\in Y$
for each $i=0,\dots,n$, and $q^\infty_k(y_i)\in N_k(y_{i+1})$ for each $i<n$.
Then there exists an $h>k$ such that $q^h_k(N_h(y_i))\i N_k(y_{i+1})$ for each
$i<n$.
Pick some $x_i\in M_h(y_i)$ for each $i=0,\dots,n$.
Note that $p^h_k(x_i)\in M_k(y_{i+1})$.
By the Ferry Construction for L$\sC_{-1}$ compacta, given an $l$, we may choose $k$
so that $x_i$ is joined to some $x_i^-\in f^{-1}(y_i)$ by a path in
$M_{[l,\infty]}(y_i)$ and to some $x_i^+\in f^{-1}(y_{i+1})$ by a path in
$M_{[l,\infty]}(y_{i+1})$.
Let $\ell_i\:(I,\partial I)\to(P_{[0,\infty]},X)$ be the product of these two
paths, connecting $x_i^-$ with $x_i^+$.
Since $Q_{[0,\infty]}$ is L$\sC_\infty$, given an $m$, we may choose $l$ so that
$f_{[0,\infty]}(\ell_i)$ and $\ell'|_{[\frac in,\frac{i+1}n]}$ are homotopic
$\rel\partial$, with values in $Q_{[m,\infty]}$.
On the other hand, by Proposition 8.1, we may choose $l=l(m)$ so that for each $t$,
each $\ell_i$ is homotopic $\rel\partial$ with values in $P_{[m,\infty]}$ to
a path in $P_{[t,\infty]}$. \qed
\enddemo

\proclaim{Theorem 8.3} {\rm (Krasinkiewicz \cite{K2})} A connected union of two
Steenrod connected compacta is Steenrod connected.
\endproclaim

This assertion is a special case of Corollary 3.14; nevertheless, since it
will be essentially used below, we include a slightly different, more
direct proof.

\demo{Proof} We may represent $X$ as $\invlim(\dots\to P_2\to P_1)$, where
each $P_i$ is a union of two compact polyhedra $Q_i$ and $R_i$ such that
$Y=\invlim Q_i$ and $Z=\invlim R_i$.
By the Ferry Construction for L$\sC_{-1}$ compacta (see the beginning of \S6),
for each $k$ there exists a $j$ such that every path
$l\:(I,\partial I)\to (P_{[j,\infty]},X)$ is homotopic $\rel\,\partial I$ and
with values in $P_{[k,\infty]}$ to a product of paths $\ell_i$ alternating
between $(Q_{[k,\infty]},Y)$ and $(R_{[k,\infty]},Z)$.
Given an $l$, by Proposition 8.1, $k>l$ can be chosen so that for an arbitrary
$h>k$, each $\ell_{2i+1}$ is homotopic $\rel\,\partial I$ and with values in
$Q_{[l,\infty]}$ to a path in $Q_{[h,\infty]}$; and each $\ell_{2i}$ is
homotopic $\rel\,\partial I$ and with values in $R_{[l,\infty]}$ to a path
in $R_{[h,\infty]}$.
If we have infinitely many of $\ell_i$'s, only finitely many of the homotopies
have to be non-identical.
Hence $\ell$ is homotopic $\rel\,\partial I$ and with values in $P_{[l,\infty]}$
to a path in $P_{[h,\infty]}$. \qed
\enddemo

\remark{Remark}
The proof of Proposition 8.1 shows (see the remark to Lemma 3.11) that the images
of $\pi_1(P_i)$ in $\pi_1(P_k)$ stabilize if and only if the images of
$\spi_1(P_{[i,\infty]},X)$ in $\spi_1(P_{[k,\infty]},X)$ stabilize (with
the same $k$).
In this event, we call the connected compactum $X$ {\it Steenrod connected
over} $P_k$ (or in more detail, {\it over} $p^\infty_k\:X\to P_k$).
It is easy to see that this property does not depend on the choice of
$P_{k+1},P_{k+2},\dots$.

By Theorem 3.1(b) and Lemma 3.3, a compactum is Steenrod connected
if and only if it is Steenrod connected over every map to a polyhedron.
By the proof of Theorem 8.3, if $X$ is a union of compacta $Y$ and $Z$ and
$f\:X\to P$ is a map to a polyhedron such that $Y$ and $Z$ are Steenrod connected
over $f|_Y$ and $f|_Z$, then $X$ is Steenrod connected over $f$.
\endremark
\medskip

An inverse sequence of groups $G_i$ is said to be {\it virtually
Mittag-Leffler} if for each $i$ there exists a $j>i$ such that for each $k>j$
the image of $G_k\to G_i$ has finite index in that of $G_j\to G_i$.
If $X$ is the limit of an inverse sequence of compact connected polyhedra
$P_i$, satisfaction of the virtual Mittag-Leffler condition for the inverse
sequence $\pi_1(P_i)$ clearly does not depend on the choice of the basepoint of
$X$ and is an invariant of the equivalence relation from Proposition 2.6(iii).
Thus we may call $X$ {\it virtually Steenrod connected} if $\pi_1(P_i)$
are virtually Mittag-Leffler; moreover, the property of being virtually
Steenrod connected is an (unpointed) shape invariant of $X$.

For example, $\dots@>f>>F_2@>f>>F_2$ is not virtually Mittag-Leffler, where
$F_2$ is the free group $\left<x,y\mid\,\right>$ and $f(x)=x^p$, $f(y)=y$.
In particular, the wedge of the $p$-adic solenoid and $S^1$ is not
virtually Steenrod connected.

\proclaim{Theorem 8.4} {\rm (Geoghegan--Krasinkiewicz \cite{GK})} If $X$ is
a connected compactum that is not virtually Steenrod connected, then
$\spi_0(X)@>\sstau>>\pi_0(X)$ is not surjective.
\endproclaim

We modify the original proof in \cite{GK}, replacing certain superfluous
constructions with an algebraic analysis of compactified mapping telescopes,
which will be useful for the rest of this section.

\demo{Proof} Suppose that $X$ is an inverse limit of compact connected
polyhedra $P_i$ and let $G_i=\pi_1(P_i)$.
By the hypothesis, there exists a $k$ such that the images $A_i$ of $G_{k+i}$
in $G_k$ do not virtually stabilize, i.e.\ the pointed sets $A_i/A_{i+1}$ are
infinite for infinitely many values of $i$.
Then if $A:=\invlim A_0/A_i$ is endowed with the inverse limit topology
(where each $A_0/A_i$ is discrete), every compactum in $A$ is nowhere
dense in $A$ (the details are similar to those in the proof of Lemma 3.7(b)).
Since $A$ is complete, by Baire's Category Theorem 7.11, $A$ is not a countable
union of compacta.
\enddemo

The remainder of the proof is presented in two ways.

\demo{Simple argument (using Steenrod homotopy of uniform spaces, see \S7)}
If $K$ is a compactum, $\spi_0(K)$ is compact since $K$ is a continuous image
of the Cantor set.
Let $\tl X\to X$ be the pullback of the universal covering $\tl P_k\to P_k$.
Then $\tl X$ is a countable union of compacta.
In particular, $\spi_0(\tl X)$ is a countable union of compacta.

Let $\tl P_{k+i}$ be the pullback of $\tl P_k$ over $P_{k+i}$.
Then the inverse sequence $\dots\to\tl P_1\to\tl P_0$ converges, with
inverse limit $\tl X$.
It is easy to see that $\pi_0(\tl P_{k+i})$ can be identified with the $A_0$-set
$A_0/A_i$, and therefore $\check\pi_0(\tl X):=\invlim\pi_0(\tl P_i)$ is
$A_0$-homeomorphic to $A$.
Since $\check\tau\:\pi_0(\tl X)\to\check\pi_0(X)$ is a surjection similarly
to Theorem 3.1(b), $\pi_0(\tl X)$ is not a countable union of compacta.

Hence $\stau\:\spi_0(\tl X)\to\pi_0(\tl X)$ is not surjective.
But from the homotopy exact sequences of a Steenrod fibration (similarly to
Theorem 3.15(b)) and fibration, this map is an $A_0$-equivariant lift of
$\stau\:\spi_0(X)\to\pi_0(X)$.
Thus the latter is not a surjection. \qed
\enddemo

\demo{Elementary argument} The basepoint $x\in X$ gives rise to the base ray
$p_{[k,\infty)}$ in $P_{[k,\infty)}$.
Since each $P_{[i,\infty)}$ deformation retracts onto $P_i$, each $A_0/A_i$
can be identified, as an $A_0$-set, with
$\pi_1(P_{[k,\infty)},P_{[i,\infty)};\,p_k)$.
Similarly to the proof of Theorem 3.1(b), this yields a surjection
$\phi\:G\to A$, where $G$ is the pointed set of proper homotopy classes of
proper maps $([0,\infty),0)\to (P_{[k,\infty)},p_k)$.%
\footnote{By Lemma 3.13, $G$ can be identified with $\pi_1(P_{[k,\infty]},X)$.}
Then the orbit space $A_0\q G$ of the left action of $A_0$ on $G$ can be
identified with $\pi_0(X)$.%
\footnote{In fact, it follows from the proof of Lemma 3.3 that $\phi$ is an
$A_0$-equivariant lift of the surjection $\derlim G_i\to\derlim A_i$ induced
by the surjections $G_i\to A_i$.}
From the homotopy exact sequence of a pair, $\spi_0(X;\,x)$ is the orbit space
$A_0\q H$ of the left action of $A_0\simeq\spi_1(P_{[k,\infty]};\,x)$ on
$H:=\spi_1(P_{[k,\infty]},X;\,x)$.
Using the homotopy $\Pi_t$ from \S2, every $h\in H$ can be represented by
a map $[-\infty,\infty]\to P_{[k,\infty]}$ sending $[-\infty,0]$
homeomorphically onto $p_{[k,\infty]}$ and $[0,\infty)$ into $P_{[k,\infty)}$;
and any two such representatives of $h$ are homotopic through such maps.
This yields an $A_0$-equivariant map $\psi\:H\to G$ such that the map
of the orbit spaces $A_0\q H\to A_0\q G$ can be identified with
$\stau\:\spi_0(X)\to\pi_0(X)$.
In particular, if $\stau$ is surjective, so is $\psi$.
Let us topologize $H$ and $G$ by declaring basic open sets to be
the point-inverses of the maps to
$\spi_1(P_{[k,\infty]},P_{[i,\infty]})\simeq\pi_1(P_{[k,i]},P_i)$,
resp.\ to $\spi_1(P_{[k,\infty)},P_{[i,\infty)})\simeq\pi_1(P_{[k,i]},P_i)$.
Then $\psi$ and $\phi$ are continuous.

Let $T_1,T_2,\dots$ be some triangulations of $P_1,P_2,\dots$ such that
each bonding map $P_{k+1}\to P_k$ is simplicial with respect to $T_{k+1}$
and some subdivision of the barycentric subdivision of $T_k$.
If $Q_k$ denotes the dual $0$-skeleton of $T_k$, by associating to each
top-dimensional simplex $\sigma$ of $T_{k+1}$ the simplex of $T_k$ containing
the image of $\sigma$, we obtain a surjection $Q_{k+1}\to Q_k$ with finite
point-inverses.
The inverse limit of $Q_k$'s is a Cantor set $C$, and the inclusion
$Q_{[0,\infty)}\i P_{[0,\infty)}$ extends to a continuous map
$Q_{[0,\infty]}\to P_{[0,\infty]}$ sending $C$ onto $X$.
Moreover, writing $K=\spi_1(Q_{[k,\infty]}\cup P_k,\,C)$, it is not hard to
see that the induced continuous map $\chi\:K\to H$ is onto.
Finally, $K$ is the inverse limit of the discrete spaces
$\pi_1(Q_{[k,i]}\cup P_k,\,Q_i)$, where the bonding maps have finite
point-inverses.
Thus $K$ is a countable union of compacta.
Since the composition $K@>\chi>>H@>\psi>>G@>\phi>>A$ is a continuous
surjection, $A$ must be a countable union of compacta, which is a
contradiction. \qed
\enddemo

\proclaim{Lemma 8.5}  {\rm (Krasinkiewicz--Geoghegan \cite{K3; p.\ 48},
\cite{GK; \S7})} Let $X$ be the limit of an inverse sequence of compact
connected polyhedra $P_i$ and maps $p_i\:P_{i+1}\to P_i$.
Let $b\in X$ be the basepoint and $b_i\in P_i$ its images.

The class $\alpha$ of $(g_0,g_1,\dots)\in\prod\pi_1(P_i)$ in
$\derlim\pi_1(P_i,b_i)=\pi_0(X,b)$ is represented by precisely those points
of $X$ that lie in the image of the projection of the inverse limit $W_\alpha$
of the universal covering spaces
$\dots@>\tl p_1>>\tl P_1@>g_1>>\tl P_1@>\tl p_0>>\tl P_0@>g_0>>\tl P_0$, where
each $\tl p_i$ sends a fixed lift $\tl b_{i+1}$ of $b_{i+1}$ into the previous
one $\tl b_i$.

In particular, $\alpha$ lies in the image of $\spi_0(X,b)$ if and only if
$W_\alpha$ is non-empty.
\endproclaim

We call the set of points of $X$ representing some $\alpha\in\pi_0(X)$
a {\it Steenrod component} of $X$.
In this terminology, Example 5.5 and Theorem 8.4 assert existence of empty
Steenrod components.

\demo{Proof}
Suppose that $(g_0,g_1,\dots)\in\prod\pi_1(P_i)$ represents the image in
$\pi_0(X,b)$ of some $[a]\in\spi_0(X,b)$, where $a\in X$ is a point.
Let us consider its images $a_i\in P_i$ and fix a path $w_i$ from $b_i$ to $a_i$.
Then $(g_0,g_1,\dots)$ represents the same element of $\derlim\pi_1(P_i,b_i)$
as $(g_0',g_1',\dots)$, where $g_i'$ is the class of the loop
$w_ip_{i+1}(\bar w_{i+1})$.
We may assume that $g_i'=g_i$ by appropriately redefining the paths $w_i$
proceeding from the definition of $\derlim$.
Consider the lift $\tl w_i$ of each $w_i$ starting at $\tl b_i$, and let
$\tl a_i$ be its other endpoint.
Then the composition $\tl P_{i+1}@>\tl p_i>>\tl P_i@>g_i>>\tl P_i$ sends
$\tl a_{i+1}$ precisely into $\tl a_i$.
Hence the inverse limit of these compositions contains the inverse limit of
the singletons $\{\tl a_i\}$, which is non-empty.

Conversely, given a point $\tl a$ in the inverse limit of the compositions
$\tl P_{i+1}@>\tl p_i>>\tl P_i@>g_i>>\tl P_i$, the preceding argument can be
reversed to show that the projection $a$ of $\tl a$ into $X$ represents
the class of $(g_0,g_1,\dots)$ in $\pi_1(X)$. \qed
\enddemo

\remark{Remark}
Similarly to Lemma 8.5 one can prove its homological version, which can be used,
for instance, to clarify the algebraic content of Examples 5.4 and 5.5.
Specifically, it says that the class $\alpha$ of
$(g_0,g_1,\dots)\in\prod H_1(P_i)$ in
$\derlim H_1(P_i)=H_0(X,b)$ is represented by precisely those points
of $X$ that lie in the image of the projection of the inverse limit $W_\alpha$
of the universal abelian covering spaces
$\dots@>\tl p_1>>\tl P_1@>g_1>>\tl P_1@>\tl p_0>>\tl P_0@>g_0>>\tl P_0$, where
each $\tl p_i$ sends a fixed lift $\tl b_{i+1}$ of $b_{i+1}$ into the previous
one $\tl b_i$.
In particular, $\alpha$ lies in the image of the composition
$\spi_0(X,b)\to\sH_0(X)\to H_0(X)$ if and only if $W_\alpha$ is non-empty.
\endremark
\medskip

The literal converse of Theorem 8.4 does not hold:

\example{Example 8.6 (Brin's compactum%
\footnote{This may or may not be the actual Brin's compactum mentioned in
\cite{GK; Remark 10.4}.}%
)}
Let $X\subset \Sigma_p\x\Sigma_p\x\Sigma_p$ be defined as the union
$[\Sigma_p\x\Sigma_p\x pt]\cup[\Sigma_p\x pt\x\Sigma_p]\cup
[pt\x\Sigma_p\x\Sigma_p]$.
Thus $X=\invlim(\dots@>f>>P@>f>>P)$, where
$P=[S^1\x S^1\x pt]\cup[S^1\x pt\x S^1]\cup[pt\x S^1\x S^1]$ and $f$
is the restriction of the self-map $p\x p\x p$ of $S^1\x S^1\x S^1$.
Since the inclusion $P\i (S^1)^3$ induces an isomorphism on $\pi_1$,
the inclusion $X\i (\Sigma_p)^3$ induces a bijection on $\pi_0$.
In particular, $X$ is virtually Steenrod connected.

Let $\alpha=\overline{\dots\alpha_1\alpha_0}\in\Z_p\but\Z$.
We claim that $g:=(\alpha+\Z,\,\alpha+\Z,\,\alpha+\Z)\in(\Z_p/\Z)^3\simeq\pi_0(X)$
is not in the image of $\spi_0(X)$.
Indeed, the universal cover $\tl P$ of $P$ can be identified with
the underlying space of the topological subgroup
$[\R\x\R\x\Z]\cup[\R\x\Z\x\R]\cup[\Z\x\R\x\R]$ of $\R^3$.
The intersection of the nested sequence of cosets
$$\dots\subset (\overline{\alpha_1\alpha_0},\overline{\alpha_1\alpha_0},
\overline{\alpha_1\alpha_0})+p^2\tl P\i(\overline{\alpha_0},\overline{\alpha_0},
\overline{\alpha_0})+p\tl P\i\tl P$$ is empty.
Therefore by Lemma 8.5(a), $g$ is not in the image of $\spi_0(X)$.
\endexample

\remark{Remarks} Sufficient conditions for the surjectivity of
$\spi_0(X)@>\sstau>>\pi_0(X)$, involving the finiteness condition FP$_n$ and
$n$-connectedness at infinity of the fundamental groups of nerves of $X$, are
found in \cite{GK; \S\S7,9}.
According to \cite{GK; Remark 10.3}, Ferry proved the following converse
of Theorem 8.4: every virtually Steenrod connected compactum (in particular,
Brin's compactum!) is Steenrod homotopy equivalent to a compactum $X$ such that
$\spi_0(X)\to\pi_0(X)$ is surjective.
\endremark

\proclaim{Theorem 8.7} {\rm (Krasinkiewicz--Minc \cite{KM})} Let $X$ be
a connected, Steenrod disconnected compactum.
Then the image of $\stau\:\spi_0(X)\to\pi_0(X)$ is uncountable (in particular,
nontrivial).
\endproclaim

We substantially simplify the original proof \cite{KM; pp.\ 147--151} (continued
from \cite{K2; pp.\ 146--147}), which depended on the theory of continua
(essentially using composants and indecomposable continua).
Note that Theorem 8.7 implies Theorems 8.2 and 8.3.

\demo{Proof} Let $\dots@>p_2>>P_2@>p_1>>P_1$ be an inverse sequence of polyhedra
with $\invlim P_i=X$ and such that $X$ is Steenrod disconnected over $P_1$
(see the definition in the remark to Theorem 8.3).
Let $G_i$ be the image of $\pi_1(P_i)$ in $\pi_1(P_1)$.
Let $BG_1$ be a classifying space for $\pi_1(P_1)$ (a locally compact polyhedron
of the homotopy type of $K(G_1,1)$), and let $BG_i$ be
the covering space of $BG_1$ corresponding to $G_i$.
Then the bonding maps $p^i_1\:P_i\to P_1$ lift to maps $f_i\:P_i\to BG_i$
commuting with the $p_i$ and with the covering maps $q_i\:BG_{i+1}\to BG_i$.
These converge to a map $f$ from $X$ into the inverse limit $BG$ of $BG_i$'s.
Let $Q_i$ be the union of all (closed) simplices of $BG_i$ that intersect
$f_i(P_i)$, where $BG_i$ is triangulated arbitrarily for $i=1$ and so that
$q_{i-1}$ is simplicial for $i>1$.
Note that each $Q_i$ is compact since $f_i(P_i)$ is.
Each inclusion induced map $\pi_1(Q_i)\to G_i$ is an epimorphism since $f_i$
factors through $Q_i$.
Hence $Y:=\invlim Q_i$ is Steenrod disconnected over $Q_1$.

By Lemma 8.5, each Steenrod component $\alpha$ of $BG$ is the image of
an inverse limit $W_\alpha$ of self-homeomorphisms of the universal cover $EG_1$
of $BG_1$.
If $EG_1$ is triangulated so that the projection $EG_1\to BG_1$ is simplicial,
we obtain that $\alpha$ is a continuous image of the triangulated polyhedron
$W_\alpha$, moreover the composition $W_\alpha\to\alpha\i BG\to BG_i$ is
simplicial.
In particular, it follows that $\alpha$ contains only countably many points of
the inverse limit $BG^{(0)}$ of the $0$-skeleta $BG_i^{(0)}$.
On the other hand, each point of the inverse limit $Y^{(0)}$ of the $0$-skeleta
$Q_i^{(0)}$ is contained in the same Steenrod component of $BG$ as some point
of $f(X)$.
Indeed, the ``closed star'' $S$ of each vertex $y\in Y^{(0)}$ (i.e.\ the image
of the star of a lift of $y$ into the triangulated polyhedron $W_\alpha$, where
$y\in\alpha$) is the inverse limit of the stars $S_i$ of $q^\infty_i(y)$.
Each $S_i\cap f_i(P_i)$ is nonempty, hence so is their inverse limit
$S\cap f(X)$.

Suppose that $Y^{(0)}$ is countable.
Then by the Baire Category Theorem 7.11 it contains an isolated point $v$.
Removing the ``open star'' of the vertex $v$ from $Y_0:=Y$, we obtain
a compactum that has finitely many connected components by the cohomological
Mayer--Vietoris
sequence (using that the ``link'' of $v$ is a compact polyhedron).
At least one of them is Steenrod disconnected over $Q_0$ by the remark
to Theorem 8.3.
This component $Y_1$ has the $0$-skeleton $Y_1^{(0)}:=Y_1\cap Y^{(0)}$, which
again contains an isolated point, so that the argument can be repeated.
Then $Y_\omega:=\bigcap_{i<\omega} Y_i$ is the inverse limit of the inclusions
$\dots\i Y_1\i Y_0$.
Hence $Y_\omega$ is connected (by the continuity of Pontryagin cohomology)
and homeomorphic to the limit of the `staircase' inverse sequence
$\dots@>p_{11}>>P_{11}@>f_{11}>>P_{01}@>p_{00}>>P_{00}$, where
each $Y_i=\invlim(\dots@>p_{i1}>>P_{i1}@>p_{i0}>>P_{i0}=P)$, and each
$f_{i,[0,\infty)}\:P_{i+1,[0,\infty)}\to P_{i,[0,\infty)}$ is the
level-preserving map representing the Steenrod homotopy class of the inclusion
$Y_{i+1}\i Y_i$.
Hence each $\Gamma_i:=\im[\pi_1(P_{ii})\to\pi_1(P)]$ properly contains
$\Gamma_{i+1}$, so $Y_\omega$ is Steenrod disconnected over $Q_0$.%
\footnote{Note that the intersection of a nested sequence of Steenrod disconnected
compacta does not have to be Steenrod disconnected.
For instance, the intersection of a nested sequence of clusters of solenoids
can be a single point.}

Finally, $Y_\omega^{(0)}:=Y_\omega\cap Y^{(0)}$ is again countable, and the process
can be continued by transfinite (countable) induction until it stops,
i.e.\ until we reach the stage $\lambda$ where $Y_\lambda$ is empty.
This is a contradiction, since the empty set $Y_\lambda$ cannot be Steenrod
disconnected over $Q_0$.

Thus $Y^{(0)}$ is uncountable.
Since each Steenrod component of $BG$ contains only countably many points of
$BG^{(0)}$, we obtain that $Y^{(0)}$ meets uncountably many Steenrod components
of $BG$.
Since each point of $Y^{(0)}$ is contained in the same Steenrod component of
$BG$ as some point of $f(X)$, the image of the composition
$\spi_1(X)\to\spi_1(BG)\to\pi_1(BG)$ is uncountable.
Since this composition factors through $\pi_1(X)$, the image
of $\stau\:\spi_1(X)\to\pi_1(X)$ is uncountable. \qed
\enddemo

\remark{Remark} The proof works if the map $f\:X\to BG$ is replaced by
a map of $Z:=\invlim\Gamma_i/G_1$ into $BG$, where each $\Gamma_i$ is
the Cayley graph of the action of a certain finite set of generators of $G_i$
on $G_1$ (which is disconnected unless $G_i=G_1$).
In more detail, $Z$ can be constructed as follows.
Let $D\i EG_1$ be a (closed) fundamental domain of the action of $G_1$.
Let $\tl X\to X$ be the pullback of $EG_1\to BG_1$.
Since $\tl X\to EG_1$ factors through the universal cover $\tl P_1$ of
$P_1$, and $\tl X\to\tl P_1$ is proper, the fundamental domain
$\hat D$ of the action of $G_1$ on $\tl X$ is compact.
Then $\check\pi_0(\hat D)$ is also compact, and so is its image $Z^{(0)}$ in
$\check\pi_0(\tl X):=\invlim\pi_0(\tl P_i)$.
As $G_1$ acts on $\check\pi_0(\tl X)$, each generator $g_i$ of $G_1$ sends
$Z^{(0)}$ to its copy $g_i Z^{(0)}$.
In particular, this yields a map $Z^{(0)}\cap g_i^{-1}Z^{(0)}\to Z^{(0)}$,
which may be viewed as a partial self-map $h_i$ of $Z^{(0)}$.
Let $Z$ be the union of the mapping tori of all $h_i$ corresponding to a fixed
collection of generators of $G_1$; since this collection may be assumed finite,
$Z$ is compact.
\endremark

\remark{Remark} Theorem 8.7 has a homological analogue: If $X$ is connected
but $H_0(X)\ne 0$, then the image of $\stau\:\sH_0(X)\to H_0(X)$ is uncountable.
Compared with Theorem 8.7, this assertion is much easier; it can be verified
by comparing $X$ with an $l$-adic solenoid ($l$ is a sequence of primes).
We leave the details to the reader.
\endremark

\medskip

We are now ready for the promised $0$-dimensional counterpart of
Corollary 6.6:

\proclaim{Corollary 8.8} A compactum $X$ is semi-LC$_0$ iff
$\pi_0(X)$ is discrete.
\endproclaim

\demo{Proof} The ``if' part follows from Lemma 3.4(b).
Conversely, let us fix some metric on $X$, and let $\eps>0$ be such that every
two $\eps$-close points of $X$ represent the same Steenrod homotopy class.
Obviously, $\check\pi_0(X)$ is discrete.
If $C$ is a component of $X$, any pair of points $p,q\in C$ is connected by
a chain $p=p_0,p_1,\dots,p_n=q$, where each $p_{i+1}$ is $\eps$-close to $p_i$.
Hence all points of $C$ represent the same Steenrod homotopy class.
By Theorem 8.7, $C$ is Steenrod connected. \qed
\enddemo

\remark{Remark} By regarding $\spi_n(X)$ as $\spi_0(\Omega^nX)$, where
$\Omega^nX$ is the iterated loop space of $X$ (which is a separable metrizable
complete uniform space), the proof of Theorem 8.7 can apparently be used to
prove that if $X$ is an LC$_{n-1}$ and UV$_1$ compactum and
$\ker[\pi_n(X)@>\check\tau>>\check\pi_n(X)]$ is nontrivial, then the image of
$\ker[\spi_n(X)@>\check\tau\stau>>\check\pi_n(X)]$ in $\pi_n(X)$ is
uncountable.
(When $n\ge 2$, this assertion also follows from Theorems 6.5 and 6.1
and the uncountability of every nontrivial derived limit of countable groups,
which in turn follows from Lemma 3.3.)
The details will hopefully appear elsewhere (see also Lemma 7.9(g)).
\endremark

\proclaim{Theorem 8.9} A connected compactum $X$ is Steenrod connected iff
every overlay space over $X$ has countably many uniform components.
\endproclaim

The proof continues the analysis begun in the proof of Theorem 8.7.

\demo{Proof} Let us fix a map $f\:X\to P$, where $P$ is a polyhedron.
We will show that $X$ is Steenrod connected over $P$ (see the definition in
the remark after Theorem 8.3) iff the pullback $\tl X$ of the universal cover
$\tl P$ of $P$ has countably many uniform components.

We may assume that $X=\invlim(\dots\to P_2\to P_1=P)$, where the $P_i$ are
compact polyhedra and the bonding maps are PL.
Let $\tl P_i$ be the pullback of $\tl P$ over $P_i$.
Then $G:=\pi_1(P)$ acts on the components of $\tl P_i$ with stabilizer
$\im(G_i\to G)$, where $G_k=\pi_1(P_i)$.
If $X$ is Steenrod connected over $P$, let $k$ be such that
$\im(G_k\to G)=\im(G_i\to G)$ for each $i>k$.
Then the $G$-equivariant map $\pi_0(\tl P_i)\to\pi_0(\tl P_k)$ is
a bijection for each $i>k$.
Hence $\check\pi_0(\tl X)$ is countable.
On the other hand, since $\dots\to P_2\to P_1$ is convergent (see \S7), so is
$\dots\to\tl P_2\to\tl P_1$.
Now $\tl X$ is the uniformly disjoint union of the preimages $\tl X^c$ of
the components $\tl P_k^c$ of $\tl P_k$.
Hence if $\tl P_i^c$ denotes the corresponding component of $\tl P_i$ for $i>k$,
the inverse sequence $\dots\to\tl P_{k+1}^c\to\tl P_k^c$ is convergent.
Then by Lemma 7.9(f), its inverse limit $\tl X^c$ is uniformly connected.
Thus the set of uniform components of $\tl X$ injects into $\check\pi_0(\tl X)$.

Conversely, suppose that $X$ is Steenrod disconnected over $P$; we shall use
the notation in the proof of Theorem 8.7.
The universal covering $\tl P\to P$ is induced from $EG_1\to BG_1$ by $f_1$.
Let $\tl Q_i$ and $\tl Y$ be the pullbacks of $EG_1$ over $Q_i$ and $Y$.
Since each $Q_i$ is embedded into the covering space $(q^i_1)^{-1}(Q_1)\i BG_i$
of $Q_1$, the map $\tl Q_i\to\tl Q_1$ embeds each component of $\tl Q_i$.
On the other hand, since $Y^{(0)}$ has been proved to be uncountable, there
exists an $x\in Q_1^{(0)}$ whose preimage $(q^\infty_1)^{-1}(x)$ in $Y^{(0)}$
is uncountable.
Pick a lift $\tl x\in\tl Q_1$; then its preimage in $\tl Y$ is uncountable,
and at the same time the points of its preimage in each $\tl Q_i$ lie in
distinct components of $\tl Q_i$.
Thus $\tl Y$ has uncountably many uniform components.
However, by the proof of Theorem 8.7, $f\:X\to Y$ induces a surjection on
path components.
Therefore so does its lift $\tl f\:\tl X\to\tl Y$. \qed
\enddemo

\remark{Remark} Theorem 8.9 implies Theorem 8.7, because if $f\:E\to B$ is
an overlaying, every point-inverse of $f_*\:\pi_0(E)\to\pi_0(B)$ is countable.
Indeed, the point-inverse of the Steenrod homotopy class of the basepoint is
countable similarly to Theorem 3.15(b); and by the proof of 3.15(b) we are
free to choose any Steenrod homotopy class $b\:pt\ssm B$ that factors into
composition of some $\tl b\:pt\to E$ and the overlaying $f$ as our
``basepoint''.
\endremark

\medskip
As remarked in \S3, an inverse limit of fibrations $E_i\to B_i$ need not be
a fibration; however, an inverse limit of compositions $f_i$ of fibrations with
the pullbacks of the preceding compositions $f_{i-1}$ is a fibration.
In particular, an inverse limit of coverings is a fibration.
(This does not depend on whether the inverse sequence of the total spaces is
convergent.)
On the other hand, by definition (see \S3), the limit of an inverse sequence
of fibrations over compact polyhedra is a Steenrod fibration provided that
the inverse sequence of the total spaces is convergent (see \S7).
The following example shows that the latter assumption cannot be dropped.

\example{Example 8.10} Let $X\subset \Sigma_p\x\Sigma_p\x I$ be defined
as the union $[\Sigma_p\x\Sigma_p\x\partial I]\cup[\Sigma_p\x pt\x I]\cup
[pt\x\Sigma_p\x I]$.
Thus $X=\invlim(\dots@>p_1>>P_1@>f>>P_0)$, where each
$P_i=[S^1\x S^1\x\partial I]\cup[S^1\x pt\x I]\cup[pt\x S^1\x I]$ and each $f_i$
is the restriction of the self-map $p\x p\x\id_I$ of $S^1\x S^1\x I$.
Let us consider the Steenrod component $h=(\alpha+\Z,\,\alpha+\Z)$ in
$(\Z_p/\Z)^2\simeq\pi_0(X)$, where
$\alpha=\overline{\dots\alpha_1\alpha_0}\in\Z_p\but\Z$ (compare Example 8.6).
Let $X^h\i X$ be the set of its points, i.e.\ the union of all path components
of $X$ that map into $h$.
Let $\phi\:\tl X\to X$ be the overlaying, induced from the universal covering
$\tl P_1\to P_1$, and let $\phi^h\:\tl X^h\to X^h$ denote its restriction to
$\tl X^h:=\phi^{-1}(X^h)$.

We claim that $\phi^h$ is not a Steenrod fibration.
Indeed, let us consider the pullbacks $\tl P_i\to P_i$ of the universal
covering of $P_1$.
We can identify $\tl X^h$ with the inverse limit $\R\x\R\x\partial I$ of
the subpolyhedra $\tl P_i^h:=[\R\x\R\x\partial I]\cup[K_i\x I]$ of the $\tl P_i$,
where $K_i=[(\overline{\alpha_{i-1}\dots\alpha_0}+p^i\Z)\x\R]\cup
[\R\x(\overline{\alpha_{i-1}\dots\alpha_0}+p^i\Z)]$.
Note that $\tl b:=(0,0,0)$ and $\tl b':=(0,0,1)$ belong to different components
of $\tl X^h$ but the corresponding proper rays $\tl b_{[0,\infty)}$ and
$\tl b_{[0,\infty)}$ are properly homotopic in the mapping telescope
$\tl P_{[0,\infty)}^h$.
This proper homotopy projects to a Steenrod path $\ell\:(I,\partial I)\x[0,\infty)
\to(P_{[0,\infty)},b_{[0,\infty)}\cup b'_{[0,\infty)})$, where
$b=\phi^h(\tl b)$ and $b'=\phi^h(\tl b')$.

If $\phi^h$ is a Steenrod fibration, there exists a convergent (see \S7) inverse
sequence $\dots\to\hat P_2^h\to\hat P_1^h$ with $\invlim\hat P_i^h=\tl X^h$
and a level-preserving map $w_{[0,\infty)}\:\hat P_{[0,\infty)}^h\to P_{[0,\infty)}$
with $\invlim w_i=\phi^h$ such that, by the proof of Theorem 3.15(b),
$p^{[0,\infty)}_{[-3,\infty-3)}\ell$ lifts to a Steenrod path
$\hat\ell\:I\x[0,\infty)\to\hat P_{[0,\infty)}^h$ starting at $\hat b_{[0,\infty)}$,
where $\hat b_i$ denotes the image of $\tl b$ in $\hat P_i^h$.
This $\tl\ell$ has to end at $\hat b''_{[0,\infty)}$, where $\hat b_i''$ denotes
the image in $\hat P_i^h$ of some lift $\tl b''\in\tl X^h$ of $b'$.
Then $\tl b''$ and $\tl b$ belong to the same Steenrod component of
$\tl X^h=\R\x\R\x\partial I$, which is a contradiction since $\tl b=(0,0,0)$
and $\tl b'=(m,n,1)$ for some $m,n\in\Z$.
\endexample

\remark{Remark} The above example also shows that the inclusion induced maps
$\pi_0(\tl X^h)\to\pi_0(\tl X)$ and $\check\pi_0(\tl X^h)\to\check\pi_0(\tl X)$
need not be injective.
This may be of interest in connection with the proof of the ``easy'' (only if)
direction in Theorem 8.9.
\endremark

\end